\documentclass[a4paper, 11pt]{article}

\usepackage{inputenc}
\usepackage{amsmath, amsthm, amssymb}
\usepackage{color}
\usepackage{graphicx}
\usepackage{hyperref}
\hypersetup{colorlinks=true,linkcolor=blue}
\usepackage{a4wide}
\usepackage{tikz}
\usepackage[font=small,labelfont=bf, margin=2ex]{caption}
\usepackage{mathtools}

\numberwithin{equation}{section}

\newtheorem{thm}{Theorem}[section]
\newtheorem{lem}[thm]{Lemma}
\newtheorem{prop}[thm]{Proposition}
\newtheorem{cor}[thm]{Corollary}

\theoremstyle{definition}

\newtheorem{defi}[thm]{Definition}

\newtheorem{rem}[thm]{Remark}
\newtheorem{prob}[thm]{Problem}

\newenvironment{prf}{\noindent\textbf{Proof.}}
{\hfill $\Box$\\[-1ex]}

\newenvironment{abst}{\begin{minipage}[c]{0.9\textwidth} \footnotesize \textbf{Abstract.}}
{\end{minipage}\\[2ex]}

\newenvironment{key}{\begin{minipage}[c]{0.9\textwidth} \footnotesize \textbf{Keywords.}}
{\end{minipage}\\[2ex]}

\newenvironment{amsclass}{\begin{minipage}[c]{0.9\textwidth} \footnotesize \textbf{AMS Subject Classification.}}
{\end{minipage}\\[5ex]}

\newenvironment{ackno}{\begin{minipage}[c]{1\textwidth} \footnotesize \textbf{Acknowledgements.}}
{\end{minipage}\\[5ex]}

\newcommand{\eps}{\varepsilon}
\newcommand{\m}{\mathbf{m}}
\newcommand{\B}{\mathbf{B}}
\newcommand{\M}{\mathbf{M}}
\newcommand{\E}{\mathbf{E}}
\newcommand{\R}{\mathbb{R}}
\newcommand{\BB}{\mathbb{B}}
\renewcommand{\SS}{\mathbb{S}}
\newcommand{\supp}{\textnormal{supp}}
\renewcommand{\H}{\mathcal{H}}
\newcommand{\D}{\mathcal{D}}
\newcommand{\DD}{\textrm{D}}
\newcommand{\N}{\mathcal{N}}
\newcommand{\G}{\mathcal{G}}
\newcommand{\T}{\mathcal{T}}

\newcommand{\Hp}{\mathcal{H}_{+,R}^2}
\newcommand{\Hm}{\mathcal{H}_{-,R}^2}
\newcommand{\Hpo}{\mathcal{H}_{+,R_1}^2}
\newcommand{\Hmo}{\mathcal{H}_{-,R_1}^2}
\newcommand{\Hpz}{\mathcal{H}_{+,R_0}^2}
\newcommand{\Hmz}{\mathcal{H}_{-,R_0}^2}

\newcommand{\id}{\textnormal{Id}}

\newcommand{\RR}{\mathbb R}

\newcommand{\sph}{\mathbb S}

\newcommand\atopp[2]{\genfrac{}{}{0pt}{}{#1}{#2}}
\newcommand{\dd}{{\mathrm{d}}}

\newcommand{\LL}{\textnormal{L}}
\DeclareMathOperator*{\esssup}{\textnormal{ess.}\,\sup}

\begin{document}
\begin{center}
\Large\bfseries On the Recovery of Core and Crustal Components of Geomagnetic Potential Fields\normalsize\mdseries
\\[3ex]
{L. Baratchart\footnotemark[1], C. Gerhards\footnotemark[2]}
\footnotetext[1]{INRIA, Project APICS, 2004 route de Lucioles, BP 93, Sophia-Antipolis F-06902 Cedex, France, 
e-mail: laurent.baratchart@sophia.inria.fr}
\footnotetext[2]{University of Vienna, Computational Science Center, Oskar-Morgenstern-Platz 1, A-1090 Vienna, Austria, e-mail: christian.gerhards@univie.ac.at}
\\[3ex]
\today
\end{center}

\begin{abst}
In Geomagnetism it is of interest to separate the Earth's core magnetic field from the crustal magnetic field. However, measurements by satellites can only sense the sum of the two contributions. In practice, the measured magnetic field is expanded in spherical harmonics and separation into crust and core contribution is achieved empirically, by a sharp cutoff in the spectral domain. In this paper, we derive a mathematical setup in which
the two contributions are modeled by harmonic potentials $\Phi_0$ and $\Phi_1$ generated on two different spheres $\SS_{R_0}$ (crust) and $\SS_{R_1}$ (core) with radii $R_1<R_0$.
Although it is not possible in general to recover 
$\Phi_0$ and $\Phi_1$ knowing their 
superposition $\Phi_0+\Phi_1$ on a sphere $\SS_{R_2}$ with radius $R_2>R_0$, 
we show that it becomes possible if the magnetization $\mathbf{m}$ generating $\Phi_0$ is localized in a strict subregion of $\SS_{R_0}$. 
Beyond unique recoverability, we show in this case 
how to numerically reconstruct 
characteristic features of $\Phi_0$ (e.g., 
spherical harmonic Fourier coefficients). 
An alternative way of phrasing the results is
that knowledge of  $\mathbf{m}$ on a nonempty open subset
of  $\SS_{R_0}$ allows one to perform separation.
\end{abst}

\begin{key}
Harmonic Potentials, Hardy-Hodge Decomposition, Separation of Sources, Geomagnetic Field, Extremal Problems
\end{key}

\begin{amsclass}
33C55, 42B37, 45Q05, 53A45, 86A22
\end{amsclass}

\section{Introduction}

The Earth's magnetic field $\B$, as measured by several satellite missions, 
is a superposition of various contributions, e.g., of iono-/magnetospheric fields, crustal magnetic field, and of the core/main magnetic field, see
\cite{kono09,hulot10,olsen15} for an overview 
and \cite{lesur10,maus08,sabaka15,thebault15} for some recent geomagnetic field models. While iono-/magnetospheric contributions can to a certain extent be filtered out due to their temporal variations, the separation of the core/main field $\B_{core}$ and the crustal field $\B_{crust}$ is typically based on the empirical observation that the power spectra of Earth magnetic field models have a sharp knee at spherical harmonic degree 15 (see, e.g., \cite{langel82,olsen15}). However, under this spectral separation, large-scale contributions (i.e., spherical harmonic degrees smaller than 15) are entirely neglected in crustal magnetic field models. In \cite{holschneider16}, a Bayesian approach has been proposed that addresses the separation of geomagnetic sources based on their correlation structure. The correlation of certain components, e.g., internally and externally produced magnetic fields, can (to some extent) be obtained from the underlying geophysical equations. But this approach does not address the problem that some of the involved separation problems, e.g., the separation into crustal and core magnetic field contributions, are generally not unique for the given data situation. The goal of this paper is to derive conditions under which a rigorous separation of the contributions $\B_{crust}$ and $\B_{core}$ is possible, as well as to formulate extremal problems whose solutions lead to approximations of these contributions or certain features thereof. The main 
assumption that we make for our approach to work 
is that the magnetization generating $\mathbf{B}_{crust}$ is localized 
in a strict subregion of the crust.  By linearity, this is equivalent to 
assuming that this magnetization is known on a spherical cap 
that may, in principle, be
arbitrary small.  For applications, this is interesting in as much as that the crustal magnetization may be estimated
in certain places of the Earth from local measurements. Thus, given such a local estimation, its contribution can be substracted from
global magnetic field measurements to yield a crustal contribution that
stems from magnetizations localized in a strict subregion of the Earth (namely the complement of those places where a local estimate of the magnetization has been 
performed), thereby allowing us to apply the separation approach indicated in this paper. Similarly, if one can identify places on the Earth which are 
only weakly magnetized as compared to others,
the separation process that we will describe
may reasonably be applied by neglecting magnetizations 
in such places.  

We assume throughout
that the overall magnetic field is of the form $\B=\B_{crust}+\B_{core}$ in $\R^3\setminus\overline{\BB_{R_0}}$, where $\BB_{R_0}=\{x\in\R^3:|x|<R_0\}$ denotes the ball of radius $R_0>0$ and overline indicates
closure (here $R_0$ can be interpreted as the radius of the Earth). Since the sources of $\B_{crust}$ and $\B_{core}$ are located inside $\BB_{R_0}$ (hence, the corresponding magnetic fields are curl-free and divergence-free in $\R^3\setminus\overline{\BB_{R_0}}$), there exist potential fields $\Phi$, $\Phi_{crust}$, $\Phi_{core}$ such that $\B=\nabla \Phi$, $\B_{crust}=\nabla \Phi_{crust}$, and $\B_{core}=\nabla \Phi_{core}$ in $\R^3\setminus\overline{\BB_{R_0}}$. Therefore, from a mathematical point of view, the problem reduces to finding unique $\Phi_{crust}$, $\Phi_{core}$ from the knowledge of $\Phi$ (but we should keep in mind that the actual measurements bear on the magnetic field $\B$).

It is known that $\B_{crust}$  is generated by a magnetization $\M$ confined in
a thin spherical shell $\BB_{R_0-d,R_0}=\{x\in\R^3:R_0-d<|x|<R_0\}$ of thickness $d>0$ (for the Earth, $d\approx 30$km is typical), therefore
the corresponding magnetic potential can be expressed as (see, e.g., \cite{blakely95, gubbins11})
\begin{align}
 \Phi_{crust}(x)=\frac{1}{4\pi}\int_{\BB_{R_0-d,R_0}}\M(y)\cdot\frac{x-y}{|x-y|^3}\,\dd\lambda(y),\quad x\in\R^3,\label{eqn:phivol}
\end{align}
where the dot indicates Euclidean scalar 
product in $\RR^3$ and $\lambda$ the Lebesgue measure. Due to the thinness of the magnetized layer relative to the Earth's radius, it is reasonable to substitute the volumetric $\M$ by a spherical magnetization $\m$ (i.e., $\mathbf{M}=\mathbf{m}\otimes\delta_{\SS_{R_0}}$ in a distributional sense). Then, the magnetic potential \eqref{eqn:phivol} becomes 
\begin{align}
 \Phi_{crust}(x)=\frac{1}{4\pi}\int_{\SS_{R_0}}\m(y)\cdot\frac{x-y}{|x-y|^3}\,\dd \omega_{R_0}(y),\quad x\in\R^3\setminus\SS_{R_0},\label{eqn:phisurf}
\end{align}
where  $\SS_{R_0}=\{x\in\R^3:|x|=R_0\}$ denotes
the sphere of radius $R_0>0$ and $\dd \omega_{R_0}$ the corresponding surface 
element. When interested in reconstructing the actual magnetization $\M$, 
substituting a spherical magnetization $\mathbf{m}$ is of course a significant restriction (however, one that is fairly frequent in Geomagnetism). But since our main focus is on $\B_{crust}$ and the corresponding potential $\Phi_{crust}$ rather than the magnetization itself, this restriction actually involves 
no loss of information: in Section \ref{sec:harmpot} we 
show that, under mild summability assumptions,
 any potential $\Phi_{crust}$ produced by a volumetric magnetization $\M$ in $\BB_{R_0-d,R_0}$ can also be generated by a spherical magnetization $\m$ on $\SS_{R_0}$. 

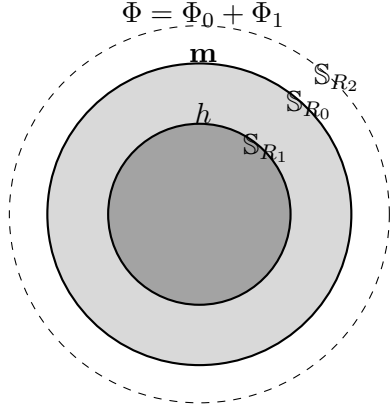
\begin{figure}
\begin{center}
\begin{tikzpicture}
\fill[gray,opacity=0.3] (0,0) circle (2);
\fill[gray,opacity=0.6] (0,0) circle (1.2);
\draw[thick] (0,0) circle (2);
\draw[thick] (0,0) circle (1.2);
\node at (0.05,2.1) {{$\m$}};
\node at (1.45,1.45) {$\SS_{R_0}$};
\node at (0.88,0.88) {$\SS_{R_1}$};
\node at (0.05,1.35) {$h$};
\draw[dashed] (0,0) circle (2.5);
\node at (0.05,2.65) {{$\Phi=\Phi_0+\Phi_1$}};
\node at (1.81,1.81) {$\SS_{R_2}$};
\end{tikzpicture}
\end{center}
\caption{Illustration of the setup of Problem \ref{prob:1}.}\label{fig:setup}
\end{figure}

The core/main contribution $\B_{core}$ is governed by the Maxwell equations (see, e.g., \cite{backus96})
\begin{align}
\label{Maxwell}
 \nabla\times \B_{core}&=\sigma(\E+\mathbf{u}\times \B_{core}),
 \\\nabla\cdot \B_{core}&=0,
 \\\nabla\times\E&=-\partial_t\B_{core},
 \\\nabla\cdot\E&=\rho,
\end{align}
where $\sigma$ denotes the conductivity, $\rho$ the charge density, and $\mathbf{u}$ the fluid velocity in the Earth's outer core (the constant permeability $\mu_0$ and permittivity $\eps_0$ have been set to $1$). The conductivity $\sigma$ is assumed to be zero outside a sphere $\SS_{R_1}$ of radius $0<R_1<R_0$. The condition $R_1<R_0$ is crucial to the forthcoming  arguments and is justified by common geophysical practice and results (see, e.g., \cite{ballani02,puethe15}). In particular it implies that $\nabla\times \B_{core}=0$ in 
$\R^3\setminus\overline{\BB_{R_1}}$, therefore, $\B_{core}=\nabla \Phi_{core}$ in $\R^3\setminus\overline{\BB_{R_1}}$ for some harmonic potential $\Phi_{core}$. Although the geophysical processes in the Earth's outer core can be extremely complex, of importance to us is only  that $\Phi_{core}$ can be expressed in $\R^3\setminus\overline{\BB_{R_1}}$ as a Poisson transform:
\begin{align}\label{eqn:phicore}
 \Phi_{core}(x)=\frac{1}{4\pi R_1}\int_{\SS_{R_1}} h(y)\frac{|x|^2-R_1^2}{|x-y|^3}\dd \omega_{R_1}(y),\quad x\in\R^3\setminus\overline{\BB_{R_1}},
\end{align}
for some scalar valued auxiliary function $h$ on $\SS_{R_1}$; this follows 
from previous considerations which imply that $ \Phi_{core}$ is harmonic in
$\R^3\setminus\overline{\BB_{R_1}}$ and continuous in
$\R^3\setminus\BB_{R_1}$. Summarizing, the problem we treat in this paper is the following (the setup is illustrated in Figure \ref{fig:setup}):

\begin{prob}\label{prob:1}
Let $\Phi\in L^2(\SS_{R_2})$ be given on a sphere $\SS_{R_2}\subset\R^3\setminus\overline{\BB_{R_0}}$ of radius $R_2>R_0$. Assume $\Phi$ is decomposable into $\Phi=\Phi_0+\Phi_1$ on $\SS_{R_2}$, where $\Phi_0=\Phi_0[\m]$ is of the form \eqref{eqn:phisurf}, with $\m\in L^2(\SS_{R_0},\R^3)$, and  $\Phi_1=\Phi_1[h]$ is of the form \eqref{eqn:phicore}, with $h\in L^2(\SS_{R_1})$ and $R_1<R_0$. Are $\Phi_0$ and $\Phi_1$ uniquely determined by the knowledge of $\Phi$ on $\SS_{R_2}$, and if yes can they be reconstructed efficiently?
\end{prob}

The answer to the uniqueness issue in Problem \ref{prob:1}
is generally negative. But under the additional assumption that $\supp(\m)\subset\Gamma_{R_0}$ for a strict subregion $\Gamma_{R_0}\subset\SS_{R_0}$
( i.e.  $\overline{\Gamma_{R_0}}\not=\SS_{R_0}$), uniqueness is guaranteed. 
This follows from results in \cite{baratchart13,lima13} and their formulation on the sphere in \cite{gerhards16a}, to be reviewed in greater 
detail in 
Section \ref{sec:seppot}. In fact, we show in this case
that $h$ and the curl-free contribution of $\m$ can be reconstructed uniquely from the knowledge of $\Phi$. Additionally, we provide a means of approximating $\langle\Phi_0,g\rangle_{L^2(\SS_{R_2})}$  knowing $\Phi$ on $\SS_{R_2}$, where $g$ is some appropriate test function (e.g., a spherical harmonic). This 
allows one to separate the crustal and the core contributions to
the Geomagnetic potential if, e.g., the crustal magnetization 
can be estimated over a small subregion on Earth by other means.

Throughout the paper, we call $\Phi_0$ the crustal contribution and $\Phi_1$ the core contribution. We should point out that the examples we provide 
at the end of the paper are not based on real Geomagnetic field data but they reflect some of the main properties of realistic scenarios (e.g., the domination of the core contribution at low spherical harmonic degrees). In Section \ref{sec:harmpot}, we take a closer look at harmonic potentials of the form \eqref{eqn:phivol} and \eqref{eqn:phisurf} and show that the balayage  onto $\SS_{R_0}$ of a volumetric potential supported in $\BB_{R_0-d,R_0}$ preserves divergence form. More precisely, if
$\mathbf{M}$ is supported in $\BB_{R_0-d,R_0}$ and its restriction to
$\SS_R$ is uniformly square-summable for
$R\in (R_0-d,R_0)$, then
there exists a spherical magnetization $\mathbf{m}$ supported on  $\SS_{R_0}$,
which is square summable and generates the same potential as 
$\mathbf{M}$ in $\R^3\setminus\overline{\BB_{R_0}}$. 
The latter property justifies the above-described modeling of the crustal magnetic field. Basic background and auxiliary material on  geometry, 
spherical decomposition of vector fields as well as Sobolev and
Hardy spaces is recapitulated in Section \ref{sec:aux}. Some parts in the beginning are described in more detail than necessary for the core part of this paper and are only required again in Appendices \ref{sec:appendix2} and \ref{sec:appendix1}. So the reader familiar with the background and notation may directly proceed to Definition \ref{def:HH}. Eventually, in Section \ref{sec:num} we provide some initial  examples of numerical 
approximation of $\Phi_0$ and $\langle\Phi_0,g\rangle_{L^2(\SS_{R_2})}$, followed
by a brief conclusion in Section \ref{sec:conc}. Some technical results 
on potentials of distributions, gradients, and divergence-free vector fields are gathered in the appendices.

\section{Auxiliary Notations and Results}\label{sec:aux}

We start with some basic definitions of function spaces and differentiation on the sphere. For $R>0$, the sphere $\SS_R$ is a smooth, compact oriented surface embedded in $\R^3$. That is, $\SS_R$ can be described by finitely many charts $\psi_j:U_j\to V_j$ (for open subsets $U_j\subset\SS_R$ and $V_j\subset \R^2$, $j=1,\ldots,N$), which allows a meaningful definition of the surface area measure $\omega_R$ on the sphere $\SS_R$ via the Lebesgue measure $\lambda$ in $\R^2$. For $x\in U_j\subset\SS_R$, the tangent space $T_x$ at $x$ is the image of the derivative
$\DD\psi_j^{-1}[\psi_j(x)]:\R^2\to\R^3$. The tangent space
may  be described intrinsically as 
$T_x=\{y\in\R^3:x\cdot y=0\}$.
A $k$-times differentiable 
or $C^{k}$-smooth function  $f:\SS_R\to\R$ is a function
such that 
$f\circ\psi_j^{-1}$ is $k$-times differentiable or has continuous partial derivatives up
to order $k$, respectively, for each $j=1,\ldots,N$. We simply say that $f$ is smooth if it is $C^{\infty}$-smooth. Due to the simple geometry of the sphere $\SS_R$,
this definition of differentiability is in fact equivalent to requiring that the radial extension 
$\bar{f}(x)=f(R\frac{x}{|x|})$ of  $f$ has the corresponding  regularity in 
$\R^3\setminus\{0\}$. This allows us to express the surface gradient $\nabla_{\SS_R}f(x)$ of a differentiable function $f:\SS_R\to\R$ at a point $x\in\SS_R$ via the relation $\nabla_{\SS_R}f(x)=\nabla\bar{f}(y)|_{y=x}$, where $\nabla$ denotes the Euclidean gradient. Formally, the surface gradient at $x$ is defined as the unique vector $v\in T_x$ such that the differential $\dd f[x]:T_x\to \R$ can be identified by the scalar product with $v$, i.e., $\dd f[x](y)=v\cdot y$ for $y\in T_x$. The differential of $f$ at $x\in U_j$ is the linear map
$\dd f(x):T_x\to\RR$ given at $v\in T_x$
by $\dd f[x](v)=\dd (f\circ \psi_j^{-1})[y](w)$,
where $y=\psi_j(x)$ and
$w\in\RR^2$ is such that $v=\dd \psi_j^{-1}[y](w)$. Here, the Euclidean 
differential  $\dd(f\circ \phi_j)[y]$ is defined as usual: $\dd (f\circ \psi_j^{-1})[y](w)=\partial_{y_1}(f\circ\psi_j^{-1})[y]w_1+
\partial_{y_1}(f\circ\psi_j^{-1})[y]w_2$, where $\partial_{y_i}$ indicates partial derivative with respect to $y_i$.

Furthermore, $L^2(\SS_{R})$ is denoted to be the space of square-integrable 
scalar valued functions $f:\SS_{R}\to\R$, while $L^2(\SS_{R},\R^3)$ denotes the space of square integrable vector valued spherical functions $\mathbf{f}:\SS_{R}\to\R^3$,  equipped with the inner products $\langle f,h\rangle_{L^2(\SS_{R})}=\int_{\SS_{R}}f(y)h(y)\dd \omega_{R}(y)$ and  $\langle \mathbf{f},\mathbf{h}\rangle_{L^2(\SS_{R},\R^3)}=\int_{\SS_{R}}\mathbf{f}(y)\cdot\mathbf{h}(y)\dd \omega_R(y)$, respectively.
A vector field $\mathbf{f}:\SS_R\to\R^3$ 
is said to be tangential  if $\mathbf{f}(x)\in T_x$ for all $x\in\SS_R$. The subspace of all tangential vector fields in $L^2(\SS_R,\R^3)$ is denoted by
$\T_R$. Note that the smooth vector fields are dense in $\T_R$. Clearly, if $f$ is smooth, then $\nabla_{\SS_R}f$ lies in $\T_R$.
The Sobolev space $W^{1,2}(\SS_R)$ may be defined as 
the completion of smooth functions 
with respect to the norm \cite{Hebey}
\begin{align*}
\label{defSobinf}
\|f\|_{W^{1,2}(\SS_R)}=\Bigl(\|f\|^2_{L^2(\SS_R)}+
\|\nabla_{\SS_R} f\|^2_{L^2(\SS_R,\RR^3)}\Bigr)^{1/2}.
\end{align*}
Since, for an appropriate set of charts $\psi_j:U_j\to V_j$, $j=1,\ldots,N$, of the sphere, the $V_j$ are bounded and the corresponding determinants of the metric tensors are bounded from above and below by 
strictly positive constants, it holds that $f\in W^{1,2}(\SS_R)$ if 
and only if the functions
$f\circ\psi_j^{-1}$ lie in the Euclidean Sobolev spaces $W^{1,2}(V_j)$ (see, e.g., \cite{lions68}). The gradient $\nabla_{\SS_R} f(x)$ at $x\in\SS_R$ of a function $f\in W^{1,2}(\SS_R)$ still satisfies the
representation $\dd f[x](y)=\nabla_{\SS_R} f(x)\cdot y$ for $y\in T_x$, where $\dd f$ has to be understood in the sense of distributional derivatives and $\nabla_{\SS_R} f(x)$
needs not be a pointwise derivative in the strong sense (see \cite[Ch.VIII]{stein70}). Let us put 
\begin{align*}
\G_R=\{\nabla_{\SS_R} f:f\in W^{1,2}(\SS_R)\}.
\end{align*}
We claim that $\G_R$ is closed in $L^2(\SS_R,\R^3)$. Indeed, if $\nabla_{\SS_R} f_n$
is a Cauchy sequence in $\G_R$, where $f_n\in W^{1,2}(\SS_R)$ 
is defined up to an additive constant, we may pick $f_n$ so that
$\int_{\SS_{R}} f_n\dd \omega_R=0$ and then it follows from the H\"older and the 
Poincar\'e 
inequalities \cite[Prop. 3.9]{Hebey} that 
$\|f_n-f_m\|_{L^2(\SS_R)}\leq C\|\nabla_{\SS_R}f_n-\nabla_{\SS_R}f_m\|_{L^2(\SS_R,\R^3)}$ for some constant $C$. Hence $f_n$ is a Cauchy sequence in 
$W^{1,2}(\SS_R)$,
therefore it converges to some $f$ there 
and consequently $\nabla_{\SS_R}f_n$ converges to
$\nabla_{\SS_R} f$ in $L^2(\SS_R,\R^3)$. Thus,
$\G_R$ is complete and therefore it is closed in $L^2(\SS_R,\R^3)$,
{which proves the claim}.

When $\mathbf{h}$ is a smooth tangential vector field on $\SS_R$, 
its surface divergence $\nabla_{\SS_R}\cdot\mathbf{h}$ is the smooth real valued function such that
\begin{equation}
\label{divsmooth}
\int_{\SS_R}f\,\nabla_{\SS_R}\cdot\mathbf{h}\,\dd \omega_R=-
\int_{\SS_R}(\nabla_{\SS_R} f)\,\cdot\,\mathbf{h}\,\dd \omega_R,\quad\textnormal{ for all }
f \in C^\infty(\SS_R).
\end{equation}
When $\mathbf{h}\in \T_R$ is not smooth, \eqref{divsmooth} must be interpreted in a weak sense,
namely $\nabla_{\SS_R}\cdot\mathbf{h}$
is the distribution on $\SS_R$ acting on smooth real-valued functions
by $\langle f\,,\nabla_{\SS_R}\cdot\mathbf{h}\rangle=-
\int_{\SS_R}\nabla_{\SS_R} f\,\cdot\, \mathbf{h}\,\dd \omega_R$, for all $f \in C^\infty(\SS_R)$. This clearly extends by density
to a linear form on $W^{1,2}(\SS_R)$, upon letting $f$ converge to a Sobolev 
function. Then it is apparent  that 
\begin{align*}
\D_R=\{\mathbf{h}\in\T_R:\nabla_{\SS_R}\cdot\mathbf{h}=0\} 
\end{align*}
is the orthogonal complement to $\G_R$ in $\T_R$. In particular,
\begin{align}
 \T_R=\G_R\oplus \D_R,\label{eqn:hd}
\end{align}
which is the so-called Helmholtz-Hodge decomposition. 
The particular geometry of $\SS_R$ makes it  easy to see that $\mathbf{f}\in\D_R$ 
if and only if its radial extension
$\bar{\mathbf{f}}(x)=\mathbf{f}(R\frac{x}{|x|})$ is  divergence free, as a $\R^3$-valued
distribution on $\R^3\setminus\{0\}$.

We now consider the operator
$J_x:T_x\to T_x$ given by $J_x(y)=\frac{x}{|x|}\times y$, for $y\in T_x$, where $\times$ indicates 
the vector product in $\R^3$; that is, $J_x$ is the rotation by $\pi/2$ in
$T_x$. We define
$J:\T_R\to\T_R$ to be the isometry acting 
pointwise as $J_x$ on  $T_x$, namely $(J\mathbf{f})(x)=J_x(\mathbf{f}(x))$ for $\mathbf{f}\in\T_R$. 
It turns out that
$J(\G_R)=\D_R$. This fact holds for more general sufficiently smooth surfaces embedded in $\R^3$. A proof seems not easy to find in the literature and is provided in Appendices \ref{sec:appendix2} and \ref{sec:appendix1} (for the special case  of continuously differentiable tangential vector fields on
 the sphere, the assertion essentially corresponds to 
\cite[Thm. 2.10]{freedenschreiner09}). This motivates the notion of a surface curl gradient $\LL_{\SS_R}=x\times\nabla_{\SS_R}$, acting at a point $x\in\SS_R$, and justifies the representation $\D_R=\{\LL_{\SS_R} f:f\in W^{1,2}(\SS_R)\}$. For convenience, we define the following ''normalized'' operators: $\nabla_\SS=R\nabla_{\SS_R}$ and $\LL_\SS=\frac{x}{|x|}\times \nabla_\SS$. The Euclidean gradient then has the expression $\nabla = \frac{x}{|x|}\partial_\nu+\frac{1}{|x|}\nabla_\SS$, acting at a point $x\in\R^3$, where $\partial_\nu=\frac{x}{|x|}\cdot\nabla$ denotes the radial derivative.

Eventually, if we let $\N_R$ indicate the space of radial vector fields in
$L^2(\SS_R,\R^3)$ (i.e., those functions whose value at $x$ is perpendicular to $T_x$ for each
$x\in\SS_R$), we get from \eqref{eqn:hd} the orthogonal decomposition
\begin{align}
 L^2(\SS_R,\R^3)=\N_R\oplus\G_R\oplus \D_R.\label{eqn:shd}
\end{align}
Related to the latter but of more relevance to our problem is the 
Hardy-Hodge decomposition that we now explain. 
For that purpose, we require the following definition.

\begin{defi}\label{def:HH}
 The Hardy space $\Hp$ of harmonic gradients in $\BB_R$ 
is defined by
 \begin{align*}
  \Hp=\left\{\mathbf{g}=\nabla g:\textnormal{ function  } g:\BB_R\to\R \textnormal{ with } \Delta g=0 \textnormal{ in }\BB_{R}\textnormal{ and }\|\nabla{g}\|_{2,+}<\infty\right\},
 \end{align*}
 where $\|\mathbf{g}\|_{2,+}=\big(\sup_{r\in[0,R)}\int_{\SS_r}|\mathbf{g}(ry)|^2\dd \omega_r(y)\big)^{\frac{1}{2}}$ and $\Delta$ is the Euclidean Laplacian in $\R^3$.  Likewise, the Hardy space $\Hm$ of harmonic gradients in $\R^3\setminus\overline{\BB_R}$ is defined by
 \small\begin{align*}
  \Hm=\left\{\mathbf{g}=\nabla g: \textnormal{ function  } g:\R^3\setminus\overline{\BB_R}\to\R \textnormal{ with } \Delta g=0 \textnormal{ in }\R^3\setminus\overline{\BB_R}\textnormal{ and }\|\nabla{g}\|_{2,-}<\infty\right\},
 \end{align*}\normalsize
 where $\|\mathbf{g}\|_{2,-}=\big(\sup_{r\in(R,\infty)}\int_{\SS_r}|\mathbf{g}(ry)|^2\dd \omega_R(y)\big)^{\frac{1}{2}}$. Note that, by Weyl's lemma
\cite[Theorem 24.9]{Forster},
it makes no difference whether
the Euclidean gradient and Laplacian are understood in 
the distributional or in the strong sense.
\end{defi}

Members of $\Hp$ and $\Hm$ have non-tangential limits a.e. on $\SS_R$,
and if $\mathbf{g}\in\H^2_{\pm,R}$, its nontangential 
limit has $L^2(\SS_R, \R^3)$-norm equal to $\|\mathbf{g}\|_{2,\pm}$, see
\cite[VII.3.1]{stein70} and \cite[VI.4]{steinweiss71}. 
We still write $\mathbf{g}$ for this non-tangential limit and we regard 
it as the trace of $\mathbf{g}$ on $\sph_R$. This way
Hardy spaces can be interpreted as function spaces on $\SS_R$ as well as on 
$\BB_{R}$ or $\RR^3\setminus\overline{\BB_{R}}$,
but the context will make it clear if the Euclidean or the spherical interpretation is meant because the argument belongs to  $\R^3\setminus\SS_R$
in the former case and to $\SS_R$ in the latter.
The Hardy-Hodge decomposition is the orthogonal sum 
\begin{align}
 L^2(\SS_R,\R^3)=\Hp\oplus \Hm\oplus \D_R.\label{eqn:hhd}
\end{align}
Projecting \eqref{eqn:hhd} onto the tangent space $\T_R$ and grouping the first two summands into a single gradient vector field yields back the Hodge decomposition \eqref{eqn:hd}. The Hardy-Hodge decomposition drops out at once from \cite{atfeh10} and \eqref{eqn:hd}. Its application to the study of inverse magnetization problems has been illustrated in \cite{baratchart13,gerhards16a,lima13}. Although not studied in mathematical detail, spherical versions of the Hardy-Hodge decomposition have previously been used to a various extent in Geomagnetic applications (see, e.g., \cite{backus96,gerhards12,gubbins11,mayer06}).

By means of the reflection $\mathcal{R}_R(x)=\frac{R^2}{|x|^2}\,x$ across $\SS_R$, we define the Kelvin transform ${K}_R[f]$ of a function $f$ defined on an open set $\Omega\subset\RR^3$ to be the function on $\mathcal{R}_R(\Omega)$ given by
\begin{align}\label{defKelvin}
{K}_R[f](x)= \frac{R}{|x|} f(\mathcal{R}_R(x)), \quad  x\in \mathcal{R}_R(\Omega).
\end{align}
A function $f$ is harmonic in $\Omega$ if and only if ${K}_R[f]$ is harmonic 
in $\mathcal{R}_R(\Omega)$ (e.g., \cite[Thm. 4.7]{axler01}).

Now, assume that   $\mathbf{f}\in \Hp$  with $\mathbf{f}=\nabla f$ and $f(0)=0$. Then $\nabla{K}_R[f]\in \Hm$.  In fact, if for $\mathbf{f}\in \Hp$ (resp. $\mathbf{f}\in \Hm$) we let $\int\mathbf{f}$ indicate the harmonic function $f$ in $\BB_R$ (resp. in $\RR^3\setminus\overline{\BB}_R$) whose gradient is $\mathbf{f}$, normalized so that $f(0)=0$ (resp. $\lim_{|x|\to\infty}f(x)=0$), then 
$\mathbf{f}\mapsto \nabla{K}_R\circ\int\mathbf{f}$ maps   $ \Hp$ continuously into  $ \Hm$ and back \cite{atfeh10}. Moreover, in view of (\ref{defKelvin}) we have that
\begin{align}
\label{gradKel}
\nabla {K}_R[f](x) = \frac{R^3\nabla f(\mathcal{R}_R(x))}{|x|^3} - 2\, x \cdot \nabla f(\mathcal{R}_R(x)) \, \frac{R^3x}{|x|^{5}} - \, f(\mathcal{R}_R(x)) \,  \frac{Rx}{|x|^3} \, .
\end{align}
Clearly $f$ and  ${K}_R[f]$ coincide on $\SS_R$, therefore the  tangential components of $\nabla f$ and $\nabla {K}_R[f]$ agree on $\SS_R$ (these are the spherical gradients $\nabla_{\SS_R}f$ and $\nabla_{\SS_R}{K}_R[f]$). The normal components $\partial_\nu f$ and $\partial_\nu {K}_R[f]$, though, are different. Indeed, we get from \eqref{gradKel} that
\begin{align}
\label{diffnK}
\partial_\nu {K}_R[f] (x) = -\partial_\nu f(x) - \frac{f(x)}{R},\quad x\in\SS_R.
\end{align}

We turn to some special systems of functions. First, let $\{Y_{n,k}\}_{n\in\mathbb{N}_0,\,k=1,\ldots,2n+1}$ be an $L^2(\SS)$-orthonormal system of spherical harmonics of degrees $n$ and orders $k$. A possible choice is 
\small\begin{align*}
 Y_{n,k}(x)=\left\{\begin{array}{ll}
                    \sqrt{\frac{2n+1}{2\pi}\frac{(k-1)!}{(2n+1-k)!}}P_{n,n+1-k}(\sin(\theta))\cos((n+1-k)\varphi)&k=1,\ldots,n,
                    \\\sqrt{\frac{2n+1}{4\pi}}P_{n,0}(t),&k=n+1,
                    \\\sqrt{\frac{2n+1}{2\pi}\frac{(2n+1-k)!}{(k-1)!}}P_{n,k-(n+1)}(\sin(\theta))\sin((k-(n+1))\varphi)&k=n+2,\ldots,2n+1,
                   \end{array}\right.
\end{align*}\normalsize
for $x=(\cos(\theta)\cos(\varphi),\cos(\theta)\sin(\varphi),\sin(\theta))^T\in\SS_1$, $\theta\in[-\frac{\pi}{2},\frac{\pi}{2}],\varphi\in[0,2\pi)$, and $P_{n,k}$ the associated Legendre polynomials of degree $n$ and order $k$ (see, e.g., \cite[Ch. 3]{freedenschreiner09} for details; another common notation is to indicate the order of the spherical harmonics by  $k=-n,\ldots,n$ rather than $k=1,\ldots,2n+1$). Then $H_{n,k}^R(x)=\big(\frac{|x|}{R}\big)^nY_{n,k}\big(\frac{x}{|x|}\big)$ is a homogeneous, harmonic polynomial of degree $n$ in $\RR^3$ (sometimes also called inner harmonic and equipped with a normalization factor $\frac{1}{R}$). In fact, every homogeneous harmonic polynomial in $\R^3$ can be expressed as a linear combination of inner harmonics. The Kelvin transform $H_{-n-1,k}^R={K}_R[H_{n,k}^R]$ is a harmonic function in $\RR^3\setminus\{0\}$ with $\lim_{|x|\to\infty}H_{-n-1,k}^R(x)=0$ (sometimes called outer harmonic). In \cite[Lemma 4]{atfeh10} the following result 
was shown.

\begin{lem}\label{lem:yndense}
The vector space $\textnormal{span} \{\nabla H_{-n-1,k}^R\}_{n\in\mathbb{N}_0,\,k=1,\ldots,2n+1}$ is dense in $\Hm$ and the vector space  $\textnormal{span} \{\nabla H_{n,k}^R\}_{n\in\mathbb{N}_0,\,k=1,\ldots,2n+1}$ is dense in $\Hp$.
\end{lem}

For each fixed $x\in\R^3\setminus\overline{\BB}_R$, the function $g_x(y)= \frac{1}{|x-y|}$ is harmonic in a neighborhood of $\overline{\BB}_R$ and, therefore, its gradient 
\begin{align*}
 \mathbf{g}_x(y)=\nabla_x\,\, g_x(y)=-\frac{x-y}{|x-y|^3}
 \end{align*}
lies in $\Hp$. As a consequence of Lemma \ref{lem:yndense}, we shall
prove the following density result.

\begin{lem}\label{dsh}
The vector space $\textnormal{span}\{\mathbf{g}_x: x\in \RR^3\setminus\overline{\BB}_R\}$ is dense in $\Hp$ and  the vector space $\textnormal{span}\{\mathbf{g}_x: x\in \BB_R\}$ is dense
in $\Hm$.
\end{lem}

\begin{prf}  
As ${{K}_R}[g_x]= \frac{1}{|x|}g_{x/|x|^2}$ and $\nabla{{K}_R}\circ \int$ is an isomorphism from $\Hm$ onto $\Hp$ (see discussion before \eqref{gradKel}), we need
only prove the second assertion. Define $g(y)=\frac{1}{|y|}$ as a function of
$y\in\RR^3\setminus\{0\}$. For  $\alpha=(\alpha_1,\alpha_2,\alpha_3)\in\mathbb{N}_0^3$ with $|\alpha|=\alpha_1+\alpha_2+\alpha_3=n$, the derivative $\partial_\alpha g(y)=\frac{\partial^{n}}{\partial^{\alpha_1}y_1 \partial^{\alpha_2}y_2\partial^{\alpha_3}y_3}g(y)$ is of the form $\frac{H_\alpha(y)}{|y|^{1+2n}}$, where $H_\alpha$ is a homogeneous harmonic polynomial of degree $n$, and
actually every homogeneous harmonic polynomial $H_\alpha$ is
a scalar multiple of $|y|^{(1+2n)}\partial_\alpha g(y)$ for some $\alpha$
\cite[Lemma 5.15]{axler01}. The discussion before Lemma 
\ref{lem:yndense} 
now implies that $\partial_\alpha g$ is an element of $\textnormal{span} \{H_{-n-1,k}^R\}_{n\in\mathbb{N}_0,\,k=1,\ldots,2n+1}$.  Thus, by
this lemma,  we are done if we can show that 
whenever  $\mathbf{f}\in \Hm$ is orthogonal in $L^2(\SS_R,\RR^3)$ to all $\mathbf{g}_x$, $ x\in \BB_R$, then it must be orthogonal to all $\nabla H_{-n-1,k}^R$. 
To this end, differentiating $\langle \mathbf{f},\mathbf{g}_x\rangle_{L^2(\SS_R,\R^3)}=0$ with respect to $x$ leads us to 
\begin{align}
\label{orthT}
0=\left\langle \mathbf{f}\,,\, \nabla \frac{H_\alpha(.-x)}{|\cdot-x|^{1+2n}}\right\rangle_{L^2(\SS_R,\R^3)}
\end{align}
for all $\alpha\in\mathbb{N}_0^3$ and $n=|\alpha|$. Setting $x=0$ yields
\begin{align*}
0=\left\langle  \mathbf{f}\,,\, \nabla \frac{H_\alpha}{|\cdot|^{1+2n}} \right\rangle_{L^2(\SS_R,\R^3)}=R^{-2n-1}\left\langle  \mathbf{f}\,,\, \nabla {{K}_R}[H_\alpha] \right\rangle_{L^2(\SS_R,\R^3)}.
\end{align*}
Since every inner harmonic $H_{n,k}^R$ can be expressed as a linear combination of $H_\alpha$, this relation and the considerations before Lemma \ref{lem:yndense} imply $\langle  \mathbf{f}\,,\, \nabla H_{-n-1,k}^R \rangle_{L^2(\SS_R,\R^3)}=0$ for all $n\in\mathbb{N}_0$, $k=1,\ldots, 2n+1$, which is the desired conclusion.
\end{prf}

\section{Harmonic Potentials in Divergence-Form}\label{sec:harmpot}
 
The potential of a measure $\mu$ on $\RR^3$ is defined by
\begin{align}
p_\mu(x)=-\frac{1}{4\pi}\int_{\R^3}\frac{1}{|x-y|}\dd\mu(y)\label{eqn:defpot}.
\end{align}
It is  the solution of $\Delta \Phi=\mu$  in $\RR^3$ which is ``smallest'' at infinity. If $\mu\geq0$, the potential $p_\mu$ is a superharmonic function and therefore it is either finite quasi-everywhere or identically $-\infty$,
see \cite{armitage01} for these properties and
the definition of ``quasi everywhere''. 
Decomposing a signed measure into its positive and negative parts 
(the Hahn decomposition) yields that $p_\mu$ is finite quasi-everywhere if $\mu$ is finite and compactly supported (i.e., if $\textnormal{supp}(\mu)$, which is closed by definition, is also bounded). If $\textnormal{supp}(\mu)\subset\overline{\BB}_R$, the Riesz representation theorem and the maximum principle for harmonic functions imply that there exists a unique measure $\hat\mu$ with $\textnormal{supp}(\hat\mu)\subset\SS_R$ such that \begin{align*}
 \int g(y) \dd\mu(y)=\int g(y) \dd\hat\mu(y)
 \end{align*}
for every continuous function $g$ in $\overline{\BB}_R$ which is harmonic in $\BB_R$. Since $y\mapsto 1/|x-y|$ is harmonic
in a neighbourhood of $\overline{\BB}_R$ when $x\notin \overline{\BB}_R$,
this entails that the potentials $p_\mu$ and $p_{\hat\mu}$ coincide in $\R^3\setminus\overline{\BB}_R$, i.e.,
\begin{align*}
 p_\mu(x)=p_{\hat\mu}(x),\quad x\in\R^3\setminus\overline{\BB}_R.
\end{align*}
The measure $\hat\mu$ is called the balayage of $\mu$ onto $\sph_R$ 
(see, e.g., \cite{armitage01}). In fact, the potentials $p_\mu$ and $p_{\hat{\mu}}$ coincide quasi-everywhere on $\sph_R$ as well. An expression for $\hat\mu$ easily follows from the Poisson representation of a function $f$ which is continuous in $\overline{\BB}_R$ and harmonic in $\BB_R$: 
\begin{align} \label{Poisson}
f(x)=\frac{1}{4\pi R}\int_{\SS_R}\frac{R^2-|x|^2}{|x-y|^3}f(y)\,\dd \omega_R(y),\quad x\in\BB_R.
\end{align}
Clearly Equation \eqref{Poisson}, Fubini's theorem  and the definition of balayage imply that
\begin{align}\label{calcbal}
d\hat\mu(x)=\dd\mu_{|\sph_R}(x)+\left(\frac{1}{4\pi R}\int_{\BB_R}\frac{R^2-|y|^2}{|x-y|^3}\dd\mu(y)\right)\dd \omega_R(x).
\end{align}

\begin{lem}\label{lem:balayage}
Let the measure $\mu$ be supported in $\overline{\BB_R}$. Furthermore, assume that $\mu$ is absolutely continuous in $\BB_R$ with  a density $h$ (i.e., $\dd\mu(y)=h(y)dy$) that satisfies the Hardy condition
\begin{align}
\esssup_{0 \le r < R} \int_{\sph_r} |h(y)|^2 \, \dd \omega_r(y) < \infty.\label{eqn:esssup}
\end{align}
Then the balayage $\hat\mu$ of $\mu$ on $\SS_R$ is absolutely continuous with respect to $\omega_R$  (i.e., $d\hat\mu(y)=\hat h(y)\dd \omega_R(y)$) and it has a density $\hat h\in L^2(\sph_R)$.
\end{lem}

\begin{prf}
Starting from \eqref{calcbal} and the assumption that $\mu$ is absolutely continuous, we find that the density $\hat h$ of $\hat \mu$ is
\begin{align*}
 \hat h(x)=\frac{1}{4\pi R}\int_{\BB_R}\frac{R^2-|y|^2}{|x-y|^3}h(y)\dd\lambda(y), \quad x\in \SS_R.
\end{align*}
Using Fubini's theorem and the identity
\begin{align*}
\left|\frac{x}{|x|}-|x|y\right|=\left|\frac{y}{|y|}-|y|x\right|, \qquad x,y\in\RR^3\setminus\{0\},
\end{align*}
together with the changes of variable 
$\eta=\frac{\xi}{r}$, $y=\frac{rx}{R^2}$, we are led to
\begin{align}
\|\hat h\|^2_{L^2(\sph_R)}=&\frac{1}{(4\pi R)^2}\int_{\sph_R}\left(\int_{\BB_R}\frac{R^2-|y|^2}{|x-y|^3}h(y)\,\dd\lambda(y)\right)^2\dd \omega_R(x)\nonumber
\\
=&\frac{1}{(4\pi R)^2}\int_{\sph_R}\left(\int_0^R\left(\int_{\sph_r}\frac{R^2-|\xi|^2}{|x-\xi|^3}h(\xi)\,\dd \omega_r(\xi)\right)dr\right)^2\dd \omega_R(x)\nonumber
\\
\leq&\frac{R}{(4\pi R)^2}\int_{\sph_R}\left(\int_0^R\left(\int_{\sph_r}\frac{R^2-|\xi|^2}{|x-\xi|^3}h(\xi)\,\dd \omega_r(\xi)\right)^2dr\right)\dd \omega_R(x)\nonumber
\\
=&\frac{1}{(4\pi R)^2}\int_{\sph_R}\left(\int_0^R\left(\int_{\sph_r}\frac{1-(\frac{r}{R})^2}{|\frac{x}{R}-\frac{\xi}{R}|^3}h(\xi)\,\dd \omega_r(\xi)\right)^2dr\right)\dd \omega_R(x)\nonumber
\\
=&\frac{1}{(4\pi R)^2}\int_{\sph_R}\left(\int_0^R\left(\int_{\sph_r}\frac{1-\left|\frac{rx}{R^2}\right|^2}{|\frac{\xi}{r}-\frac{rx}{R^2}|^3}h(\xi)\,\dd \omega_r(\xi)\right)^2dr\right)\dd \omega_R(x)\nonumber
\\
=&\frac{1}{(4\pi R)^2}\int_0^Rr^4\left(\int_{\sph_R}\left(\int_{\sph_1}\frac{1-\left|\frac{rx}{R^2}\right|^2}{|\eta-\frac{rx}{R^2}|^3}h(r\eta)\,\dd \omega_1(\eta)\right)^2 \dd \omega_R(x)\right)dr\nonumber
\\
=&\int_0^Rr^4\left(\frac{1}{4\pi (\frac{r}{R})^2}\int_{\sph_{\frac{r}{R}}}\left(\frac{1}{4\pi}\int_{\sph_1}\frac{1-|y|^2}{|\eta-y|^3}h(r\eta)\,\dd \omega_1(\eta)\right)^2 \dd \omega_{\frac{r}{R}}(y)\right)dr.\label{eqn:hatgcomp}
\end{align}
Now, the function 
\begin{align*}
f(y)=\frac{1}{4\pi}\int_{\sph_1}\frac{1-|y|^2}{|\eta-y|^3}h(r\eta)\,\dd \omega_1(\eta)
\end{align*}
is the Poisson integral of $h(r\cdot)$ over the unit sphere $\SS_1$ (and represents the middle integral on the right hand side of \eqref{eqn:hatgcomp}). Thus, $f$ is harmonic in $\BB_1$ and its square $|f|^2$ is subharmonic there. 
The latter implies that the mean of $|f|^2$ over the sphere $\sph_{\frac{r}{R}}$, $r<R$, is not greater than its mean over $\sph_1$, i.e., 
\begin{align*}
\frac{1}{4\pi (\frac{r}{R})^2}\int_{\sph_{\frac{r}{R}}}|f(y)|^2 \dd \omega_{\frac{r}{R}}(y)&\leq \lim_{\frac{s}{R}\to 1-}\frac{1}{4\pi (\frac{s}{R})^2}\int_{\sph_{\frac{s}{R}}}|f(y)|^2 \dd \omega_{\frac{s}{R}}(y)= \frac{1}{4\pi}\int_{\sph_1}|h(r\eta)|^2\,\dd \omega_1(\eta) 
\\&=\frac{1}{4\pi r^2}\int_{\sph_r}|h(y)|^2\,\dd \omega_r(y)\leq \frac{M}{4\pi r^2},
\end{align*}
where the constant $M>0$ comes from the Hardy condition \eqref{eqn:esssup}. Together with \eqref{eqn:hatgcomp}, we find that
\begin{align*}
\|\hat h\|^2_{L^2(\sph_R)}\leq \frac{MR^3}{12\pi},
\end{align*}
eventually showing that $\hat{h}\in L^2(\sph_R)$ and that $\hat\mu$ is absolutely continuous with respect to $\omega_R$ with density $\hat h$.
\end{prf}

More generally, an arbitrary distribution $D$ with compact support
has a potential $p_D$ given outside of $\textnormal{supp}(D)$ by
\begin{align}\label{eqn:p_Ddef}
 p_D(x)=D\left(-\frac{1}{4\pi}\frac{1}{|x-\cdot|}\right), \quad x\in\R^3\setminus \textnormal{supp} D.
\end{align}
Compactness of $\textnormal{supp} (D)$ easily implies that $D$ indeed acts on
$-1/(4\pi|x-\cdot|)$ when $x\notin \textnormal{supp} D$ so that $P_D$ is well-defined, see Appendix \ref{sec:appendix}
for details.  
If $D$ is supported in $\overline{\BB}_R$ (in particular, if it is supported in some shell $\overline{\BB_{R-d,R}}$), we define the balayage of $D$ onto $\sph_R$ to be the distribution $\hat D$ on 
$\SS_R$ that satisfies 
\begin{align*}
 p_{\hat D}(x)=p_D(x),\quad x\in\R^3\setminus\overline{\BB}_R.
\end{align*}
Strictly speaking, $\hat{D}$ is a distribution on $\SS_R$, so that
$p_{\hat D}$ should rather be denoted by
$p_{\hat D\otimes \delta_{\SS_R}}$, where
$\hat D\otimes\delta_{\SS_R}$ is the distribution
on $\RR^3$ which is the tensor product of $\hat D$ with
the measure $\delta_{\SS_R}$, corresponding in spherical coordinates 
to  a Dirac mass at $r=R$, see \cite{Schwartz}. Nevertheless, to alleviate 
notation, we do write $p_{\hat D}$. Thus, what is meant  in  
\eqref{eqn:p_Ddef} when $D=\hat D$ is that $\hat D$ is applied to the restriction to
$\SS_R$ of $ -1/(4\pi|x-\cdot|)$.

We briefly comment on the existence and uniqueness of such a balayage in Appendix \ref{sec:appendix}. 
If $D$ is (associated with) a measure $\mu$, then \eqref{eqn:p_Ddef} coincides with \eqref{eqn:defpot} and the balayage was given in  \eqref{calcbal}. 
The main difference between the case of a finite compactly supported measure
$\mu$ and the case of a general compactly supported distribution $D$ is that 
usually $p_D(x)$
cannot be assigned a meaning when $x\in \textnormal{supp}(D)$ whereas $p_\mu$ is well-defined quasi everywhere on $\textnormal{supp}(\mu)$.
We say that $D$ is in divergence form if 
\begin{align}\label{eqn:divform}
D=\nabla\cdot \mathbf{M},
\end{align}
where $\nabla\cdot$ is to be understood as the distributional divergence 
and $\mathbf{M}$  is a $\RR^3$-valued distribution. If, e.g., $\mathbf{M}\in L^2(\BB_{R-d,R},\R^3)$ and $\textnormal{supp}(\mathbf{M})\subset\overline{\BB_{R-d,R}}$, then the corresponding potential $p_D$ coincides with $\Phi_{crust}$ in \eqref{eqn:phivol}. Now we can formulate the main result of this section, namely, that balayage preserves divergence form for those $\mathbf{M}$ satisfying a Hardy condition.

\begin{lem}\label{balayageH}
Let $D=\nabla\cdot \mathbf{M}$, where $\mathbf{M}\in L^2(\BB_R,\RR^3)$ satisfies the Hardy condition
\begin{align*}
\esssup_{0 \leq r < R} \int_{\sph_r} |\mathbf{M}(y)|^2 \, \dd \omega_r(y) < \infty.
\end{align*}
Then there exists ${\mathbf{m}}\in L^2(\sph_R,\RR^3)$ such that $\hat D=\nabla\cdot ({\mathbf{m}}\otimes \delta_{\sph_R})$ is the  balayage of $D$ onto $\sph_R$. 
\end{lem}

\begin{prf}
Let $\mathbf{M}=(M_1,M_2,M_3)^T$ denote the components of $\mathbf{M}$. The definition of $p_D$ yields
\begin{align}
p_D(x)&=\frac{1}{4\pi}\int_{\BB_{R}}\M(y)\cdot\frac{x-y}{|x-y|^3}\,\dd\lambda(y)\nonumber
\\&=\frac{1}{4\pi}\sum_{j=1}^3\int_{\BB_R} M_j(y)\frac{x_j-y_j}{|x-y|^3}  \, \dd\lambda(y),\quad x\in\R^3\setminus\overline{\BB}_R.\label{eqn:pD}
\end{align}
If we choose the measure $\mu_j$ such that $\dd\mu_j(y)=M_j(y)dy$, we get from Lemma \ref{lem:balayage} and the Hardy condition on $\mathbf{M}$ that there exists a  $m_j\in L^2(\sph_R)$ such that balayage of $\mu_j$ onto $\SS_R$ is given by the measure $\hat \mu_j$ with $d\hat\mu_j=m_j \dd \omega_R$, $j=1,2,3$. Setting ${\mathbf{m}}=(m_1,m_2,m_3)^T$ and observing that ${g}_{x,j}(y)=\frac{x_j-y_j}{|x-y|^3}=-\partial_{x_j}\frac{1}{|x-y|}$ is harmonic in $\BB_R$ and continuous in $\overline{\BB_R}$, for fixed $x\in \R^3\setminus\overline{\BB}_R$, then the definition of balayage yields together with \eqref{eqn:pD} that
\begin{align}
p_D(x)&=\frac{1}{4\pi}\sum_{j=1}^3\int_{\sph_R}m_j(t)\frac{x_j-y_j}{|x-y|^3}\, \dd \omega_R(y)\nonumber
\\&=\frac{1}{4\pi}\int_{\SS_{R}}\m(y)\cdot\frac{x-y}{|x-y|^3}\,\dd \omega_R(y)=p_{\hat D}(x),\quad x\in\R^3\setminus\overline{\BB}_R.\label{eqn:spherepot}
\end{align}
The latter implies that $\hat D=\nabla\cdot ({\mathbf{m}}\otimes \delta_{\sph_R})$, as announced.
\end{prf}

\begin{rem}
Lemma \ref{balayageH} eventually justifies the statement made in
the introduction that, to every square summable
volumetric magnetization $\mathbf{M}$ 
in the Earth's crust $\BB_{R-d,R}$ that satisfies the Hardy condition, there exists a spherical magnetization $\mathbf{m}$ on $\SS_{R}$ that produces the same magnetic potential and therefore also the same magnetic field in the exterior of the Earth. 
\end{rem}

\section{Separation of Potentials}\label{sec:seppot}
We are now in a position to approach Problem \ref{prob:1}. 
For this we study the nullspace of the potential operator 
$\Phi^{R_1,R_0,R_2}$ (cf. Definition \ref{def:ops}), mapping a
magnetization $\mathbf{m}$ on $\SS_{R_0}$ and an auxiliary 
function $h\in L^2(\SS_{R_1})$ to the sum of the potentials 
\eqref{eqn:phisurf} and \eqref{eqn:phicore}
on $\SS_{R_2}$.
First, we show in Section \ref{sec:unique}
that uniqueness holds in Problem \ref{prob:1} if 
$\textnormal{supp}\,\mathbf{m}\not= \SS_{R_0}$. 
Similar results hold for the  magnetic field operator $\mathbf{B}^{R_1,R_0,R_2}=\nabla\Phi^{R_1,R_0,R_2}$ (cf. Theorem \ref{thm:unique3}), 
and also for a modified potential field operator 
$\Psi^{R_1,R_0,R_2}$ (cf. Definition \ref{def:ops2} and Theorem \ref{thm:unique2}). The operator $\Psi^{R_1,R_0,R_2}$ 
reflects the potential of two magnetizations $\mathbf{m}$ and $\mathbf{\mathfrak{m}}$ supported on two different spheres (i.e., at different depths),
therefore it does not apply to the separation of the crustal and core 
contributions since the latter does not arise from a magnetization (cf.\eqref{Maxwell}).  Still it is of  interest on its own,
moreover we get it at no extra cost.
In Section \ref{sec:extremal}, we discuss how the previous results can be used
to approximate quantities like the Fourier coefficients $\langle \Phi_0,Y_{n,k}\rangle_{L^2(\SS_{R_2})}$ of $\Phi_0$. Finally, in Section \ref{sec:gammas},
we show that  $\Phi=\Phi_0+\Phi_1$ 
may well vanish though $\Phi_0,\Phi_1\neq0$. This follows from 
Lemma \ref{NG} and answers 
the uniqueness issue of Problem \ref{prob:1} in the negative when 
$\textnormal{supp}\,\mathbf{m}= \SS_{R_0}$.

\subsection{Uniqueness Issues}\label{sec:unique}
In accordance with the notation from Problem \ref{prob:1}, we define 
two operators: one mapping a spherical magnetization $\mathbf{m}$ to the potential $p_{\hat D}$ with $\hat D=\nabla\cdot (\mathbf{m}\otimes \delta_{\SS_{R_0}})$, and the other mapping  an auxiliary function $h\in L^2(\SS_{R_1})$ to its Poisson integral,
both evaluated on $\SS_{R_2}$.

\begin{defi}\label{def:ops}
 Let $0<R_1<R_0<R_2$ be fixed radii and $\Gamma_{R_0}$ a 
closed subset of $\SS_{R_0}$.
Let
\begin{align*}
  \Phi_0^{R_0,R_2}:L^2(\Gamma_{R_0},\R^3)\to L^2(\SS_{R_2}),\quad \mathbf{m}\mapsto\frac{1}{4\pi}\int_{\Gamma_{R_0}}\m(y)\cdot\frac{x-y}{|x-y|^3}\,\dd \omega_{R_0}(y),\quad x\in\SS_{R_2},
 \end{align*}
and
\begin{align*}
  \Phi_1^{R_1,R_2}:L^2(\SS_{R_1})\to L^2(\SS_{R_2}),\quad h\mapsto\frac{1}{4\pi R_1}\int_{\SS_{R_1}} h(y)\frac{|x|^2-R_1^2}{|x-y|^3}\dd \omega_{R_1}(y),\quad x\in\SS_{R_2}.
 \end{align*}
 The superposition of the two operators above is denoted by 
 \begin{align*}
  \Phi^{R_1,R_0,R_2}:L^2(\Gamma_{R_0},\R^3)\times L^2(\SS_{R_1})\to L^2(\SS_{R_2}),\quad(\mathbf{m},h)\mapsto  \Phi_0^{R_0,R_2}[\mathbf{m}]+ \Phi_1^{R_1,R_2}[h].
 \end{align*}
\end{defi}

We start by characterizing the potentials $p_{\hat D}$, with $\hat D$ in divergence-form, which are zero in $\R^3\setminus\overline{\BB_R}$.

\begin{lem}\label{lem:nullpd}
Let $\mathbf{m}\in L^2(\SS_R,\R^3)$ and ${\hat D}=\nabla\cdot(\mathbf{m}\otimes\delta_{\SS_R})$ be in divergence-form. Let further $\mathbf{m}=\mathbf{m}_{+}+\mathbf{m}_{-}+\mathbf{d}$ be the Hardy-Hodge decomposition of $\mathbf{m}$, i.e., $\mathbf{m}_{+}\in \Hp$, $\mathbf{m}_{-}\in \Hm$, and $\mathbf{d}\in \D_R$. Then $p_{\hat D}(x)=0$, for all $x\in \R^3\setminus\overline{\BB_R}$, if and only if $\mathbf{m}_+\equiv0$. Analogously, $p_{\hat D}(x)=0$, for all $x\in \BB_R$, if and only if $\mathbf{m}_-\equiv0$.
\end{lem}

\begin{prf}
 We already know that $\mathbf{g}_x(y)=\frac{x-y}{|x-y|^3}$ lies in $\Hp$ for every fixed $x\in\R^3\setminus\overline{\BB_R}$. The orthogonality of the Hardy-Hodge decomposition and the representation \eqref{eqn:spherepot} of $p_{\hat D}$ yield that $\mathbf{m}_-$ and $\mathbf{d}$ do not change $p_{\hat D}$ in $\R^3\setminus\overline{\BB_R}$. Conversely, if $p_{\hat D}(x)=0$ for all $x\in \R^3\setminus\overline{\BB_R}$, then 
 \begin{align*}
  p_{\hat D}(x)=\langle \mathbf{g}_x,\mathbf{m}\rangle_{L^2(\SS_R,\R^3)}=\langle \mathbf{g}_x,\mathbf{m}_+\rangle_{L^2(\SS_R,\R^3)}=0,\quad x\in\R^3\setminus\overline{\BB_R}.
 \end{align*}
 Since Lemma \ref{dsh} asserts that $\textnormal{span}\{\mathbf{g}_x: x\in \RR^3\setminus\overline{\BB}_R\}$ is dense in $\Hp$, the above relation implies $\mathbf{m}_+\equiv0$. The assertion  for the case where 
$p_{\hat D}(x)=0$, for all $x\in \BB_R$ likewise follows by observing that $\mathbf{g}_x(y)=\frac{x-y}{|x-y|^3}$ lies in $\Hm$ for fixed $x\in\BB_R$.
\end{prf}

Since  $\Phi_0^{R_0,R_2}[\mathbf{m}]=p_{\hat D}$, we may use 
Lemma \ref{lem:nullpd} to characterize the nullspace of 
$\Phi_0^{R_0,R_2}$ (extending the magnetization $\mathbf{m}\in L^2(\Gamma_{R_0},\R^3)$ by zero on $\SS_{R_0}\setminus\Gamma_{R_0}$ if the latter 
is nonempty). 
As to $\Phi^{R_1,R_2}_1$, we know its  nullspace 
reduces to zero because the Poisson integral  \eqref{eqn:phicore}
yields the unique harmonic 
extension of $h\in L^2(\SS_{R_1})$ to  
 $\R^3\setminus\overline{\BB_{R_1}}$ which is zero at infinity
( i.e.  $h$ is the nontangential limit of its Poisson extension
a.e. on $\SS_{R_1}$, see \cite[Thm. 6.13]{axler01}).
This motivates the following statement on the nullspace 
$N(\Phi^{R_1,R_0,R_2})$ of $\Phi^{R_1,R_0,R_2}$. 

\begin{thm}\label{thm:unique}
 Let the setup be as in Definition \ref{def:ops} and assume that
$\Gamma_{R_0}\not=\SS_{R_0}$. Then the nullspace of $\Phi^{R_1,R_0,R_2}$ is given by
 \begin{align*}
  N(\Phi^{R_1,R_0,R_2})=\{(\mathbf{d},0):\mathbf{d}\in \mathcal{D}_{R_0},\ 
\textnormal{supp}(\mathbf{d})\subset\Gamma_{R_0}\}.
 \end{align*}
\end{thm}

\begin{prf}
Clearly 
$\Phi^{R_1,R_0,R_2}[(\mathbf{m},h)]$ is harmonic in 
$\R^3\setminus\{\Gamma_{R_0}\cup \,\SS_{R_1}\}$ and vanishes at infiity.
If  
$\Phi^{R_1,R_0,R_2}[(\mathbf{m},h)](x)=0$ for $x\in\SS_{R_2}$,  
then it 
follows from the maximum principle 
that 
$\Phi^{R_1,R_0,R_2}[(\mathbf{m},h)](x)=0$ for all 
$x\in \R^3\setminus\BB_{R_2}$.  Subsequently, by real analyticity,
$\Phi^{R_1,R_0,R_2}[(\mathbf{m},h)]$
must vanish identically in 
$\R^3\setminus\{\Gamma_{R_0}\cup \overline{\BB_{R_1}}\}$ which is 
connected because $\Gamma_{R_0}\not=\SS_{R_0}$. Thus,
$\Phi^{R_1,R_0,R_2}[(\mathbf{m},h)]$ extends harmonically 
(by the zero function) across $\Gamma_{R_0}$:
 \begin{align}\label{eqn:zeropot}
  \Phi^{R_1,R_0,R_2}[(\mathbf{m},h)](x)=0,\quad x\in \R^3\setminus\overline{\BB_{R_1}}. 
 \end{align}
Since $\Phi^{R_1,R_0,R_2}[(\mathbf{m},h)]= 
\Phi_0^{R_0,R_2}[\mathbf{m}]+ \Phi_1^{R_1,R_2}[h]$, where $\Phi_1^{R_1,R_2}[h]$ is harmonic on $\R^3\setminus\overline{\BB_{R_1}}$, we find that
$\Phi_0^{R_0,R_2}[\mathbf{m}]$ in turn extends harmonically
 across $\Gamma_{R_0}$, therefore it is harmonic
in all of $\R^3$. Additionally $\Phi_0^{R_0,R_2}[\mathbf{m}]$ vanishes 
at infinity, hence $\Phi_0^{R_0,R_2}[\mathbf{m}](x)=0$ for all $x\in\R^3$
by Liouville's theorem. Since  $\Phi_0^{R_0,R_2}[\mathbf{m}]=p_{\hat D}$ for $\hat D=\nabla\cdot (\mathbf{m}\otimes \delta_{\SS_{R_0}})$, Lemma \ref{lem:nullpd} now implies that $\mathbf{m}=\mathbf{d}\in\D_{R_0}$ with
$\textnormal{supp}\,\mathbf{d}\subset\Gamma_{R_0}$.
Next, as $\Phi_0^{R_0,R_2}[\mathbf{m}]$ vanishes identically on $\R^3$,
we get from \eqref{eqn:zeropot} that $\Phi_1^{R_1,R_2}[h](x)=0$ for all 
$x\in \R^3\setminus \overline{\BB_{R_1}}$. Then, 
injectivity of the Poisson transform entails that $h\equiv0$,
hence $N(\Phi^{R_1,R_0,R_2})\subset\{(\mathbf{m},0):\mathbf{m}\in \mathcal{D}_{R_0}\,\textrm{supp}(\m)\subset\Gamma_{R_0}\}$.

The reverse inclusion $N(\Phi^{R_1,R_0,R_2})\supset\{(\mathbf{m},0):\mathbf{m}\in \mathcal{D}_{R_0},\,\textrm{supp}(\m)\subset\Gamma_{R_0}\}$ is clear because Lemma \ref{lem:nullpd} yields 
that $\Phi^{R_1,R_0,R_2}[(\mathbf{m},0)](x)=
\Phi_0^{R_0,R_2}[\mathbf{m}](x)=0$, for all $x\in\R^3\setminus\Gamma_{R_0}$
 if $\mathbf{m}\in \mathcal{D}_{R_0}$.
\end{prf}

\begin{cor}\label{cor:unique}
Notation being as in Definition \ref{def:ops} with $\Gamma_{R_0}\not=\SS_{R_0}$, let $\Phi=\Phi^{R_1,R_0,R_2}[(\mathbf{m},h)]$ for some
$\mathbf{m}\in L^2(\Gamma_{R_0},\R^3)$ and some $h\in L^2(\SS_{R_1})$.
Then, a pair of  potentials of the form $\bar{\Phi}_0=\Phi_0^{R_0,R_2}[\bar{\mathbf{m}}]$ and $\bar{\Phi}_1=\Phi_1^{R_1,R_2}[\bar{h}]$, with  $\bar{\mathbf{m}}\in L^2(\Gamma_{R_0},\R^3)$ and $\bar{h}\in L^2(\SS_{R_1})$, is
uniquely determined by the condition $\Phi(x)=\bar{\Phi}_0(x)+\bar{\Phi}_1(x)$, $x\in\SS_{R_2}$.
\end{cor}

\begin{prf}
 From  Theorem \ref{thm:unique} we get that $h$ is uniquely determined
by the values of $\Phi$ on $\SS_{R_2}$, and also that
the components $\mathbf{m}_+\in\Hpz$ and $\mathbf{m}_-\in\Hmz$ of the Hardy-Hodge decomposition of $\mathbf{m}$ are uniquely determined. The former implies $\bar{h}\equiv h$ and the latter $\bar{\mathbf{m}}\equiv\mathbf{m}+\bar{\mathbf{d}}$, for some $\bar{\mathbf{d}}\in\D_{R_0}$. By Lemma \ref{lem:nullpd}  we have that  $\Phi_0^{R_0,R_2}[\mathbf{m}](x)=\Phi_0^{R_0,R_2}[\mathbf{m}+\bar{\mathbf{d}}](x)$ for $x\in\R^3\setminus\SS_{R_0}$, so
we eventually find that $\bar{\Phi}_0$ and $\bar{\Phi}_1$ are uniquely determined.
\end{prf}

Corollary \ref {cor:unique}  answers the uniqueness issue of 
Problem \ref{prob:1} in the positive provided that
$\textnormal{supp}(\mathbf{m})\not=\SS_{R_0}$. In other words, assuming
a locally supported magnetization, it is possible to separate the contribution of the Earth's crust from the contribution of the Earth's core if only the superposition of both magnetic potentials is known on some external orbit $\SS_{R_2}$. Of course, in Geomagnetism, it is 
the magnetic field $\mathbf{B}=\nabla \Phi$ which is measured
rather than the magnetic potential $\Phi$. However,
the result carries over at once to this setting. More in fact is true:
if $\textnormal{supp}(\mathbf{m})\not=\SS_{R_0}$, separation is possible if
only the normal component of  $\mathbf{B}$
is known on $\SS_{R_2}$. Indeed, we have the following theorem.

\begin{thm}\label{thm:unique3}
 Let the setup be as in Definition \ref{def:ops} with $\Gamma_{R_0}\not=\SS_{R_0}$, and consider  the operator 
 \begin{align*}
  \mathbf{B}^{R_1,R_0,R_2}:L^2(\Gamma_{R_0},\R^3)\times L^2(\SS_{R_1},\R^3)\to L^2(\SS_{R_2},\R^3),\quad (\mathbf{m},h)\mapsto \nabla\Phi_0^{R_0,R_2}[\mathbf{m}]+ \nabla\Phi_1^{R_1,R_2}[h].
 \end{align*}
Define further the normal operator:
\begin{align*}
  &\mathbf{B}_\nu^{R_1,R_0,R_2}:L^2(\Gamma_{R_0},\R^3)\times L^2(\SS_{R_1},\R^3)\to L^2(\SS_{R_2}),\quad (\mathbf{m},h)\mapsto \partial_{\nu}\left(\Phi_0^{R_0,R_2}[\mathbf{m}]+ \Phi_1^{R_1,R_2}[h]\right).
 \end{align*}
Then the nullspaces of $\mathbf{B}^{R_1,R_0,R_2}$ and $\mathbf{B}_\nu^{R_1,R_0,R_2}$ are all  given by
 \begin{align*}
  N(\mathbf{B}^{R_1,R_0,R_2})= N(\mathbf{B}_\nu^{R_1,R_0,R_2})
=\{(\mathbf{d},0):\mathbf{d}\in \mathcal{D}_{R_0},\textnormal{  supp}(\mathbf{d})\subset\Gamma_{R_0}\}.
 \end{align*}
\end{thm}

\begin{prf}
Let
$\mathbf{B}_\nu^{R_1,R_0,R_2}[(\mathbf{m},h)](x)=0$ for $x\in\SS_{R_2}$.
Then $\Phi^{R_1,R_0,R_2}[(\mathbf{m},h)]$ has vanishing normal derivative on
$\SS_{R_2}$, and is otherwise harmonic in  
$\R^3\setminus\overline{\BB_{R_2}}$. Note that $\Phi^{R_1,R_0,R_2}[(\mathbf{m},h)]$ is even harmonic across $\SS_{R_2}$ onto a slightly larger open set,
hence there is no issue of smoothness to define derivatives everywhere on
$\SS_{R_2}$. Since $\Phi^{R_1,R_0,R_2}[(\mathbf{m},h)]$ vanishes at 
infinity, its Kelvin transform $u=K_{R_2}[\Phi^{R_1,R_0,R_2}[(\mathbf{m},h)]]$ is harmonic in $\BB_{R_2}$ with $u(0)=0$ \cite[Thm. 4.8]{axler01}, and by
\eqref{diffnK} it holds that $\partial_\nu u(x)+ u(x)/R_2=0$ for
$x\in\SS_{R_2}$. Now, if $u$ is nonconstant and
$x$ is a maximum place for $u$ on $\SS_{R_2}$,
then $\partial_\nu u(x)>0$ by the Hopf lemma 
\cite[Ch. 1, Ex. 25]{axler01}. Hence  $ u(x)<0$, implying that $u<0$ on
$\BB_{R_2}$, which contradicts the maximum principle because $u(0)=0$.
Therefore $u$ vanishes identically and so
does $\Phi^{R_1,R_0,R_2}[(\mathbf{m},h)]$ on $\SS_{R_2}$. 
Appealing to Theorem \ref{thm:unique} now achieves the proof.
\end{prf}

The next corollary follows in the exact same manner as Corollary \ref{cor:unique}. To state it, we indicate  with a subscript $\nu$ the normal component of a field in $L^2(\SS_{R_2},\R^3)$ while a subscript $\tau$ denotes the tangential component.

\begin{cor}\label{cor:unique3}
 Let the setup be as in Definition \ref{def:ops} with $\Gamma_{R_0}\not=\SS_{R_0}$,
and let the operator $\mathbf{B}^{R_1,R_0,R_2}$,
be as in Theorem \ref{thm:unique3}. Define further the operators
 \begin{align*}
  \mathbf{B}_0^{R_0,R_2}:L^2(\Gamma_{R_0},\R^3)\to L^2(\SS_{R_2},\R^3),\quad \mathbf{m}\mapsto\nabla\Phi_0^{R_0,R_2}[\mathbf{m}],
 \end{align*}
and
\begin{align*}
  \mathbf{B}_1^{R_1,R_2}:L^2(\SS_{R_1})\to L^2(\SS_{R_2},\R^3),\quad h\mapsto\nabla\Phi_1^{R_1,R_2}[h].
 \end{align*}
 Let further $\mathbf{B}=\mathbf{B}^{R_1,R_0,R_2}[(\mathbf{m},h)]$, with  $\mathbf{m}\in L^2(\Gamma_{R_0},\R^3)$ and $h\in L^2(\SS_{R_1})$.
A pair of fields of the form $\bar{\mathbf{B}}_0=\mathbf{B}_0^{R_0,R_2}[\bar{\mathbf{m}}]$ and $\bar{\mathbf{B}}_1=\mathbf{B}_1^{R_1,R_2}[\bar{h}]$, with  $\bar{\mathbf{m}}\in L^2(\Gamma_{R_0},\R^3)$ and $\bar{h}\in L^2(\SS_{R_1})$, is uniquely determined by the condition $\mathbf{B}_\nu(x)=(\bar{\mathbf{B}}_0)_\nu(x)+(\bar{\mathbf{B}}_1)_\nu(x)$
and thus,
a fortiori, by the condition $\mathbf{B}(x)=\bar{\mathbf{B}}_0(x)+\bar{\mathbf{B}}_1(x)$ for  $x\in\SS_{R_2}$.
\end{cor}

\begin{rem}
 Opposed to the normal component, it does not suffice to know the tangential component $\B_\tau$ on $\SS_{R_2}$ in order to obtain uniqueness of $\B_0$ and $\B_1$. Namely, letting $\m\equiv0$ and $h$ be any nonzero constant function on $\SS_{R_1}$, then $\B_\tau(x)=(\B_0)_\tau(x)+(\B_1)_\tau(x)=\nabla_{\SS_{R_2}}\Phi_0^{R_0,R_2}[\m](x)+\nabla_{\SS_{R_2}}\Phi_1^{R_1,R_2}[h](x)=0$ and $\B_0(x)=\nabla \Phi_0^{R_0,R_2}[\m](x)=0$ but $\B_1(x)=\nabla\Phi_1^{R_1,R_2}[h](x)=-\frac{hR_1}{|x|^3}x\not=0$ for $x\in\SS_{R_2}$.
\end{rem}

Analogously to the previous considerations, one can separate two potentials 
produced by two magnetizations located on two distinct spheres of radii $R_1<R_0$ (of which the outer magnetization again has to be supported
on a strict subset of $\SS_{R_0}$).
We need only slightly change the setup of Definition \ref{def:ops}:

\begin{defi}\label{def:ops2}
 Let $0<R_1<R_0<R_2$ be fixed radii and $\Gamma_{R_0}\subset\SS_{R_0}$ 
a closed subset. We define
\begin{align*}
  \Psi_0^{R_0,R_2}:L^2(\Gamma_{R_0},\R^3)\to L^2(\SS_{R_2}),\quad \mathbf{m}\mapsto\frac{1}{4\pi}\int_{\Gamma_{R_0}}\m(y)\cdot\frac{x-y}{|x-y|^3}\,\dd \omega_{R_0}(y),\quad x\in\SS_{R_2},
 \end{align*}
and
\begin{align*}
  \Psi_1^{R_1,R_2}:L^2(\SS_{R_1},\R^3)\to L^2(\SS_{R_2}),\quad \mathbf{m}\mapsto\frac{1}{4\pi}\int_{\SS_{R_1}}\m(y)\cdot\frac{x-y}{|x-y|^3}\,\dd \omega_{R_0}(y),\quad x\in\SS_{R_2}.
 \end{align*}
 The superposition of these two operators is denoted by 
 \begin{align*}
  \Psi^{R_1,R_0,R_2}:L^2(\Gamma_{R_0},\R^3)\times L^2(\SS_{R_1},\R^3)\to L^2(\SS_{R_2}),\quad(\mathbf{m},\mathfrak{m})\mapsto  \Psi_0^{R_0,R_2}[\mathbf{m}]+ \Psi_1^{R_1,R_2}[\mathfrak{m}].
 \end{align*}
\end{defi}

\begin{thm}\label{thm:unique2}
 Let the setup be as in Definition \ref{def:ops2} and
assume that $\Gamma_{R_0}\not=\SS_{R_0}$. Then the nullspace of $\Psi^{R_1,R_0,R_2}$ is given by
 \begin{align}
\label{decnomag}
  N(\Psi^{R_1,R_0,R_2})=\{(\mathbf{d},\mathfrak{m}_-+\mathfrak{d}):\mathfrak{m}_-\in\Hmo,\,\mathbf{d}\in \mathcal{D}_{R_0}, \textnormal{ supp}(\mathbf{d})\subset
\Gamma_{R_0},\,\mathfrak{d}\in \mathcal{D}_{R_1}\}.
 \end{align}
\end{thm}

\begin{prf}
 Let $\Psi^{R_1,R_0,R_2}[(\mathbf{m},\mathfrak{m})](x)=0$ for all $x\in\SS_{R_2}$. The same argument as in the proof of Theorem \ref{thm:unique} 
then leads us to $\Psi_0^{R_0,R_2}[\mathbf{m}](x)=0$, $x\in\R^3$, and  $\Psi_1^{R_1,R_2}[\mathfrak{m}](x)=0$, $x\in \R^3\setminus \overline{\BB_{R_1}}$. The
former yields $\mathbf{m}=\mathbf{d}\in\D_{R_0}$, like in Theorem \ref{thm:unique}.  As to the latter, we observe that $\Psi_1^{R_1,R_2}[\mathfrak{m}]=p_{\hat D}$ with $\hat D=\nabla\cdot (\mathfrak{m}\otimes \delta_{\SS_{R_1}})$, 
so Lemma \ref{lem:nullpd} yields that $\mathfrak{m}=\mathfrak{m}_-+\mathfrak{d}$, where $\mathfrak{m}_-\in\Hmo$ and  $\mathfrak{d}\in\D_{R_1}$. Thus, 
the left hand side of \eqref{decnomag} is included in the right hand side.
The reverse inclusion is a direct consequence of Lemma \ref{lem:nullpd}.
\end{prf}

\begin{cor}\label{cor:unique2}
 Let the setup be as in Definition \ref{def:ops2} and 
assume that $\Gamma_{R_0}\not=\SS_{R_0}$. Let further $\Psi=\Psi^{R_1,R_0,R_2}[(\mathbf{m},\mathfrak{m})]$, with  $\mathbf{m}\in L^2(\Gamma_{R_0},\R^3)$ and $\mathfrak{m}\in L^2(\SS_{R_1},\R^3)$. A pair of potentials of the form 
$\bar{\Psi}_0=\Psi_0^{R_0,R_2}[\bar{\mathbf{m}}]$ and $\bar{\Psi}_1=\Psi_1^{R_1,R_2}[\bar{\mathfrak{m}}]$, with  $\bar{\mathbf{m}}\in L^2(\Gamma_{R_0},\R^3)$ and $\bar{\mathfrak{m}}\in L^2(\SS_{R_1},\R^3)$, is
uniquely determined by the condition $\Psi(x)=\bar{\Psi}_0(x)+\bar{\Psi}_1(x)$, $x\in\SS_{R_2}$. 
\end{cor}

\begin{prf}
 From  Theorem \ref{thm:unique} we get that the components $\mathbf{m}_+\in\Hpz$ and $\mathbf{m}_-\in\Hmz$ of the Hardy-Hodge decomposition of $\mathbf{m}$ are uniquely determined by the knowledge of $\Psi$ on $\SS_{R_2}$, while
for $\mathfrak{m}$ only the component $\mathbf{\mathfrak{m}}_+\in\Hpo$ is uniquely determined. The former implies $\bar{\mathbf{m}}\equiv\mathbf{m}+\bar{\mathbf{d}}$, for some $\bar{\mathbf{d}}\in\D_{R_0}$, and the latter yields $\bar{\mathfrak{m}}\equiv\mathfrak{m}+\bar{\mathfrak{m}}_-+\bar{\mathfrak{d}}$, for some $\bar{\mathfrak{m}}_-\in\Hmo$ and $\bar{\mathfrak{d}}\in\D_{R_1}$. Since Lemma \ref{lem:nullpd} yields that $\Psi_0^{R_0,R_2}[\mathbf{m}](x)=\Psi_0^{R_0,R_2}[\mathbf{m}+\mathbf{d}](x)$ and $\Psi_1^{R_1,R_2}[\mathfrak{m}](x)=\Psi_1^{R_1,R_2}[\mathfrak{m}+\bar{\mathfrak{m}}_-+\bar{\mathfrak{d}}](x)$ for $x\in\SS_{R_2}$, we eventually find that $\bar{\Psi}_0$ and $\bar{\Psi}_1$ are uniquely determined.
\end{prf}

\begin{rem}
Analogs of Theorems \ref{thm:unique}, \ref{thm:unique2} and Corollaries
\ref{cor:unique2}, \ref{cor:unique2} are easily seen to hold
for the case of  finitely many
magnetizations $\mathbf{m}_{1},\ldots,\mathbf{m}_n$, and $\mathbf{\mathfrak{m}}$ supported respectively on spheres $\SS_{R_{0,1}},\ldots,\SS_{R_{0,n}}$ 
and $\SS_{R_1}$ of radii $R_1<R_{0,1}<\ldots<R_{0,n}<R_2$,
under the localization assumptions that $\textnormal{supp}(\mathbf{m}_{i})$ 
is a strict subset of  $\SS_{R_{0,i}}$, $i=1,\ldots,n$.
The corresponding separation properties may be of interest when 
investigating the depth profile of crustal magnetizations.
\end{rem}

\subsection{Reconstruction Issues}\label{sec:extremal}

In this section, we discuss how quantities such as the Fourier coefficients $\langle \Phi_0,Y_{n,k}\rangle_{L^2(\SS_{R_2})}$ of $\Phi_0$ can be approximated 
knowing $\Phi$, without having to reconstruct $\Phi_0$ itself. 
Such Fourier coefficients are of interest, e.g., when looking at the power spectra of $\Phi$ and $\Phi_0$ (cf. the empirical way of separating the crustal and the core magnetic fields mentioned in the introduction).
As an extra piece of notation, given $\Gamma_R\subset\SS_R$ and $f:\SS_R\to \R^k$,
we let $f_{|\Gamma_R}:\Gamma_R\to\R^k$ designate the restriction of $f$ to $\Gamma_R$.

\begin{thm}\label{thm:errestfunc}
 Let the setup be as in Definition \ref{def:ops} and assume that $\Gamma_{R_0}\not=\SS_{R_0}$. Then, for every $\eps>0$ and every function $\mathbf{g}\in\Hpz\oplus\Hmz$, 
there exists $f\in L^2(\SS_{R_2})$ (depending on $\eps$ and $\mathbf{g}$) such that
 \begin{align*}
  \left|\langle \Phi^{R_1,R_0,R_2}[\mathbf{m},h],f\rangle_{L^2(\SS_{R_2})}-\langle \mathbf{m},\mathbf{g}_{|\Gamma_{R_0}}\rangle_{L^2(\Gamma_{R_0}, \R^3)}\right|\leq \eps\|(\mathbf{m},h)\|_{L^2(\Gamma_{R_0},\R^3)\times L^2(\SS_{R_1})},
 \end{align*}
 for all $\mathbf{m}\in L^2(\Gamma_{R_0},\R^3)$ and $h\in L^2(\SS_{R_1})$.
\end{thm}

\begin{prf}
 According to Theorem \ref{thm:unique}
and the orthogonality of the Hardy-Hodge decomposition, $(\mathbf{g}_{|\Gamma_{R_0}},0)$ is orthogonal to the nullspace $N(\Phi^{R_1,R_0,R_2})$ of
$\Phi^{R_1,R_0,R_2}$,
for  if $\textnormal{supp}\,\mathbf{d}\subset\Gamma_{R_0}$
then $\langle \mathbf{g}_{|\Gamma_{R_0}},\mathbf{d}\rangle_{L^2(\Gamma_{R_0},\R^3)}=
\langle \mathbf{g},\mathbf{d}\rangle_{L^2(\SS_{R_0},\R^3)}=0$. Therefore,  $(\mathbf{g}_{|\Gamma_{R_0}},0)$ lies in the closure of the range of the adjoint operator $\big(\Phi^{R_1,R_0,R_2}\big)^*$, i.e., to each $\eps>0$ there is  $f\in L^2(\SS_{R_2})$ with
 \begin{align}\label{eqn:adjest}
  \left\|\left(\Phi^{R_1,R_0,R_2}\right)^*[f]-(\mathbf{g}_{|\Gamma_{R_0}},0)\right\|_{L^2(\Gamma_{R_0},\R^3)\times L^2(\SS_{R_1})}\leq \eps.
 \end{align}
 Taking the scalar product with $(\mathbf{m},h)$, we get from \eqref{eqn:adjest} and the Cauchy-Schwarz inequality:
 \begin{align*}
  &\left|\langle \Phi^{R_1,R_0,R_2}[\mathbf{m},h],f\rangle_{L^2(\SS_{R_2})}-\langle \mathbf{m},\mathbf{g}_{|\Gamma_{R_0}}\rangle_{L^2(\Gamma_{R_0}, \R^3)}\right|
  \\&=\left|\left\langle(\mathbf{m},h),\left(\Phi^{R_1,R_0,R_2}\right)^*[f]-(\mathbf{g}_{|\Gamma_{R_0}},0)\right\rangle_{L^2(\Gamma_{R_0}, \R^3)\times L^2(\SS_{R_2})}\right|
  \\&\leq\left\|\left(\Phi^{R_1,R_0,R_2}\right)^*[f]-(\mathbf{g}_{|\Gamma_{R_0}},0)\right\|_{L^2(\Gamma_{R_0},\R^3)\times L^2(\SS_{R_1})}\|(\mathbf{m},h)\|_{L^2(\Gamma_{R_0},\R^3)\times L^2(\SS_{R_1})}
  \\&\leq \eps\|(\mathbf{m},h)\|_{L^2(\Gamma_{R_0},\R^3)\times L^2(\SS_{R_1})},
 \end{align*}
 which is the desired result.
\end{prf}

\begin{cor}\label{cor:approxcoeffs}
 Let the setup be as in Definition \ref{def:ops} with $\Gamma_{R_0}\not=\SS_{R_0}$. Then, for every $\eps>0$ and every function $g\in L^2(\SS_{R_2})$, there exists  $f\in L^2(\SS_{R_2})$ (depending on $\eps$ and $g$) such that
 \begin{align*}
  \left|\langle \Phi^{R_1,R_0,R_2}[\mathbf{m},h],f\rangle_{L^2(\SS_{R_2})}-\langle \Phi_0^{R_0,R_2}[\mathbf{m}],g\rangle_{L^2(\SS_{R_2})}\right|\leq \eps\|(\mathbf{m},h)\|_{L^2(\Gamma_{R_0},\R^3)\times L^2(\SS_{R_1})},
 \end{align*}
 for all $\mathbf{m}\in L^2(\Gamma_{R_0},\R^3)$ and $h\in L^2(\SS_{R_1})$.
\end{cor}

\begin{prf}
 First observe that
\begin{align}\label{eqn:adjrel}
 \left\langle \Phi_0^{R_0,R_2}[\mathbf{m}],g\right\rangle_{L^2(\SS_{R_2})}=\left\langle\mathbf{m},\left( \Phi_0^{R_0,R_2}\right)^*[g]\right\rangle_{L^2(\Gamma_{R_0},\R^3)},
\end{align}
where the adjoint operator of $\Phi_0^{R_0,R_2}$ is given by
\begin{align}
\label{adp}
 \left( \Phi_0^{R_0,R_2}\right)^*:L^2(\SS_{R_2})\to L^2(\Gamma_{R_0},\R^3),\quad g\mapsto \mathbf{H}[g]_{|\Gamma_{R_0}},\nonumber\\
\mathbf{H}[g](x)=-\frac{1}{4\pi}\int_{\SS_{R_2}}g(y)\frac{x-y}{|x-y|^3}\dd \omega_{R_2}(y),\quad x\in
\SS_{R_0}.
\end{align}
Clearly $\mathbf{H}[g]\in \H^2_{+,R_0}$ whenever $g\in L^2(\SS_{R_2})$,
therefore, \eqref{eqn:adjrel} together with Theorem \ref{thm:errestfunc} yield the desired result. 
\end{prf}

\begin{rem}\label{rem:empsep}
The interest of Corollary \ref{cor:approxcoeffs} from the Geophysical viewpoint
lies with the fact that 
$\Phi^{R_1,R_0,R_2}[\mathbf{m},h]$ (more specifically: its gradient) 
corresponds to the measurements on $\SS_{R_2}$ of the superposition of
the core and crustal contributions, whereas $\Phi_0^{R_0,R_2}[\mathbf{m}]$
corresponds to the crustal contribution alone. Thus, if we can compute 
$f$ knowing $g$, we shall in principle be able to get information 
on the crustal contribution up to arbitrary small error.
Note also  that $(\mathbf{g},0)\not\in\textnormal{Ran}\,\bigl(\Phi^{R_1,R_0,R_2}\big)^*$ unless  $\mathbf{g}\equiv0$,
due to the injectivity of the 
adjoint of the Poisson transform (which is again a Poisson transform).
Therefore we can only hope for an approximation of $\langle \Phi_0^{R_0,R_2}[\mathbf{m}],g\rangle_{L^2(\SS_{R_2})}$ in
Corollary \ref{cor:approxcoeffs}, up to a relative error of $\eps>0$,
but not for an exact reconstruction.
\end{rem}

Results analogous to Theorem \ref{thm:errestfunc} and Corollary \ref{cor:approxcoeffs} mechanically hold in the setup of Theorem \ref{thm:unique3} and 
Corollary \ref{cor:unique3} (i.e., separation of the crustal and core magnetic fields $\mathbf{B}_0$ and $\mathbf{B}_1$ instead of the potentials) and  in
the setup of Theorem \ref{thm:unique2} and Corollary \ref{cor:unique2} 
(i.e., separation of the potentials $\Psi_0$ and $\Psi_1$ due to 
magnetizations on $\SS_{R_0}$ and $\SS_{R_1}$). 
Below we state the corresponding results
but we  
omit the proofs for they are similar 
to the previous ones.

\begin{thm}\label{thm:errestfunc3}
 Let the setup be as in 
Theorem \ref{thm:unique3}. Then, for every $\eps>0$ and every field $\mathbf{g}\in\Hpz\oplus\Hmz$, 
there exists  $\mathbf{f}\in L^2(\SS_{R_2},\R^3)$ (depending on $\eps$ and $\mathbf{g}$) such that
 \begin{align*}
  \left|\langle \mathbf{B}^{R_1,R_0,R_2}[\mathbf{m},h],\mathbf{f}\rangle_{L^2(\SS_{R_2},\R^3)}-\langle \mathbf{m},\mathbf{g}_{|\Gamma_{R_0}}\rangle_{L^2(\Gamma_{R_0}, \R^3)}\right|\leq \eps\|(\mathbf{m},h)\|_{L^2(\Gamma_{R_0},\R^3)\times L^2(\SS_{R_1})},
 \end{align*}
 for all $\mathbf{m}\in L^2(\Gamma_{R_0},\R^3)$ and $h\in L^2(\SS_{R_1})$.
The same holds if $\mathbf{B}^{R_1,R_0,R_2}[\mathbf{m},h]$ gets replaced by
$\mathbf{B}_\nu^{R_1,R_0,R_2}[\mathbf{m},h]$, this time with $\mathbf{f}\in L^2(\SS_{R_2})$.
\end{thm}

\begin{cor}\label{cor:approxcoeffs3}
 Let the setup be as in 
Theorem \ref{thm:unique3} and Corollary \ref{cor:unique3}. Then, for every $\eps>0$ and every field $\mathbf{g}\in L^2(\SS_{R_2},\R^3)$, there exists  $\mathbf{f}\in L^2(\SS_{R_2},\R^3)$ (depending on $\eps$ and $\mathbf{g}$) such that
 \begin{align*}
  \left|\langle \mathbf{B}^{R_1,R_0,R_2}[\mathbf{m},h],\mathbf{f}\rangle_{L^2(\SS_{R_2},\R^3)}-\langle \mathbf{B}_0^{R_0,R_2}[\mathbf{m}],\mathbf{g}\rangle_{L^2(\SS_{R_2},\R^3)}\right|\leq \eps\|(\mathbf{m},h)\|_{L^2(\Gamma_{R_0},\R^3)\times L^2(\SS_{R_1})},
 \end{align*}
 for all $\mathbf{m}\in L^2(\Gamma_{R_0},\R^3)$ and $h\in L^2(\SS_{R_1})$.
The same holds if $\mathbf{B}^{R_1,R_0,R_2}[\mathbf{m},h]$ gets replaced by
$\mathbf{B}_\nu^{R_1,R_0,R_2}[\mathbf{m},h]$, this time with $\mathbf{f}\in L^2(\SS_{R_2})$.
\end{cor}

\begin{thm}\label{thm:errestfunc2}
 Let the setup be as in Definition \ref{def:ops2} with $\Gamma_{R_0}\not=\SS_{R_0}$. Then, for every $\eps>0$ and every $\mathbf{g}\in\Hpz\oplus\Hmz$, 
there exists  $f\in L^2(\SS_{R_2})$ (depending on $\eps$ and $\mathbf{g}$) such that
 \begin{align*}
  \left|\langle \Psi^{R_1,R_0,R_2}[\mathbf{m},\mathfrak{m}],f\rangle_{L^2(\SS_{R_2})}-\langle \mathbf{m},\mathbf{g}_{|\Gamma_{R_0}}\rangle_{L^2(\Gamma_{R_0}, \R^3)}\right|\leq \eps\|(\mathbf{m},\mathfrak{m})\|_{L^2(\Gamma_{R_0},\R^3)\times L^2(\SS_{R_1},\R^3)},
 \end{align*}
 for all $\mathbf{m}\in L^2(\Gamma_{R_0},\R^3)$ and $\mathfrak{m}\in L^2(\SS_{R_1},\R^3)$.
\end{thm}

\begin{cor}\label{cor:approxcoeffs2}
 Let the setup be as in Definition \ref{def:ops2} with $\Gamma_{R_0}\not=\SS_{R_0}$. Then, for every $\eps>0$ and every function $g\in L^2(\SS_{R_2})$, there exists  $f\in L^2(\SS_{R_2})$ (depending on $\eps$ and $g$) such that
 \begin{align*}
  \left|\langle \Psi^{R_1,R_0,R_2}[\mathbf{m},\mathfrak{m}],f\rangle_{L^2(\SS_{R_2})}-\langle \Psi_0^{R_0,R_2}[\mathbf{m}],g\rangle_{L^2(\SS_{R_2})}\right|\leq \eps\|(\mathbf{m},\mathfrak{m})\|_{L^2(\Gamma_{R_0},\R^3)\times L^2(\SS_{R_1},\R^3)},
 \end{align*}
 for all $\mathbf{m}\in L^2(\Gamma_{R_0},\R^3)$ and $\mathfrak{m}\in L^2(\SS_{R_1},\R^3)$.
\end{cor}

\subsection{The Case $\Gamma_{R_0}=\SS_{R_0}$}\label{sec:gammas}

We turn to the case where $\Gamma_{R_0}=\SS_{R_0}$. Then, uniqueness 
no longer  holds in Problem \ref{prob:1}, but one can 
obtain the singular value decomposition of $\Phi^{R_1,R_0,R_2}$ 
fairly explicitly and thereby quantify non-uniqueness. 
Indeed basic computations using spherical harmonics yield:
\begin{align}
 &(\Phi_0^{R_0,R_2})^*[Y_{n,k}](x)\nonumber
 \\&=\frac{1}{4\pi}\int_{\SS_{R_2}}Y_{n,k}\left(\frac{y}{|y|}\right)\nabla_x\frac{1}{|x-y|}\dd \omega_{R_2}(y)\nonumber
 \\&=\frac{1}{4\pi}\sum_{m=0}^\infty\nabla_x\int_{\SS_{R_2}}\frac{1}{|y|}\left(\frac{|x|}{|y|}\right)^{m}Y_{n,k}\left(\frac{y}{|y|}\right)P_{m}\left(\frac{x}{|x|}\cdot\frac{y}{|y|}\right)\dd \omega_{R_2}(y)\nonumber
  \\&=\frac{1}{4\pi}\sum_{m=0}^\infty\sum_{l=1}^{2m+1}\frac{4\pi}{2m+1}\frac{1}{R_2^{m+1}}\nabla_x\left(|x|^mY_{m,l}\left(\frac{x}{|x|}\right)\right)\int_{\SS_{R_2}}Y_{n,k}\left(\frac{y}{|y|}\right)Y_{m,l}\left(\frac{y}{|y|}\right)\dd \omega_{R_2}(y)\nonumber
 \\&=\frac{R_2}{2n+1}\nabla H_{n,k}^{R_2}(x)=\frac{R_2}{2n+1}\left(\frac{R_0}{R_2}\right)^n\nabla H_{n,k}^{R_0}(x),\quad x\in\SS_{R_0},\label{eqn:phi0coeff}
 \end{align}
 and
 \begin{align}
&(\Phi_1^{R_1,R_2})^*[Y_{n,k}](x)\nonumber
\\&=\frac{1}{4\pi R_1}\int_{\SS_{R_2}}Y_{n,k}\left(\frac{y}{|y|}\right)\frac{|y|^2-R_1^2}{|x-y|^3}\dd \omega_{R_2}(y)\nonumber
 \\&=\frac{1}{4\pi R_1}\sum_{m=0}^\infty(2m+1)\int_{\SS_{R_2}}\frac{1}{|y|}\left(\frac{|x|}{|y|}\right)^mP_m\left(\frac{x}{|x|}\cdot \frac{y}{|y|}\right)Y_{n,k}\left(\frac{y}{|y|}\right)\dd \omega_{R_2}(y)\nonumber
  \\&=\frac{1}{R_1R_2}\sum_{m=0}^\infty\sum_{l=1}^{2m+1}\left(\frac{R_1}{R_2}\right)^mY_{m,l}\left(\frac{x}{|x|}\right)\int_{\SS_{R_2}}Y_{n,k}\left(\frac{y}{|y|}\right)Y_{m,l}\left(\frac{y}{|y|}\right)\dd \omega_{R_2}(y)\nonumber
  \\&=\left(\frac{R_1}{R_2}\right)^{n-1}Y_{n,k}\left(\frac{x}{|x|}\right),\quad x\in\SS_{R_1},\label{eqn:phi1coeff}
 \end{align}
 where $H_{n,k}^{R_2}$, $H_{n,k}^{R_0}$ are the inner harmonic from Section \ref{sec:aux} and $P_m$ the Legendre polynomial of degree $m$ (see, e.g., \cite[Ch. 3]{freeden98,freedenschreiner09} for details). So,  we get for the adjoint operator $(\Phi^{R_1,R_0,R_2})^*$ that
 \begin{align}\label{eqn:adjynk}
  (\Phi^{R_1,R_0,R_2})^*[Y_{n,k}]=\left(\frac{R_2}{2n+1}\left(\frac{R_0}{R_2}\right)^n\nabla H_{n,k}^{R_0}\,,\left(\frac{R_1}{R_2}\right)^{n-1}Y_{n,k}\right)^T.
 \end{align}
  Similar calculations also yield that
\begin{align*}
 \Phi_0^{R_0,R_2}[\nabla H_{n,k}^{R_0}](x)= \frac{n}{R_2}\left(\frac{R_0}{R_2}\right)^n Y_{n,k}\left(\frac{x}{|x|}\right),\quad x\in \SS_{R_2},
 \end{align*}
 and
 \begin{align*}
 \Phi_1^{R_1,R_2}[Y_{n,k}](x)=\left(\frac{R_1}{R_2}\right)^{n+1}Y_{n,k}\left(\frac{x}{|x|}\right),\quad x\in \SS_{R_2},
 \end{align*}
 so we obtain for $\Phi^{R_1,R_0,R_2}$ that
 \begin{align}
  &\Phi^{R_1,R_0,R_2}[\alpha\nabla H_{n,k}^{R_0},\beta Y_{m,l}]=\alpha \frac{n}{R_2}\left(\frac{R_0}{R_2}\right)^n Y_{n,k}+\beta \left(\frac{R_1}{R_2}\right)^{m+1} Y_{m,l},\label{eqn:phiynk}
 \end{align}
 with $\alpha,\beta\in\R$. 
 Based on the representations \eqref{eqn:adjynk} and \eqref{eqn:phiynk}, 
further computation leads us  to a characterization of the nullspace of $\Phi^{R_1,R_0,R_2}$ in Lemma \ref{lem:gammaes}. Note that $\Phi^{R_1,R_0,R_2}:L^2(\Gamma_{R_0},\R^3)\times L^2(\SS_{R_1})\to L^2(\SS_{R_2})$ is a compact operator, being the sum of two compact operators (for $\Phi_0^{R_0,R_2}$ and
$\Phi_1^{R_1,R_2}$ have continuous kernels).

 \begin{lem}\label{lem:gammaes}
  Let $\Gamma_{R_0}=\SS_{R_0}$, then the nullspace of $\Phi^{R_1,R_0,R_2}$ is 
given by 
  \begin{align*}
   N(\Phi^{R_1,R_0,R_2})=&\{(\mathbf{m}_-+\mathbf{d},0):\mathbf{m}_-\in\Hmz,\,\mathbf{d}\in \mathcal{D}_{R_0}\}
   \\&\cup\,\overline{\textnormal{span}\left\{\left(\nabla H_{n,k}^{R_0}\,,-\frac{n}{R_1}\left(\frac{R_0}{R_1}\right)^nY_{n,k}\right)^T:n\in\mathbb{N},k=1,\ldots,2n+1\right\}},
  \end{align*}
  while the orthogonal complement reads 
  \begin{align*}
    N(\Phi^{R_1,R_0,R_2})^\perp=\overline{\textnormal{span}\left\{\left(\nabla H_{n,k}^{R_0}\,,\frac{2n+1}{R_1}\left(\frac{R_1}{R_0}\right)^nY_{n,k}\right)^T:n\in\mathbb{N},k=1,\ldots,2n+1\right\}}.
  \end{align*}
  All non-zero eigenvalues values of $(\Phi^{R_1,R_0,R_2})^*\Phi^{R_1,R_0,R_2}$ are of the form
  \begin{align*}
   \sigma_{n}=\frac{n}{2n+1}\left(\frac{R_0}{R_2}\right)^{2n}+\left(\frac{R_1}{R_2}\right)^{2n},\quad n\in\mathbb{N},
  \end{align*}
  and the corresponding eigenvectors in $L^2(\SS_{R_0},\R^3)\times L^2(\SS_{R_1})$ are
  \begin{align*}
   \left(\nabla H_{n,k}^{R_0}\,,\frac{2n+1}{R_1}\left(\frac{R_1}{R_0}\right)^nY_{n,k}\right)^T,\quad n\in\mathbb{N},\,k=1,\ldots,2n+1.
  \end{align*}
\label{NG}
  \end{lem}
  
Lemma \ref{NG} entails that the nullspace of  $\Phi^{R_1,R_0,R_2}$ contains
elements of the form $(\mathbf{m},h)$ with $h\not=0$, hence 
$\Phi^{R_1,R_0,R_2}[(\mathbf{m},h)]$ may well vanish on
$\SS_{R_2}$ even though $\Phi_1^{R_1,R_2}[h]$ is nonzero there, by injectivity of the Poisson representation.
In other words, separation of the potentials $\Phi_0^{R_0,R_2}$ and
$\Phi_1^{R_1,R_2}$ knowing their sum on $\SS_{R_2}$ 
is no longer possible in general
if $\Gamma_{R_0}=\SS_{R_0}$. 


\section{Extremal Problems and Numerical Examples}\label{sec:num}

In this section, we provide some first approaches on how the results from the previous sections can be used to approximate the Fourier coefficients of $\Phi_0$ (cf. Section \ref{sec:fouriernum}), as well as $\Phi_0$ itself via the  reconstruction of $\mathbf{m}$ and $h$ (cf. Section \ref{sec:phi0num}). For brevity, we treat only separation of the crustal and core magnetic potentials
 (underlying operator $\Phi^{R_1,R_0,R_2}$) and not the separation of the crustal and core magnetic fields  (underlying operator $\mathbf{B}^{R_1,R_0,R_2}$) 
nor the separation of potentials generated by two magnetizations on different spheres (underlying operator $\Psi^{R_1,R_0,R_2}$). The procedure
in such cases is of course similar.

\subsection{Reconstruction of Fourier Coefficients of $\Phi_0$}\label{sec:fouriernum}

To get a feeling of how functions $f$ in 
Corollary \ref{cor:approxcoeffs} 
behave, let us derive some of
their basic properties.  Recall they where isentified to be those
$f\in L^2(\SS_{R_2})$ satisfying \eqref{eqn:adjest} with 
$\mathbf{g}=(\Phi_0^{R_0,R_2})^*[g]$.

 \begin{lem}\label{lem:feps}
  Let $0\not=g\in L^2(\SS_{R_2})$ and set $\mathbf{g}=(\Phi_0^{R_0,R_2})^*[g]$. To each $\eps>0$, let $f_\eps\in L^2(\SS_{R_2})$ satisfy $\|(\Phi^{R_1,R_0,R_2})^*[f_\eps]-(\mathbf{g}_{|\Gamma_{R_0}},0)\|_{L^2(\Gamma_{R_0},\R^3)\times L^2(\SS_{R_1})}\leq \eps$. Then:
  \begin{itemize}
   \item[(a)] $\lim\limits_{\eps\to0}\|f_\eps\|_{L^2(\SS_{R_2})}=\infty$,
   \item[(b)]$\lim\limits_{\eps\to0}\|(\Phi_1^{R_1,R_2})^*[f_\eps]\|_{L^2(\SS_{R_1})}=0$, 
   \item[(c)] $\lim\limits_{\eps\to0}\langle f_\eps,Y_{n,k}\rangle_{L^2(\SS_{R_2})}=0$, for fixed $n\in\mathbb{N}_0$, $k=1,\ldots,n$. 
  \end{itemize}
 \end{lem}

 \begin{prf}
  From the considerations in Remark \ref{rem:empsep} we know that $(\mathbf{g}_{|\Gamma_{R_0}},0)\in\overline{\textnormal{Ran}\big(\big(\Phi^{R_1,R_0,R_2}\big)^*\big)}$ but $(\mathbf{g}_{|\Gamma_{R_0}},0)\not\in\textnormal{Ran}\big(\big(\Phi^{R_1,R_0,R_2}\big)^*\big)$. Thus, $\|f_\eps\|_{L^2(\SS_{R_2})}$ cannot 
remain bounded as $\eps\to0$, otherwise a weak limit point $f_0\in L^2(\SS_{R_2})$
would meet $(\Phi^{R_1,R_0,R_2})^*[f_0]=(\mathbf{g}_{|\Gamma_{R_0}},0)$,
a contradiction which proves (a). Next, the relation
  \begin{align*}
  &\|(\Phi^{R_1,R_0,R_2})^*[f_\eps]-(\mathbf{g}_{|\Gamma_{R_0}},0)\|_{L^2(\Gamma_{R_0},\R^3)\times L^2(\SS_{R_1})}^2
  \\&=\|(\Phi_0^{R_0,R_2})^*[f_\eps]-\mathbf{g}_{|\Gamma_{R_0}}\|_{L^2(\Gamma_{R_0},\R^3)}^2+\|(\Phi_1^{R_1,R_2})^*[f_\eps]\|_{L^2(\SS_{R_1})}^2\leq \eps^2
  \end{align*}
immediately implies that $\lim_{\eps\to0}\|(\Phi_1^{R_1,R_2})^*[f_\eps]\|_{L^2(\SS_{R_1})}=0$ which is (b). Finally,
expanding $f_\eps$ in spherical harmonics, one readily verifies that
\eqref{eqn:phi1coeff} together with (b) yields part (c). 
 \end{prf}
 
Next, we give a quantitative appraisal of the
fact that the Fourier coefficients of $\Phi_0^{R_0,R_2}$ 
on $\SS_{R_2}$, to be estimated up to relative precision 
$\eps$ by choosing $g= Y_{p,q}$ in
Corollary \ref{cor:approxcoeffs}, can be approximated directly by those of 
 $\Phi^{R_1,R_0,R_2}$ (i.e.,  neglecting entirely the core contribution)
when $\frac{R_1}{R_2}$ is small enough (i.e.,  the core is far from 
the measurement orbit) and the degree $p$ is large enough.  We also 
give a quantitative version of Lemma \ref{lem:feps} point $(c)$.
This
provides us with bounds  on the validity of the separation technique
consisting merely of a sharp cutoff in the frequency domain.

 \begin{lem}\label{lem:fpq}
  Let $\eps>0$ and choose $\mathbf{g}=(\Phi_0^{R_0,R_2})^*[Y_{p,q}]$ for some $p\in\mathbb{N}_0$ and $q\in\{1,\ldots 2p+1\}$. The the following assertions hold true.
  \begin{itemize}
  \item[(a)]  If $R_1^2\big(\frac{R_1}{R_2}\big)^{p-1}\leq\eps$, then $f=Y_{p,q}$ satisfies 
\begin{equation}
\label{estims}
\|(\Phi^{R_1,R_0,R_2})^*[f]-(\mathbf{g}_{|\Gamma_{R_0}},0)\|_{L^2(\Gamma_{R_0},\R^3)\times L^2(\SS_{R_1})}\leq \eps.
\end{equation}
  \item[(b)] If $f\in  L^2(\SS_{R_2})$ satisfies 
$\|(\Phi^{R_1,R_0,R_2})^*[f]-(\mathbf{g}_{|\Gamma_{R_0}},0)\|_{L^2(\Gamma_{R_0},\R^3)\times L^2(\SS_{R_1})}\leq \eps$,
then, for all $n\in\mathbb{N}_0$, $k=1,\ldots, 2n+1$,
\begin{equation}
\label{qlfepsc}
|\langle f,Y_{n,k}\rangle_{L^2(\SS_{R_2})}|\leq\eps
\frac{R_2^{n-1}}{R_1^{n+1}}.
\end{equation} 
  \end{itemize}
\end{lem}

\begin{prf}
  To prove $(a)$, note that
$(\Phi^{R_1,R_0,R_2})^*=
\left((\Phi_0^{R_0,R_2})^*,(\Phi_1^{R_1,R_2})^*\right)$ and
by \eqref{eqn:phi1coeff} that
\[  \|(\Phi_1^{R_1,R_2})^*[f]\|_{L^2(\SS_{R_1})}=R_1^2\left(\frac{R_1}{R_2}\right)^{p-1}\leq \eps,
\]
while $\|(\Phi_0^{R_0,R_2})^*[f]-\mathbf{g}_{|\Gamma_{R_0}}\|_{L^2(\Gamma_{R_0},\R^3)}=0$ 
if $f=Y_{p,q}$. Hence \eqref{estims} holds.

As to $(b)$, any $f\in L^2(\SS_{R_2})$ satisfying $\|(\Phi^{R_1,R_0,R_2})^*[f]-(\mathbf{g},0)\|_{L^2(\Gamma_{R_0},\R^3)\times L^2(\SS_{R_1})}\leq \eps$ satisfies in particular, in view of \eqref{eqn:phi1coeff}:
 \begin{align*}
  \|(\Phi_1^{R_1,R_2})^*[f]\|_{L^2(\SS_{R_1})}^2=\sum_{n=0}^\infty\sum_{k=1}^{2n+1}R_1^4\left(\frac{R_1}{R_2}\right)^{2(n-1)}|\langle f,Y_{n,k}\rangle_{L^2(\SS_{R_2})}|^2\leq \eps^2,
 \end{align*}
from which \eqref{qlfepsc} follows at once.
\end{prf}

We turn to the computation of
a function $f$ as in Corollary \ref{cor:approxcoeffs}, regardless of 
assumptions on $\frac{R_1}{R_2}$ or on the degree of a spherical harmonics $Y_{n,k}$ for which we want to estinate $\langle\Phi_0^{R_0,R_2},Y_{n,k}\rangle_{L^2(\SS_{R_2})}$.  One way  is 
to solve the following extremal problem. Note that finding $f$ requires no
data on the potential $\Phi$ that we eventually want to separate into $\Phi_0+\Phi_1$.   

\begin{prob}\label{prob:2}
Let the setup be as in Definition \ref{def:ops} with 
$\Gamma_{R_0}\not=\SS_{R_0}$. Fix $g\in L^2(\SS_{R_2})$ as well as  $\eps>0$,
and  set $\mathbf{g}=(\Phi_0^{R_0,R_2})^*[g]$. Then, find $f\in W^{1,2}(\SS_{R_2})$ such that
\begin{align}\label{eqn:prob2}
 \|f\|_{W^{1,2}(\SS_{R_2})}=\inf_{\atopp{\bar{f}\in W^{1,2}(\SS_{R_2}),}{\|(\Phi^{R_1,R_0,R_2})^*[\bar{f}]-(\mathbf{g}_{|\Gamma_{R_0}},0)\|_{L^2(\Gamma_{R_0},\R^3)\times L^2(\SS_{R_1})}\leq \eps }} \|\bar{f}\|_{W^{1,2}(\SS_{R_2})}.
\end{align}
\end{prob}
It may look strange to seek $f\in W^{1,2}(\SS_{R_2})$ whereas 
Corollary \ref{cor:approxcoeffs} merely deals with scalar products
in $L^2(\SS_{R_2})$. This extra-smoothness 
requirement, though, helps regularizing the 
problem.

\begin{lem}
\label{exun}
 Let the setup be as in  Problem \ref{prob:2} and $g\in L^2(\SS_{R_2})$  
with $\|\mathbf{g}_{|\Gamma_{R_0}}\|_{L^2(\Gamma_{R_0},\R^3)}>\eps$. Then, 
there exists a unique solution $0\not\equiv f\in W^{1,2}(\SS_{R_2})$ to Problem  \ref{prob:2}. Moreover, the constraint in \eqref{eqn:prob2} is saturated, i.e. $\|(\Phi^{R_1,R_0,R_2})^*[f]-(\mathbf{g}_{|\Gamma_{R_0}},0)\|_{L^2(\Gamma_{R_0},\R^3)\times L^2(\SS_{R_1})}= \eps$. 
\end{lem}

\begin{prf}
 Since $\mathbf{H}[g]$ given by \eqref{adp} 
lies in $\Hpz$, the same argument as in the proof of 
Theorem \ref{thm:errestfunc} and the density of 
$W^{1,2}(\SS_{R_2})$ in $L^2(\SS_{R_2})$ together imply 
the existence of $\bar{f}\in W^{1,2}(\SS_{R_2})$ such that 
$\|(\Phi^{R_1,R_0,R_2})^*[\bar{f}]-(\mathbf{g}_{|\Gamma_{R_0}},0)\|_{L^2(\Gamma_{R_0},\R^3)\times L^2(\SS_{R_1})}\leq\eps$ is satisfied,  which ensures
that the closed convex subset of $W^{1,2}(\SS_{R_2})$ defined by
 \begin{align*}
  \mathcal{C}_\eps=\left\{\bar{f}\in W^{1,2}(\SS_{R_2}):\|(\Phi^{R_1,R_0,R_2})^*[\bar{f}]-(\mathbf{g}_{|\Gamma_{R_0}},0)\|_{L^2(\Gamma_{R_0},\R^3)\times L^2(\SS_{R_1})}\leq \eps\right\}
 \end{align*}
 is non-empty. Existence and uniqueness of a minimizer $f$ now follows from 
that of a projection of minimum norm on any nonempty convex set in a 
Hilbert space. 
 From the assumption that $\|\mathbf{g}\|_{L^2(\Gamma_{R_0},\R^3)}>\eps$, we 
get that $f\not\equiv0$ because $0\notin\mathcal{C}_\eps$.
If the constraint is not saturated,
then there is $\delta>0$ such that, for every $\bar{f}\in  W^{1,2}(\SS_{R_2})$ with $\|\bar{f}\|_{W^{1,2}(\SS_{R_2})}\leq 1$, also $f+t\bar{f}$  satisfies the constraint $\|(\Phi^{R_1,R_0,R_2})^*[f+t\bar{f}]-(\mathbf{g},0)\|_{L^2(\Gamma_{R_0},\R^3)\times L^2(\SS_{R_1})}\leq \eps$ for
$t\in(-\delta,\delta)$. Since $f$ is a minimizer, this implies 
 \begin{align*}
  0&=\partial_t \|f+t\bar{f}\|_{ W^{1,2}(\SS_{R_2})}^2\Big|_{t=0}=2\left\langle f,\bar{f}\right\rangle_{ W^{1,2}(\SS_{R_2})},
 \end{align*}
 for every $\bar{f}\in  W^{1,2}(\SS_{R_2})$ with $\|\bar{f}\|_{ W^{1,2}(\SS_{R_2})}\leq 1$. Thus $f\equiv 0$, contradicting what precedes.
\end{prf}

 \begin{rem}\label{rem:shdisc}
  Lemma \ref{lem:feps} and the exponential decay of the eigenvalues of $(\Phi_1^{R_1,R_2})^*$ in \eqref{eqn:phi1coeff} suggest that most
of the relevant information regarding a solution $f\in W^{1,2}(\SS_{R_2})$ 
of Problem \ref{prob:2} must be contained in Fourier coefficients $\langle f,Y_{n,k}\rangle_{L^2(\SS_{R_2})}$ of increasingly high degrees $n$ as $\eps\to0$.
  \color{black} Lemma \ref{lem:fpq} provides a hint at the range of 
accuracies $\eps$ for which numerical solutions of Problem \ref{prob:2} 
with $\mathbf{g}=(\Phi_0^{R_0,R_2})^*[Y_{p,q}]$
behave differently for  small and  large $p$.
 \end{rem}

\subsubsection*{Discretization}

For the actual solution of Problem \ref{prob:2}, we assume that
$\|\mathbf{g}_{|\Gamma_{R_0}}\|_{L^2(\Gamma_{R_0},\R^3)}>\eps$,
 hence the constraint is saturated by Lemma \ref{exun}, and
 we use a Lagrangian formulation and obtain from \cite[Thm. 2.1]{chalendar03} that $f\in W^{1,2}(\SS_{R_2})$ 
solves for
 \begin{align}
  \Big(\id+\lambda\,\big(\Phi^{R_1,R_0,R_2}\big)^{**}\,\big(\Phi^{R_1,R_0,R_2}\big)^*\Big)[f]=\lambda\, \big(\Phi^{R_1,R_0,R_2}\big)^{**}[(\mathbf{g}_{|\Gamma_{R_0}},0)],\label{eqn:lag}
 \end{align}
where $\lambda>0$ is such that $\|(\Phi^{R_1,R_0,R_2})^*[f]-(\mathbf{g}_{|\Gamma_{R_0}},0)\|_{L^2(\Gamma_{R_0},\R^3)\times L^2(\SS_{R_1})}= \eps$. Here, 
the operator $\big(\Phi^{R_1,R_0,R_2}\big)^{**}$ stands for
the adjoint of the restriction of $\big(\Phi^{R_1,R_0,R_2}\big)^{*}$ to 
the domain $W^{1,2}(\SS_{R_2})$. In order to avoid computing $\big(\Phi^{R_1,R_0,R_2}\big)^{**}$, we rewrite \eqref{eqn:lag} in variational form: to
\begin{align}
 &\left\langle f,\varphi\right\rangle_{W^{1,2}(\SS_{R_2})}+\lambda \left\langle \big(\Phi^{R_1,R_0,R_2}\big)^*[f],\big(\Phi^{R_1,R_0,R_2}\big)^*[\varphi]\right\rangle_{L^2(\Gamma_{R_0},\R^3)\times L^2(\SS_{R_1})}\nonumber
 \\&=\lambda\left\langle (\mathbf{g}_{|\Gamma_{R_0}},0), \big(\Phi^{R_1,R_0,R_2}\big)^*[\varphi]\right\rangle_{L^2(\Gamma_{R_0},\R^3)\times L^2(\SS_{R_1})},\label{eqn:varform}
 \end{align}
for all $\varphi\in W^{1,2}(\SS_{R_2})$. Remark \ref{rem:shdisc} indicates that a discretization of $f$ in terms of finitely many spherical harmonics is generally not advisable. As a remedy, we use a discretization in terms of the Abel-Poisson kernels
\begin{align}
 K_\gamma(t)=\frac{1}{4\pi}\frac{1-\gamma^2}{(1+\gamma^2-2\gamma t)^{\frac{3}{2}}},\quad t\in[-1,1].\label{eqn:APkernel}
\end{align}
More precisely, we expand $f$ as 
 \begin{align}
 f(x)&=\sum_{m=1}^M\alpha_{m} K_{\gamma,m}(x) =\sum_{m=1}^M\alpha_{m}\sum_{n=0}^\infty\sum_{k=1}^{2n+1}\gamma^nY_{n,k}\left(\frac{x}{|x|}\right)Y_{n,k}(x_m),\quad x\in\SS_{R_2},\label{eqn:discf}
\end{align}
where $K_{\gamma,m}(x)=K_\gamma(\frac{x}{|x|}\cdot x_m)$. The parameter $\gamma\in(0,1)$ is fixed and controls the spatial localization of $K_{\gamma,m}$ \color{black} (a parameter $\gamma$ close to one means a strong localization) \color{black} while $x_{m}\in\SS_1$, $m=1,\ldots M,$ denote the spatial centers of the kernels $K_{\gamma,m}$. Furthermore, one can see from \eqref{eqn:discf} that $\gamma$  relates to the influence of higher spherical harmonic degrees  in the discretization of $f$. Some general properties of the Abel-Poisson kernel $K_\gamma$ can be found, e.g., in \cite[Ch. 5]{freeden98}. Computations based on the representations in Section \ref{sec:gammas} yield
\begin{align}
  &(\Phi^{R_1,R_0,R_2})^*[K_{\gamma,m}]\nonumber
  \\&=\sum_{n=0}^{\infty}\sum_{k=1}^{2n+1}Y_{n,k}(x_m)\gamma^n\left(\frac{R_2}{2p+1}\left(\frac{R_0}{R_2}\right)^n\nabla H_{n,k}^{R_0}\,,\left(\frac{R_1}{R_2}\right)^{n-1}Y_{n,k}\right)^T\nonumber
  \\&=\left(\nabla\sum_{n=0}^{\infty}\sum_{k=1}^{2n+1}\gamma^n\frac{R_2}{2p+1}\left(\frac{R_0}{R_2}\right)^n \left(\frac{|\cdot|}{R_0}\right)^nY_{n,k}(x_m)Y_{n,k}\left(\frac{\cdot}{|\cdot|}\right)\,,\left(\frac{R_2}{R_1}\right)K_{\frac{\gamma R_1}{R_2},m}\right)^T\nonumber
 \\&=\left(\frac{R_2}{4\pi}\nabla F_{\frac{\gamma|\cdot|}{R_2},m}\,,\left(\frac{R_2}{R_1}\right)K_{\frac{\gamma R_1}{R_2},m}\right)^T,\label{eqn:phikm}
 \end{align}
 where $F_{\gamma,m}(x)=F_{\gamma}(\frac{x}{|x|}\cdot x_m)$, with $F_\gamma(t)=(1+\gamma^2-2\gamma t)^{-\frac{1}{2}}$ for $t\in[-1,1]$. Inserting \eqref{eqn:discf} and \eqref{eqn:phikm} into \eqref{eqn:varform}, fixing $\mathbf{g}=(\Phi_0^{R_0,R_2})^*[Y_{p,q}]$ and choosing $\varphi=K_{\gamma,n}$ for $n=1,\ldots,M,$ as test functions, we are lead to the following system of linear equations
 \begin{align}
  \mathbf{M}\boldsymbol{\alpha}=\mathbf{d},\label{eqn:lineqfpq}
 \end{align}
 where 
 \small\begin{align*}
 &\mathbf{M}=\begin{pmatrix*}[l]\displaystyle\frac{1}{\lambda}\left\langle K_{\gamma,m} ,K_{\gamma,n}\right\rangle_{W^{1,2}(\SS_{R_2})}+\left(\frac{R_2}{4\pi}\right)^2\left\langle\nabla  F_{\frac{\gamma|\cdot|}{R_2},m},\nabla  F_{\frac{\gamma|\cdot|}{R_2},n}\right\rangle_{L^2(\Gamma_{R_0},\R^3)}
 \\[1.75ex]\displaystyle+\left(\frac{R_2}{R_1}\right)^2\left\langle K_{\frac{\gamma R_1}{R_2},m},K_{\frac{\gamma R_1}{R_2},n}\right\rangle_{L^2(\SS_{R_1})}\end{pmatrix*}_{n,m=1,\ldots,M},
 \\&\boldsymbol{\alpha}=\begin{pmatrix*} \alpha_m\end{pmatrix*}_{m=1,\ldots,M},
  \\&\mathbf{d}=\begin{pmatrix*} \displaystyle\frac{R_2^2}{4\pi(2p+1)}\left(\frac{R_0}{R_2}\right)^p\langle\nabla H_{p,q}^{R_0},\nabla F_n\rangle_{L^2(\Gamma_{R_0},\R^3)}\end{pmatrix*}_{n=1,\ldots,M}.
 \end{align*}\normalsize
A function $f$ of the form \eqref{eqn:discf}, determined by coefficients $\alpha_m$, $m=1,\ldots,M$, which solve \eqref{eqn:lineqfpq} will from now on be denoted as $f_{p,q}$. We use $f_{p,q}$ as an approximation of the solution to \eqref{eqn:varform} for the choice $\mathbf{g}=(\Phi_0^{R_0,R_2})^*[Y_{p,q}]$.

\begin{figure}
\begin{center}\footnotesize
\color{black} Input data $\Phi$\qquad\qquad\qquad\qquad\qquad\qquad\qquad\quad Input data $\Phi$\\
\includegraphics[scale=0.45]{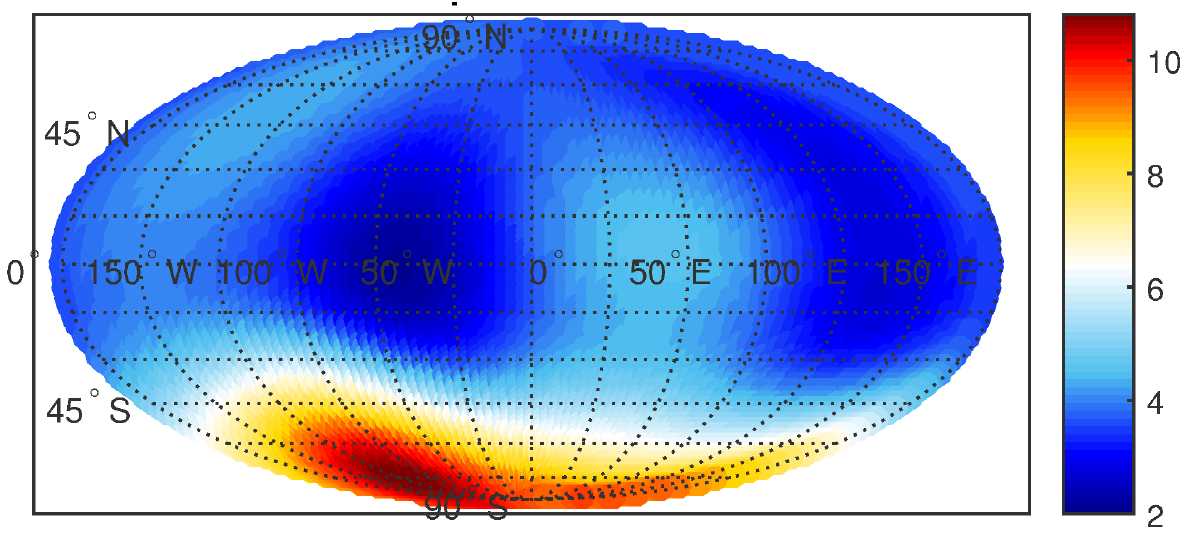}\qquad\qquad\includegraphics[scale=0.45]{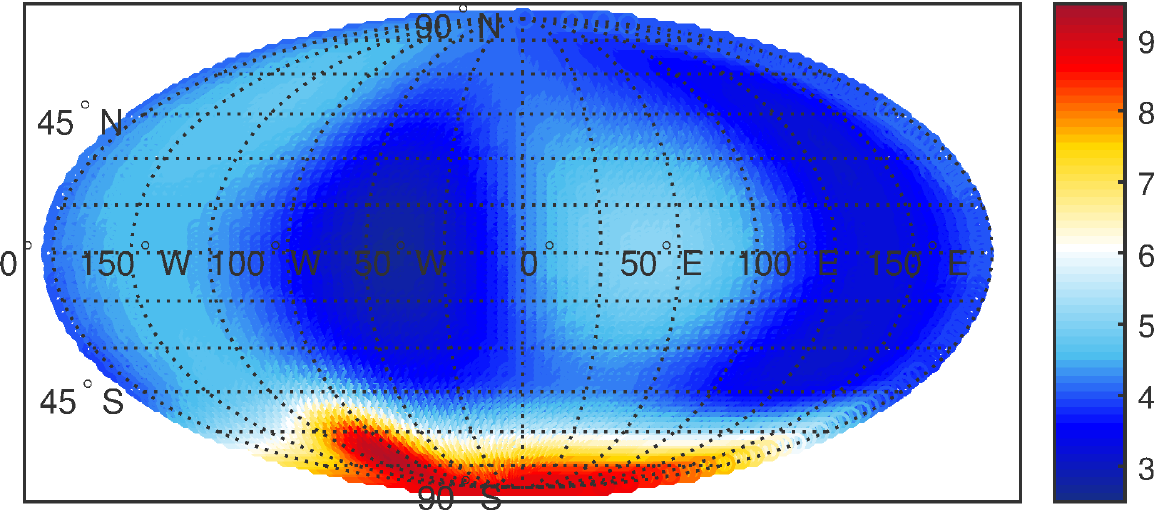}
\caption{Spatial plot of the input data $\Phi$ with parameters $\gamma_1=\frac{1}{20}$, $\gamma_2=\frac{1}{2}$ (\textit{left}) and $\gamma_1=\frac{3}{5}$, $\gamma_2=\frac{3}{5}$ (\textit{right}) for the magnetization $\m$ from \eqref{eqn:trueh}.\\[2ex]}\label{fig:spect1} \color{black} 
\color{black} \qquad \qquad\quad \qquad True $R_p$, $R_p^0$\qquad\quad\qquad\qquad\qquad\qquad\qquad True $R_p^0$ and reconstructions $\overline{R_p^0}$\\
\includegraphics[scale=0.45]{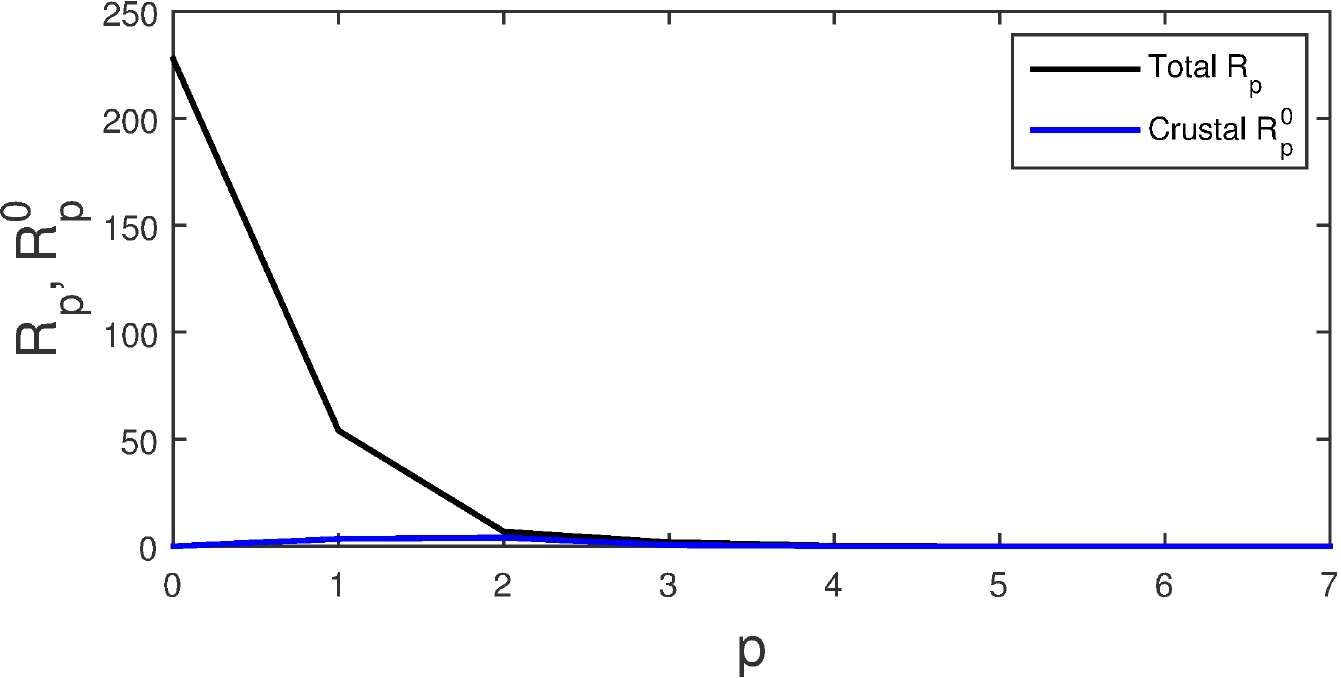}\qquad\qquad\quad\includegraphics[scale=0.45]{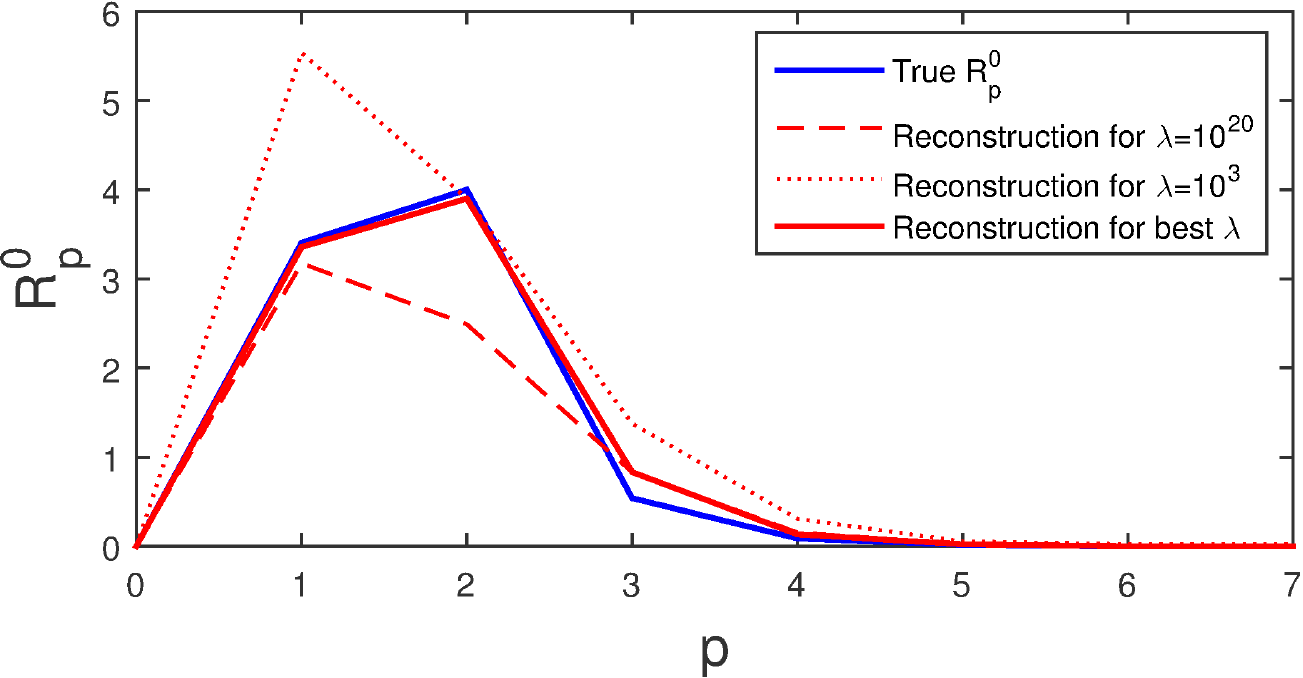}\\[2ex]
\color{black} \qquad \qquad\quad \qquad True $R_p$, $R_p^0$\qquad\quad\qquad\qquad\qquad\qquad\qquad True $R_p^0$ and reconstructions $\overline{R_p^0}$\\
\includegraphics[scale=0.45]{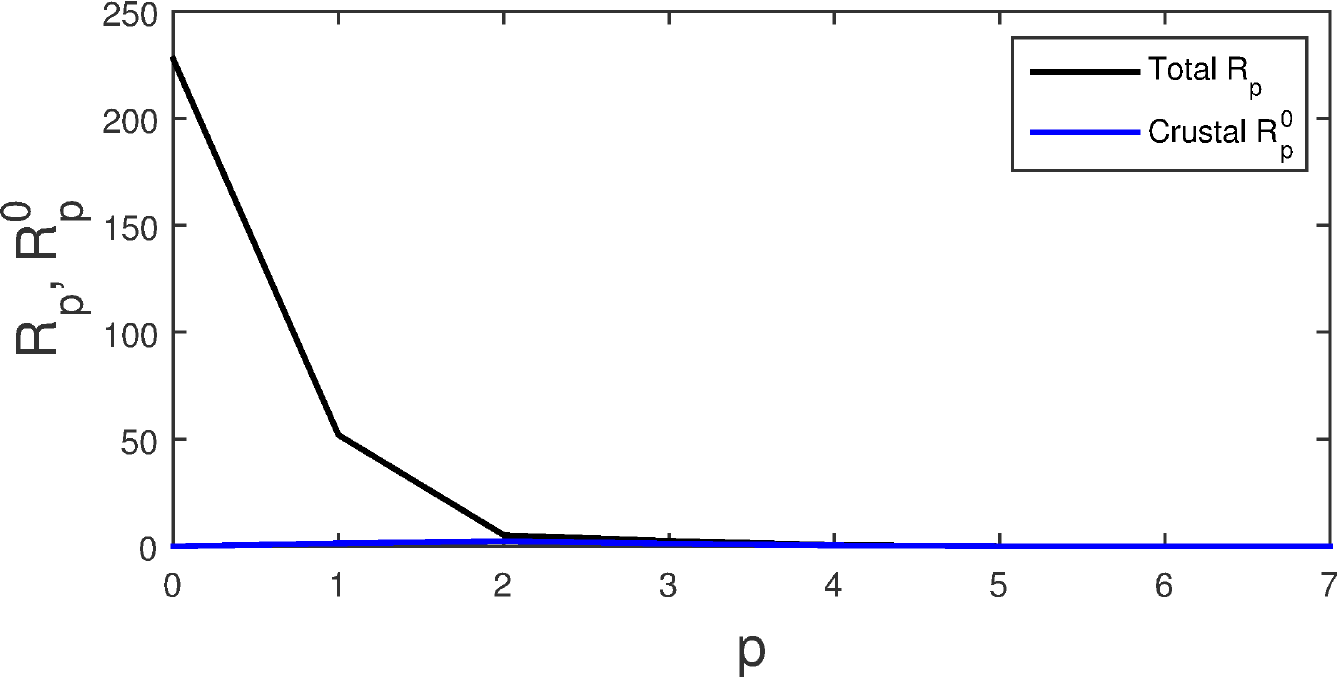}\qquad\qquad\quad\includegraphics[scale=0.45]{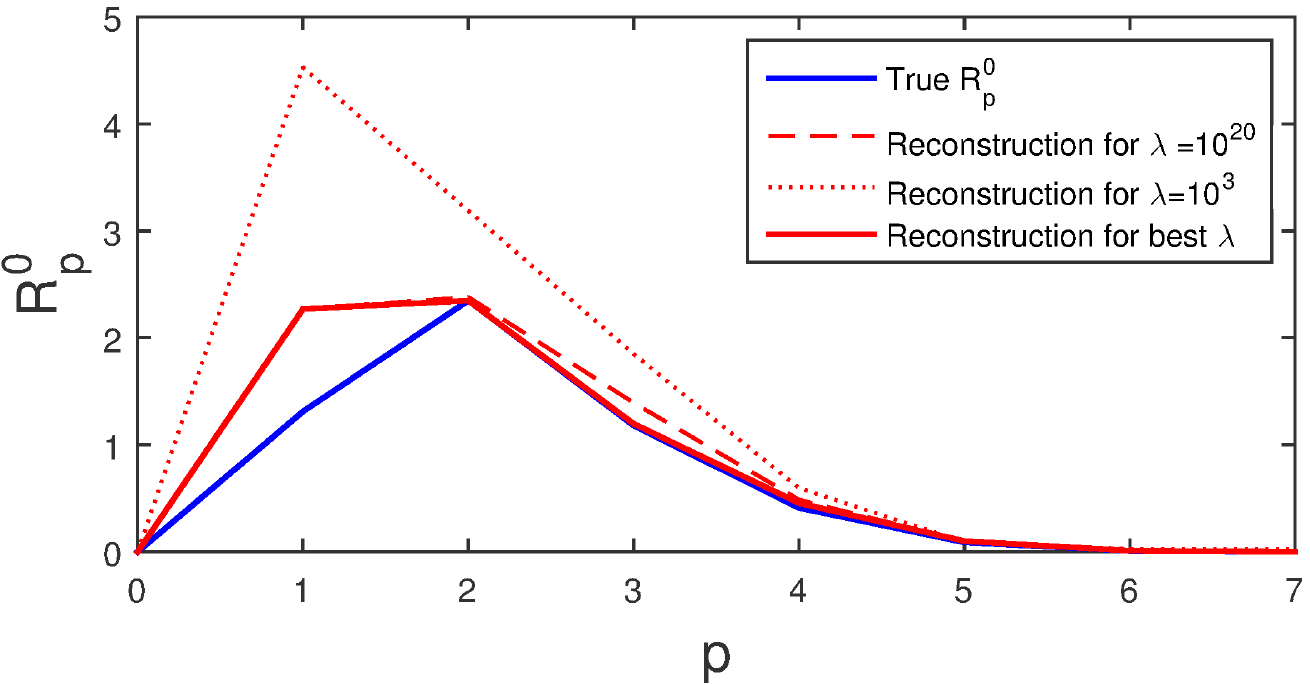}
\caption{\textit{Left}: Power spectrum $R_p$ of the input data $\Phi$ and power spectrum $R_p^0$ of the crustal contribution $\Phi_0$. \textit{Right}: True crustal power spectrum $R_p^{0}$ (blue) and reconstructed power spectrum $\overline{R_p^{0}}$ (red) for different parameters $\lambda$. The \textit{top row} shows the results for the parameters $\gamma_1=\frac{1}{20},\gamma_2=\frac{1}{2}$ and the \textit{bottom row} for $\gamma_1=\frac{3}{5},\gamma_2=\frac{3}{5}$. \\[2ex]}\label{fig:spect2}
\color{black} Power Spectrum of $f_{1,1}$\qquad\qquad\qquad\qquad\qquad\qquad\qquad\quad Power Spectrum of $f_{50,1}$\\
\includegraphics[scale=0.46]{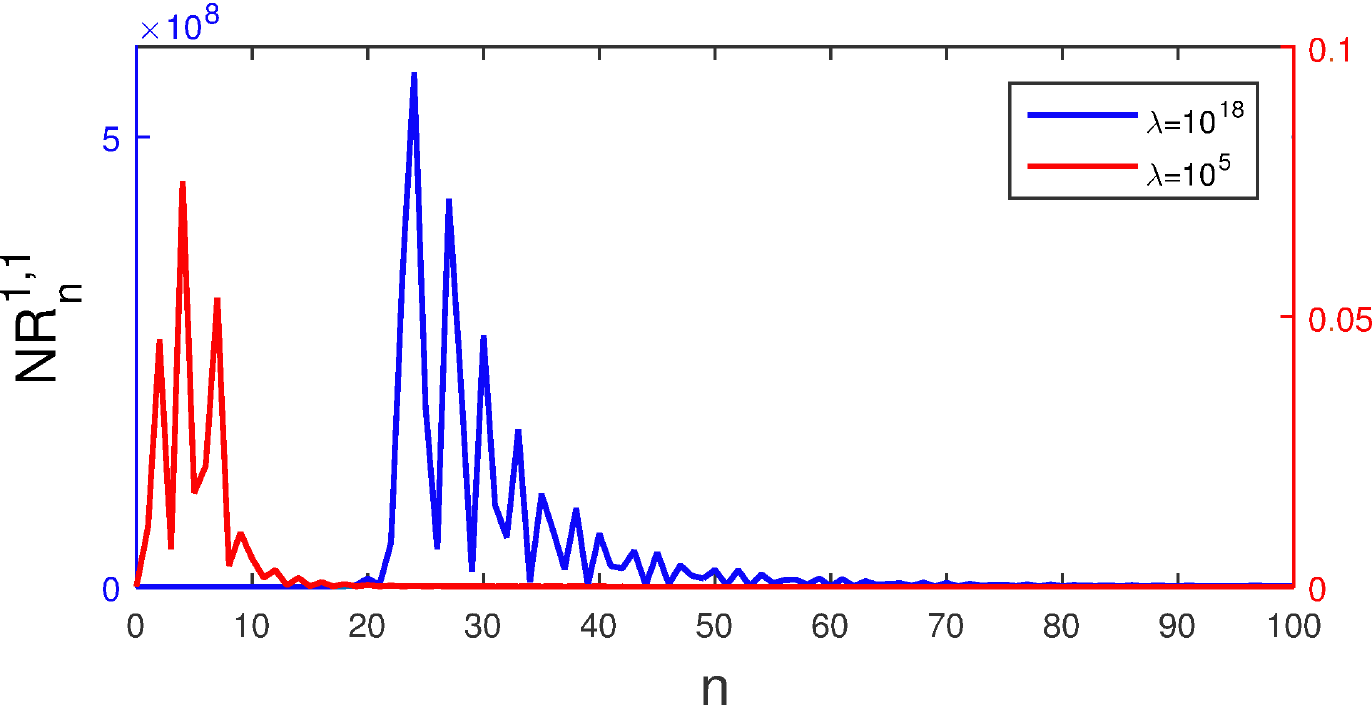}\qquad\qquad\includegraphics[scale=0.465]{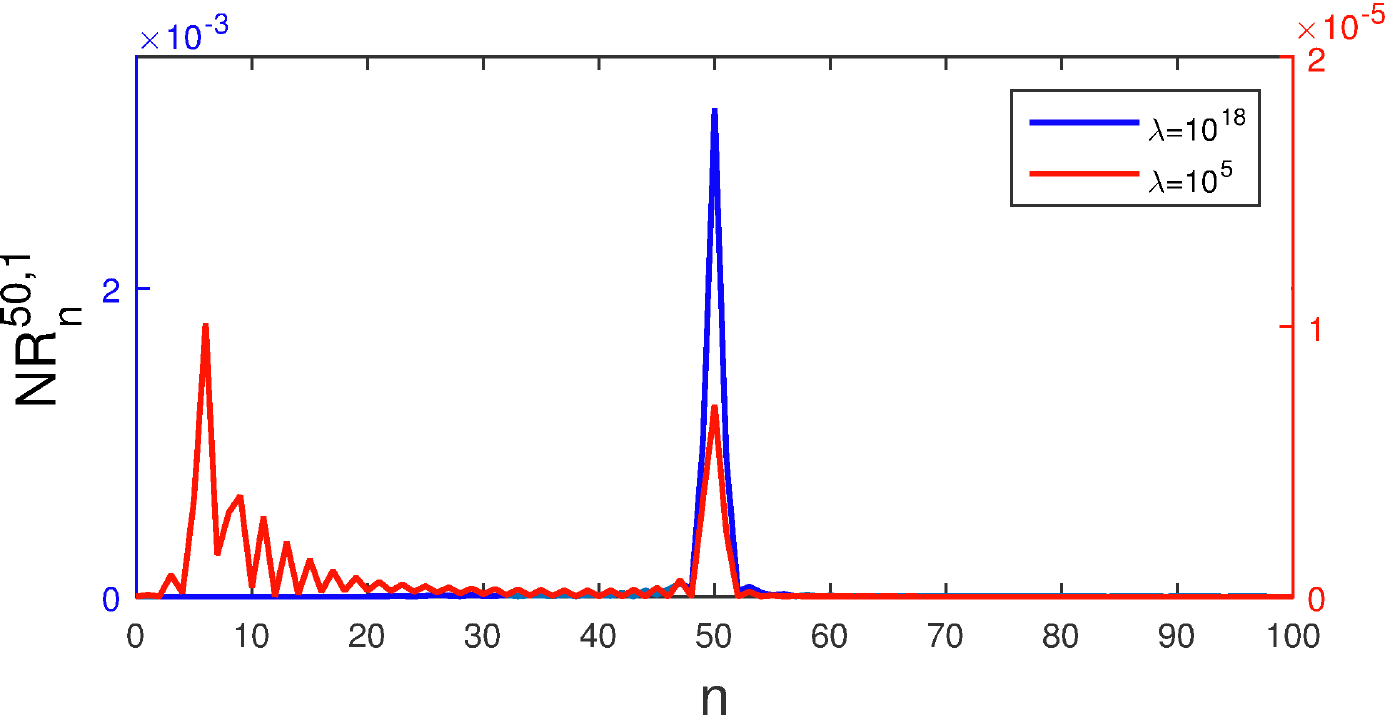}
\caption{\color{black} Scaled power spectrum $NR_n^{p,q}$ for $p=1$, $q=1$ (\textit{left}) and $p=50$, $q=1$ (\textit{right}).\\[4ex]}\label{fig:spect3}
 \end{center}
\end{figure}

\subsubsection*{A Numerical Example}

In order to generate input data $\Phi=\Phi^{R_1,R_0,R_2}[\mathbf{m},h]=\Phi_0^{R_0,R_2}[\mathbf{m}]+\Phi_1^{R_1,R_2}[h]$ for a test example, we choose
\begin{align}
 \mathbf{m}(x)&=b_1\frac{x}{|x|} L_{\gamma_1}\left(\frac{x}{|x|}\cdot y_1\right)+b_2 \frac{x}{|x|}  L_{\gamma_2}\left(\frac{x}{|x|}\cdot y_2\right), \quad b_1=15,\,b_2=10,\nonumber
\\h(x)&=\sum_{n=0}^5\sum_{k=1}^{2n+1}a_{n,k} Y_{n,k}\left(\frac{x}{|x|}\right),\quad a_{0,1}=a_{1,1}=2^5,\,a_{2,5}=a_{3,5}=a_{4,5}=2^{4},\,a_{5,5}=2^3,\nonumber
\\& \qquad\qquad\qquad\qquad\qquad\qquad\quad a_{n,k}=0 \textnormal{ else},\label{eqn:trueh}
\end{align}\normalsize
with $y_1=(0,0,-1)^T$ and $y_2=(0,\frac{1}{2},-\frac{\sqrt{3}}{2})^T$. The functions $L_{\gamma_i}$ are chosen as follows:
\begin{align}
 L_{\gamma_i}(t)=\left\{\begin{array}{ll}
                         0,&t\in[-1,\gamma_i),
                         \\\frac{(t-\gamma_i)^k}{(1-\gamma_i)^k},&t\in[\gamma_i,1],
                        \end{array}\right.
\end{align}
for $k=3$. These functions have been studied in more detail in \cite{schreiner97} and are suited for our purposes since they are compactly supported and allow a recursive computation of the Fourier coefficients of $\mathbf{m}$. \color{black}The parameters $\gamma_i\in (-1,1)$ reflect the localization of $L_{\gamma_i}$ (a parameter $\gamma_i$ close to one means a strong localization). In our test examples, we investigate the two setups $\gamma_1=\frac{1}{20},\gamma_2=\frac{1}{2}$ and $\gamma_1=\frac{3}{5},\gamma_2=\frac{3}{5}$, where latter reflects a slightly stronger localization of the underlying magnetization. \color{black} The (unknown) crustal contribution is then denoted by $\Phi_0=\Phi_0^{R_0,R_2}[\mathbf{m}]$ and the  (unknown) core contribution by $\Phi_1=\Phi_1^{R_1,R_2}[h]$. For the involved radii, we choose $R_0=1$ and $R_2=1.06$ (at scales of the Earth, the latter indicates a realistic satellite altitude of about 380km above the Earth's surface) and $R_1=0.5$ (at scales of the Earth, this is a rough approximation of the radius of the outer core). The subregion $\Gamma_{R_0}=\{x\in\SS_{R_0}:x\cdot(0,0,1)^T\leq0\}$ is set to be the Southern hemisphere and the chosen magnetizations of the form \eqref{eqn:trueh} satisfy $\supp(\m)\subset\Gamma_{R_0}$. \color{black} For our computations, we use the localization parameter $\gamma=0.95$ and choose $M=8,499$ uniformly distributed centers $x_m\in\SS_1$, $m=1,\ldots,M,$ for the kernels $K_{\gamma,m}$. All numerical integrations 
necessary during the procedure are performed via the methods of \cite{driscoll94} (when the integration region comprises the entire sphere $\SS_{R_0}$, $\SS_{R_1}$, or $\SS_{R_2}$, respectively) and \cite{hesse12} (when the integration is only performed over the spherical cap $\SS_{R_0}\setminus\Gamma_{R_0}$).
The input data for the two different setups associated with $\gamma_1,\gamma_2$ are shown in Figure \ref{fig:spect1}. These setups are not based on real geomagnetic data but they reflect a typical geomagnetic situation in the sense that the core contribution clearly dominates the crustal contribution at low spherical harmonic degrees. Figure \ref{fig:spect2} shows that an empirical separation by a sharp cut-off at degree $p=2$ or $p=3$ would neglect relevant information in the crustal contribution.

According to  Corollary \ref{cor:approxcoeffs}, an approximation of the Fourier coefficient $\langle \Phi_0,Y_{p,q}\rangle_{L^2(\SS_{R_2})}$ of the crustal contribution $\Phi_0$ is now given by $\langle \Phi,f_{p,q}\rangle_{L^2(\SS_{R_2})}$, with $f_{p,q}$ of the form described in the previous subsection. We do this for various degrees $p$ and orders $q$ and we illustrate the results in terms of power spectra: The crustal power spectrum is defined as
\begin{align*}
R_p^0=R_p[\Phi_0]=\sum_{q=1}^{2p+1}\left|\langle \Phi_0,Y_{p,q}\rangle_{L^2(\SS_{R_2})}\right|^2,\quad p\in\mathbb{N}_0.
\end{align*}
Our approximated power spectrum is then of the form
\begin{align*}
\overline{R_p^{0}}=\sum_{q=1}^{2p+1}\left|\langle \Phi,f_{p,q}\rangle_{L^2(\SS_{R_2})}\right|^2,\quad p\in\mathbb{N}_0.
\end{align*}
The power spectrum of the input signal $\Phi$ (i.e., the superposition of the crustal and core contribution) is analogously defined by $R_p=R_p[\Phi]=\sum_{q=1}^{2p+1}|\langle \Phi,Y_{p,q}\rangle_{L^2(\SS_{R_2})}|^2$. 

Figure \ref{fig:spect2} shows the reconstructed power spectra and we see that they yield good results (for a well-chosen  parameter $\lambda$),
in both  setups under investigation. Stronger deviations mainly occur at lower spherical harmonic degrees $p$. The solid red spectrum in Figure \ref{fig:spect2} indicated as 'Reconstruction for best $\lambda$' does not reflect the result for a
single choice of $\lambda$ but rather for (possibly different) 
best $\lambda$ in each degree $p$ of the spectrum. The setup for magnetizations $\m$ with parameters $\gamma_1=\frac{3}{5},\gamma_2=\frac{3}{5}$ was chosen to investigate magnetizations with a slightly stronger localization, meaning that the corresponding potential $\Phi_0$ has slightly stronger contributions at higher spherical harmonic degrees than for the setup $\gamma_1=\frac{1}{20},\gamma_2=\frac{1}{2}$ (compare the right hand images in Figure \ref{fig:spect2}). 
\color{black} In Figure \ref{fig:spect3}, we illustrate the effects mentioned in Remark \ref{rem:shdisc} by observing the scaled power spectrum $N_n^{p,q}=NR_n[f_{p,q}]=\frac{1}{2n+1}R_n[f_{p,q}]=\frac{1}{2n+1}\sum_{k=1}^{2n+1}|\langle f_{p,q},Y_{n,k}\rangle_{L^2(\SS_{R_2})}|^2$, $n\in\mathbb{N}_0$, for $p=1$, $q=1$, and $p=50$, $q=1$ (we scaled by a factor $\frac{1}{2n+1}$ solely to get a better idea of the 
average strength of the Fourier coefficients $|\langle f_{p,q},Y_{n,k}\rangle_{L^2(\SS_{R_2})}|$, $k=1,\ldots, 2n+1$, for  fixed degree $n$). As expected 
from Remark \ref{rem:shdisc}, larger Lagrange parameters $\lambda$ (which 
correspond to smaller $\eps$) result in a shift of the major contributions 
of the power spectrum towards higher spherical harmonic degrees. However, for $p=50$, $q=1$, the major spike around $n=50$ remains, somewhat motivating a different behaviour of the Fourier coefficients of $f_{p,q}$ for larger degrees $p$ compared to smaller $p$. 

\color{black}

\subsection{Approximate Reconstruction of $\Phi_0$}\label{sec:phi0num}

While the previous section aimed at the reconstruction of the Fourier coefficients of $\Phi_0$, we are now concerned with the reconstruction the magnetization $\mathbf{m}$ that generates $\Phi_0$. Actually, the goal is still an approximation of $\Phi_0$, but instead of solving multiple extremal problems like
Problem \ref{prob:2}  we rather  solve  a single 
least-squares problem to get an 
approximation $\bar{\mathbf{m}}$ of $\mathbf{m}$, and then we compute $\bar{\Phi}_0=\Phi_0^{R_0,R_2}[\bar{\mathbf{m}}]$ to approximate  $\Phi_0$. 
Beyond the instrumental parametrizations from the previous section, the only
input we retain from the rest of the paper is that, since we apply the 
technique on an example where $\Gamma_{R_0}\not=\SS_{R_0}$, we know that 
separation of the core and crustal potentials is possible by 
Corollary \ref{cor:unique}. Still, we gather from Theorem \ref{thm:unique}
that $\mathbf{m}$ is not uniquely determined though  $\Phi_0$ is. So, in order to regularize the problem, we use standard penalization term 
to compute a candidate $\bar{\mathbf{m}}$ of small norm (weighted by $\alpha$).
More precisely, we consider the following extremal problem.

\begin{prob}\label{prob:3}
Let the setup be as in Definition \ref{def:ops} with
$\Gamma_{R_0}\not=\SS_{R_0}$, and let $\Phi\in L^2(\SS_{R_2})$ be given. Then, for  fixed parameters $\alpha,\beta>0$, find $\bar{\mathbf{m}}\in W^{2,2}(\SS_{R_0},\R^3)$ and $\bar{h}\in W^{2,2}(\SS_{R_1})$ to minimize 
\small\begin{align*}
\inf_{\atopp{\bar{\mathbf{m}}\in W^{2,2}(\SS_{R_0},\R^3),}{\bar{h}\in W^{2,2}(\SS_{R_1})}}&\left\|\Phi - \Phi^{R_1,R_0,R_2}[\bar{\mathbf{m}},\bar{h}]\right\|_{L^2(\SS_{R_2})}^2+\alpha\|(\bar{\mathbf{m}},\bar{h})\|_{W^{2,2}(\SS_{R_0},\R^3)\times W^{2,2}(\SS_{R_1})}^2+\beta \|\bar{\mathbf{m}}\|_{L^2(\SS_{R_0}\setminus\Gamma_{R_0},\R^3)}^2.
\end{align*}\normalsize
Note that in this particular setup, the integration in the definition of the
operator $\Phi_0^{R_0,R_2}$ is meant over the entire sphere 
$\SS_{R_0}$ and not just over $\Gamma_{R_0}$. 
\end{prob}

\begin{rem}
Another (more natural) choice to obtain approximations of $\mathbf{m}$ and $h$ would be to minimize
\begin{align}
\inf_{\atopp{\bar{\mathbf{m}}\in W^{2,2}(\Gamma_{R_0},\R^3),}{\bar{h}\in W^{2,2}(\SS_{R_1})}}&\left\|\Phi - \Phi^{R_1,R_0,R_2}[\bar{\mathbf{m}},\bar{h}]\right\|_{L^2(\SS_{R_2})}^2+\alpha\|(\bar{\mathbf{m}},\bar{h})\|_{W^{2,2}(\Gamma_{R_0},\R^3)\times W^{2,2}(\SS_{R_1})}^2,\label{eqn:minfunc}
\end{align}
where this time the integration defining $\Phi_0^{R_0,R_2}$ 
is only over $\Gamma_{R_0}$ (as always in this paper, with the exception of Problem \ref{prob:3} and Section \ref{sec:gammas}). Solving \eqref{eqn:minfunc} leads to magnetizations $\bar{\mathbf{m}}$ that are of class $W^{2,2}(\Gamma_{R_0},\R^3)$, while solving Problem \ref{prob:3} leads to magnetizations $\bar{\mathbf{m}}$ that are of class $W^{2,2}(\SS_{R_0},\R^3)$ and localization in $\Gamma_{R_0}$ has to be enforced by adding a penalty term (weighted by $\beta$). However, for the upcoming example, the minimization proposed in Problem \ref{prob:3} yielded slightly better results. Furthermore, it allowed an easier illustration of the effect of the localization constraint by simply dropping the penalty term (i.e., setting $\beta=0$). Existence of minimizers is guaranteed
in both cases by standard arguments. The typically difficult choice of  parameters $\alpha, \beta$ will not be discussed here. In the provided examples, 
we simply chose those parameters that seemed to yield the best results when compared to the ground truth.
\end{rem}

\subsubsection*{Discretization}
In order to discretize Problem \ref{prob:3}, we expand $\bar{\mathbf{m}}$ and $\bar{h}$ in terms of Abel-Poisson kernels the way indicated in Section \ref{sec:fouriernum}:
\begin{align*}
 \bar{\mathbf{m}}(x)&=\sum_{i=1}^3\sum_{n=1}^N\bar{\alpha}_{i,n} \,o^{(i)} K_{\gamma,n}\left(x\right),\quad x\in\SS_{R_0},
 \\\bar{h}(x)&=\sum_{n=1}^N\bar{\beta}_{n}\,K_{\gamma,n}\left(x\right), \quad x\in\SS_{R_1}.
\end{align*}
For brevity, the vectorial operators $o^{(i)}$ have been introduced to denote $o^{(1)}=\nu\,\textnormal{Id}$, $o^{(2)}=\nabla_\SS$, and $o^{(3)}=\textnormal{L}_\SS$ (with $\nu$ denoting the unit normal vector). Such localized kernels are suitable here since we know/assume in advance that the sought-after magnetization $\m$ is localized in some subregion $\Gamma_{R_0}$. Using this discretization, the minimization of Problem \ref{prob:3} reduces to solving the following set of linear equations for the coefficients $\bar{\alpha}_{i,n}$ and $\bar{\beta}_n$:
\begin{align}
 \mathbf{M}\boldsymbol{\gamma}=\mathbf{d},\label{eqn:linsolve}
\end{align}
where
\begin{align*}
 \mathbf{M}&=\left(\begin{array}{c|c}
              \mathbf{A}&\mathbf{B}^T
              \\\hline\mathbf{B}&\mathbf{C}
             \end{array}\right)\in\R^{4N\times4N},\quad \boldsymbol{\gamma}&\!\!\!\!\!\!\!\!=(\overline{\boldsymbol{\beta}}\,\vline \,\overline{\boldsymbol{\alpha}}_j)_{j=1,2,3}^T\in\mathbb{R}^{4N},\quad\mathbf{d}&=(\mathbf{a}\,|\,\mathbf{b}_i)_{i=1,2,3}^T\in\mathbb{R}^{4N},
 \end{align*}
 with
 \begin{align*}
 \mathbf{A}&=\begin{pmatrix*}[l]\langle\Phi^1_{n},\Phi^1_{k}\rangle_{L^2(\SS_{R_2})}+\alpha\langle K_{\gamma,n},K_{\gamma,k}\rangle_{W^{2,2}(\SS_{R_1})} \end{pmatrix*}_{n,k=1,\ldots,N},\nonumber
\\\mathbf{B}&=(\mathbf{B}_{i})_{i=1,2,3},\quad \mathbf{B}_i=\begin{pmatrix*}\langle\Phi^0_{i,n},\Phi^1_{k}\rangle_{L^2(\SS_{R_2})}\end{pmatrix*}_{n,k=1,\ldots,N},
 \\\mathbf{C}&=(\mathbf{C}_{i,j})_{i,j=1,2,3},
\\\mathbf{C}_{i,j}&=\begin{pmatrix*}[l]\langle\Phi^0_{i,n,}\Phi^0_{j,k}\rangle_{L^2(\SS_{R_2})}+\alpha\left\langle o^{(i)} K_{\gamma,n},o^{(j)} K_{\gamma,k}\right\rangle_{W^{2,2}(\SS_{R_0},\R^3)}
\\[1.75ex]+\beta\left\langle o^{(i)} K_{\gamma,n},o^{(j)} K_{\gamma,k}\right\rangle_{L^2(\SS_{R_0}\setminus\Gamma_{R_0},\R^3)}\end{pmatrix*}_{n,k=1,\ldots,N},\nonumber
\\[1.75ex]
{\overline{\boldsymbol{\beta}}}&=(\overline{\beta}_k)_{k=1,\ldots,N},\quad\overline{\boldsymbol{\alpha}}_j=(\overline{\alpha}_{j,k})_{k=1,\ldots,N},
 \\[1.75ex]\mathbf{a}&=\left(\langle\Phi^1_{n},\Phi\rangle_{L^2(\SS_{R_2})}\right)_{n=1,\ldots,N},\quad \mathbf{b}_i=\left(\langle\Phi^0_{i,n},\Phi\rangle_{L^2(\SS_{R_2})}\right)_{n=1,\ldots,N},
  \end{align*}
and
\begin{align*}
\Phi^0_{i,n}(x)&=\frac{1}{4\pi}\int_{\SS_{R_0}} \left(o^{(i)} K_{\gamma,n}(y)\right)\cdot\frac{x-y}{|x-y|^3}\dd \omega_{R_0}(y),
\\\Phi^1_{n}(x)&=\frac{1}{4\pi R_1}\int_{\SS_{R_1}} K_{\gamma,n}(y)\frac{|x|^2-R_1^2}{|x-y|^3}\dd \omega_{R_1}(y).
\end{align*}
Again, all necessary numerical integrations are performed via the methods of \cite{driscoll94} (when the integration region comprises the entire sphere $\SS_{R_0}$, $\SS_{R_1}$, or $\SS_{R_2}$, respectively) and \cite{hesse12} (when the integration is only performed over the spherical cap $\SS_{R_0}\setminus\Gamma_{R_0}$).

\subsubsection*{A Numerical Example}

\begin{figure}
\begin{center}\footnotesize
Input data $\Phi$\\
\includegraphics[scale=0.45]{phiinput-eps-converted-to.pdf}
\begin{tabbing}
 \qquad\qquad\qquad\quad True $\Phi_0$\qquad\qquad\qquad\qquad\quad \= Reconstructed $\bar{\Phi}_0$ \qquad\qquad\qquad\= Reconstructed $\bar{\Phi}_0$ 
 \\\>($\alpha=5\cdot10^{-16}$, $\beta=1$)\>($\alpha=5\cdot10^{-16}$, $\beta=0$)\\[-5ex]
 \end{tabbing}
\includegraphics[scale=0.37]{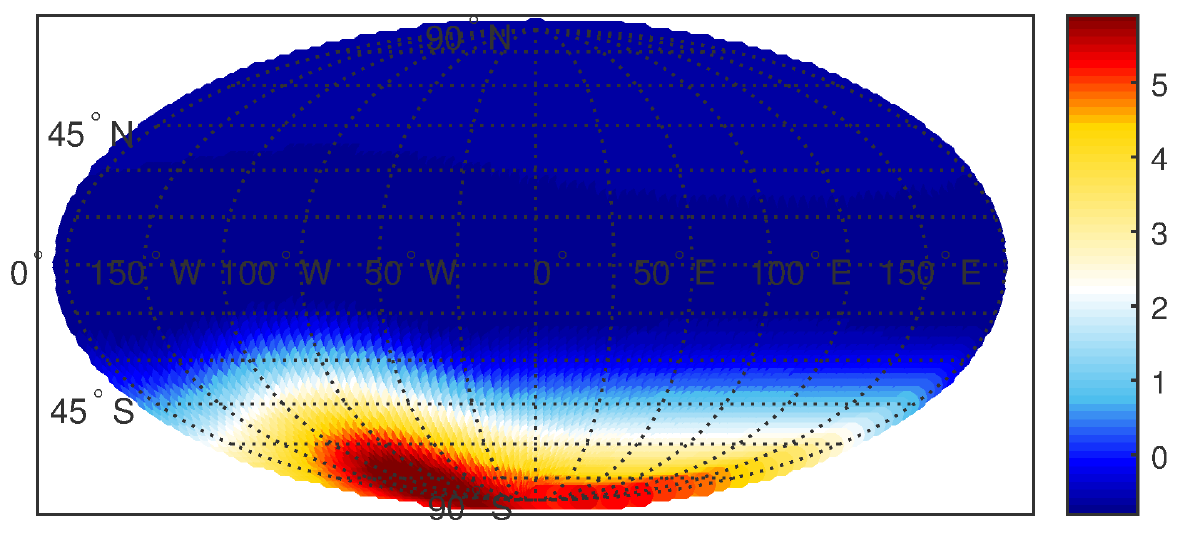}\quad\,\includegraphics[scale=0.37]{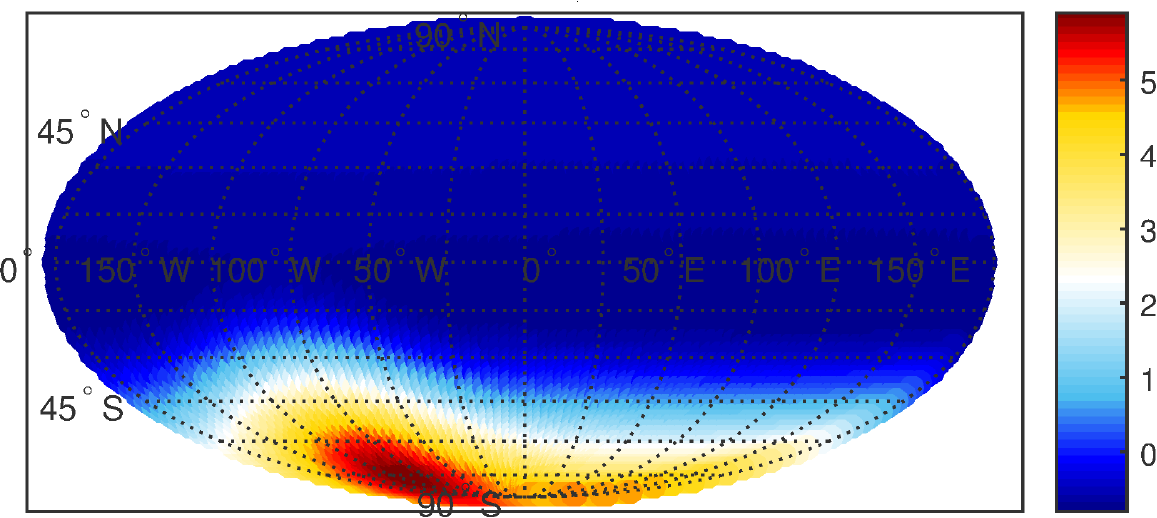}\quad\,\includegraphics[scale=0.37]{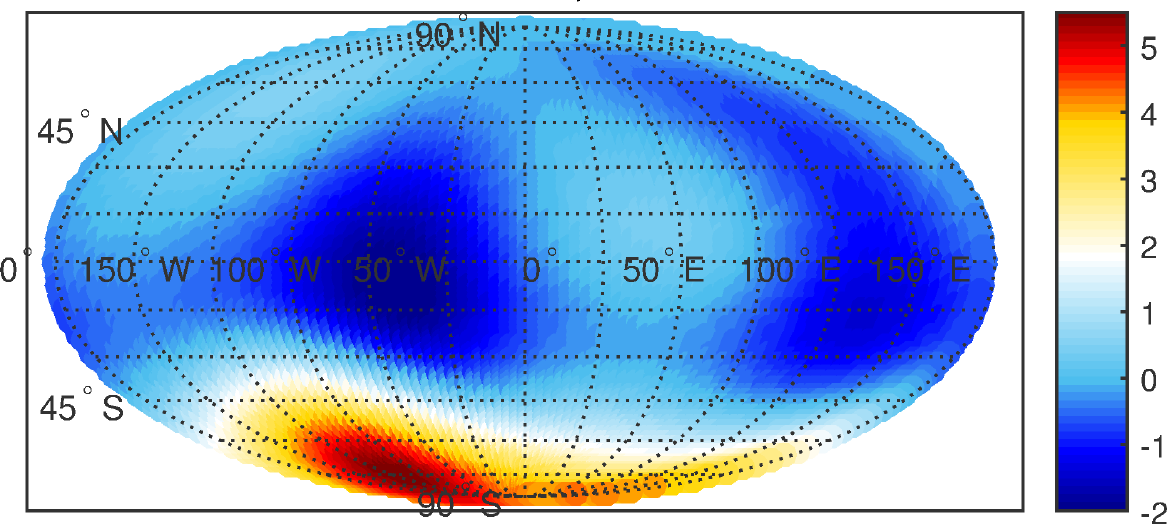}
 \begin{tabbing}
 \qquad\qquad\qquad\quad True $\Phi_1$\qquad\qquad\qquad\qquad\quad \= Reconstructed $\bar{\Phi}_1$ \qquad\qquad\qquad\= Reconstructed $\bar{\Phi}_1$ 
 \\\>($\alpha=5\cdot10^{-16}$, $\beta=1$)\>($\alpha=5\cdot10^{-16}$, $\beta=0$)\\[-5ex]
 \end{tabbing}
 \includegraphics[scale=0.37]{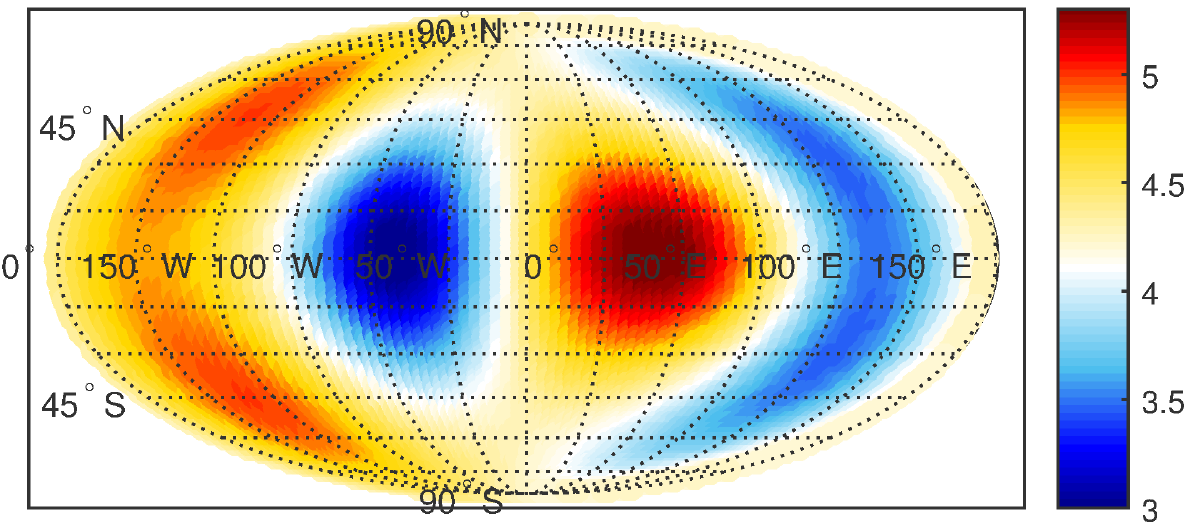}\quad\includegraphics[scale=0.37]{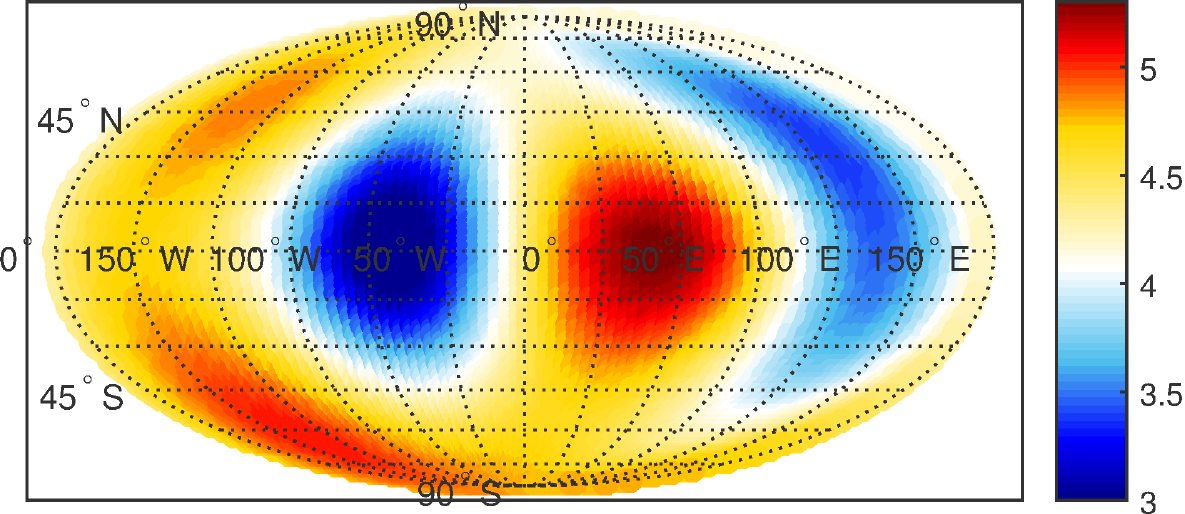}\quad\includegraphics[scale=0.37]{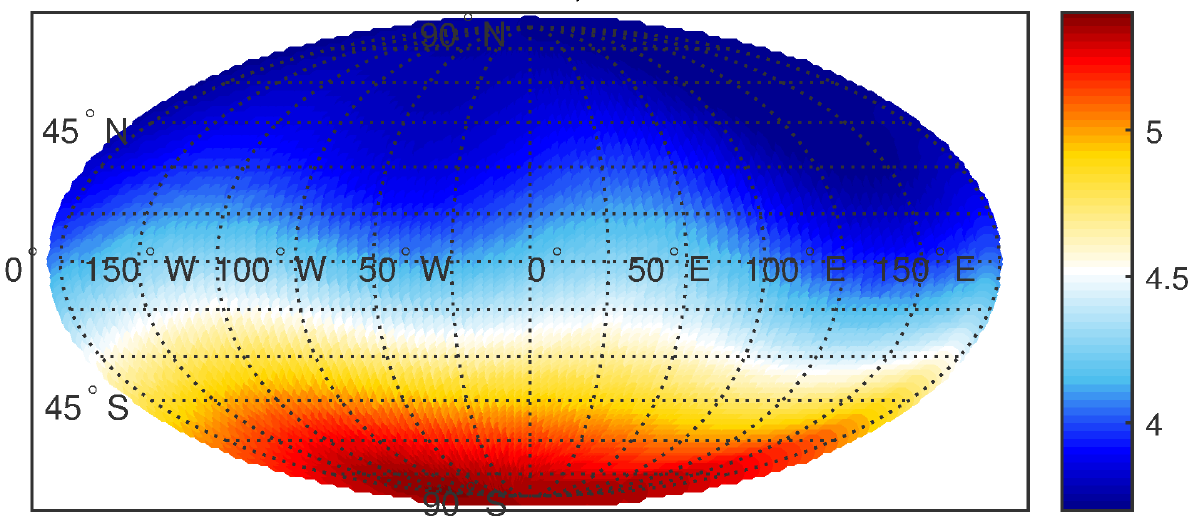}
 \caption{Results for radii $R_1=0.5$, $R_0=1$, and $R_2=1.06$: Input data $\Phi=\Phi_0+\Phi_1$ (\emph{top}), ground truth $\Phi_0$, $\Phi_1$ (\emph{bottom left}), reconstructed $\bar{\Phi}_0$, $\bar{\Phi}_1$ with localization constraint (\emph{bottom center}), and reconstructed $\bar{\Phi}_0$, $\bar{\Phi}_1$ without localization constraint (\emph{bottom right}). \\[2ex]}\label{fig:ex1a}

 Input data $\Phi$\\
\includegraphics[scale=0.45]{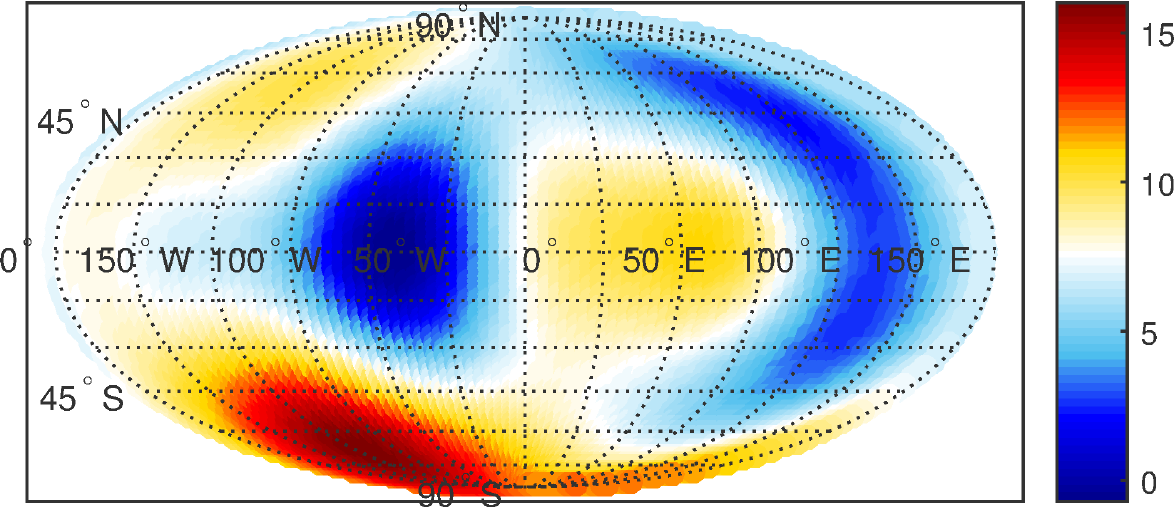}
\begin{tabbing}
 \qquad\qquad\qquad\qquad\qquad\qquad\qquad\quad True $\Phi_0$\qquad\qquad\qquad\qquad \= Reconstructed $\bar{\Phi}_0$ 
 \\\>($\alpha=5\cdot10^{-15}$, $\beta=1$)\\[-5ex]
 \end{tabbing}
\includegraphics[scale=0.37]{phi0true-eps-converted-to.pdf}\quad\,\includegraphics[scale=0.37]{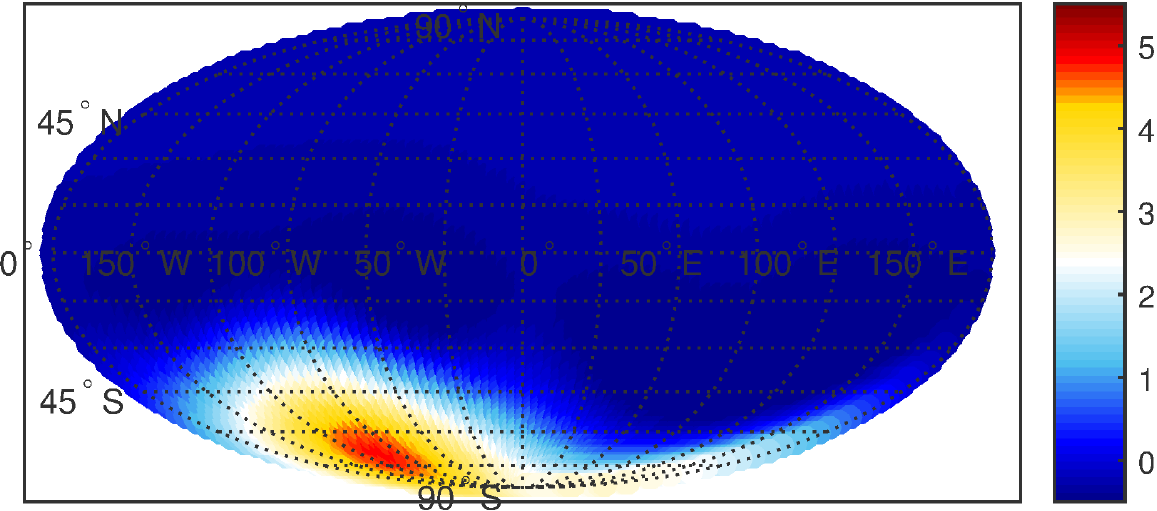}
 \begin{tabbing}
 \qquad\qquad\qquad\qquad\qquad\qquad\qquad\quad True $\Phi_1$\qquad\qquad\qquad\qquad \= Reconstructed $\bar{\Phi}_1$
 \\\>($\alpha=5\cdot10^{-15}$, $\beta=1$)\\[-5ex]
 \end{tabbing}
 \includegraphics[scale=0.37]{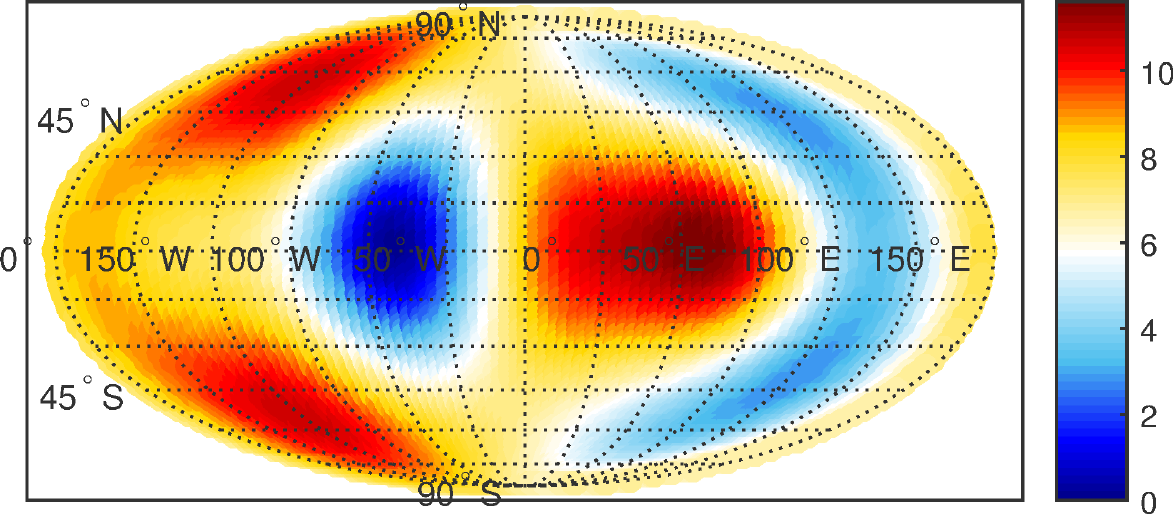}\quad\includegraphics[scale=0.37]{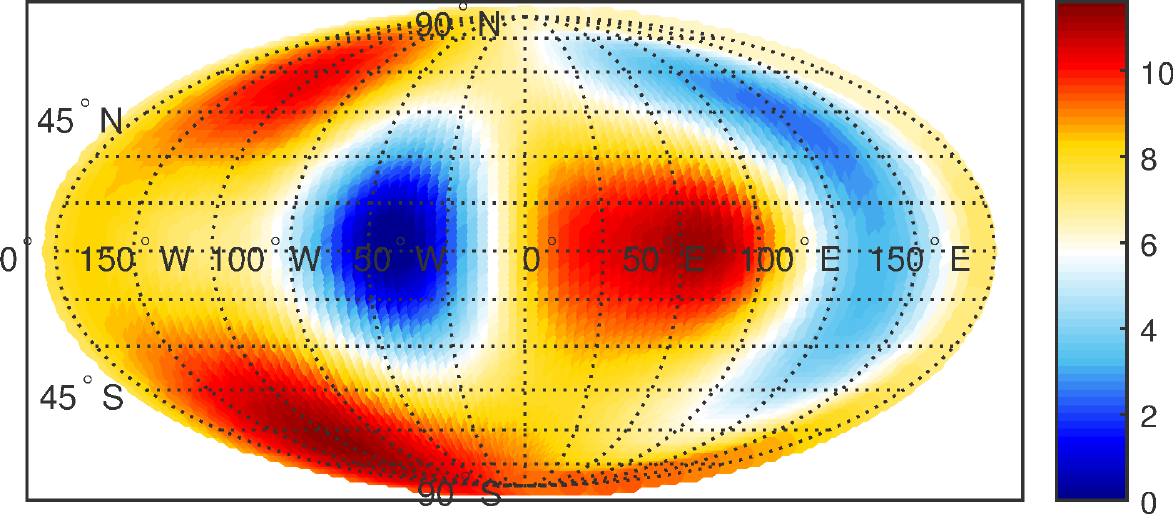}
 \caption{Results for radii $R_1=0.8$, $R_0=1$, and $R_2=1.06$: Input data $\Phi=\Phi_0+\Phi_1$ (\emph{top}), ground truth $\Phi_0$, $\Phi_1$ (\emph{bottom left}), and reconstructed $\bar{\Phi}_0$, $\bar{\Phi}_1$ with localization constraint (\emph{bottom right}).}\label{fig:ex1b}
\end{center}
\end{figure}

We use the same setup as in 
Section \ref{sec:fouriernum} \color{black} (with parameters $\gamma_1=\frac{1}{20}$, $\gamma_2=\frac{1}{2}$) \color{black} to generate $\Phi=\Phi^{R_1,R_0,R_2}[\mathbf{m},h]$, $\Phi_0=\Phi_0^{R_0,R_2}[\mathbf{m}]$, and $\Phi_1=\Phi_1^{R_1,R_2}[h]$. In the discretization above, we choose $\gamma=0.9$ and take $N=10,235$ uniformly distributed centers $x_n\in \SS_1$, $n=1,\ldots,N$. 
As in the previous example, we choose radii $R_0=1$, $R_2=1.06$, and now additionally vary $R_1$ between $0.5$ and $0.8$. The subregion $\Gamma_{R_0}$ is again the Southern hemisphere $\{x\in\SS_{R_0}:x\cdot (0,0,1)^T<0\}$. Approximations of $\bar{\m}$ and $\bar{h}$ are obtained by solving \eqref{eqn:linsolve}.

In Figure \ref{fig:ex1a}, we illustrate the potentials $\bar{\Phi}_0=\Phi_0^{R_0,R_2}[\bar{\mathbf{m}}]$ and  $\bar{\Phi}_1=\Phi_1^{R_1,R_2}[\bar{h}]$ corresponding to the reconstructed $\bar{\mathbf{m}}$ and $\bar{h}$ for radius $R_1=0.5$, while in Figure \ref{fig:ex1b} we set  $R_1=0.8$. In the first case, we see that the reconstructions yield good approximations of the ground truths $\Phi_0=\Phi_0^{R_0,R_2}[\mathbf{m}]$ and $\Phi_1=\Phi_1^{R_1,R_2}[h]$. However, Figure \ref{fig:ex1b} suggests that the reconstruction of the potential $\Phi_0$ becomes numerically more critical as the spheres $\SS_{R_1}$ and $\SS_{R_0}$ get closer. 
The influence of the localization constraint on the reconstruction can be seen on the right set of images in Figure \ref{fig:ex1a}: neglecting the localization constraint (i.e., choosing $\beta=0$) leads to a wrong separation of the contributions $\bar{\Phi}_0$ and $\bar{\Phi}_1$.

\section{Conclusion}\label{sec:conc}
In this paper, we set up a geophysically reasonable model of the core and crustal magnetic field potentials $\Phi_1$ and $\Phi_0$ respectively, for which we showed that each single potential can be recovered uniquely if only the superposition $\Phi=\Phi_0+\Phi_1$ is known on an external sphere $\SS_{R_2}$. Furthermore, we supplied first approaches to the reconstruction of $\Phi_0$ 
and of its Fourier coefficients. The latter is particularly interesting as it would allow a comparison with  the empirical approach to separation based on  a sharp cut-off in the power spectrum of $\Phi$. Two main directions 
call for further study: (1) the geophysical post-processing of real geomagnetic data in order to back up (or deny) the assumption that $\m$ is supported in a subregion $\Gamma_{R_0}$ of the Earth's surface; (2) improving numerical schemes allowing reconstruction of $\Phi_0$ or its  Fourier coefficients when the core contribution $\Phi_1$ is clearly dominating (as is expected at lower spherical harmonic degrees in realistic geomagnetic field models) and when 
$\SS_{R_1}$ is close to $\SS_{R_0}$. The domination of the core contribution has been simulated to some extent in the presented examples but is expected to be stronger in real scenarios.
\\[3ex]
\begin{ackno}
 The work of CG was partly supported by DFG GE 2781/1-1.
\end{ackno}


\appendix
\section{Appendix: Balayage of Distributions}
\label{sec:appendix}
Since potentials of distributions do not seem to be widely treated in the literature, let us briefly justify the statements made in Section \ref{sec:harmpot}. For any distribution $D$ supported in a compact set  $\Omega\subset\R^3$, the corresponding potential $p_D$ has been formally defined in \eqref{eqn:p_Ddef} via 
\begin{align}\label{eqn:ap_Ddef}
 p_D(x)=D\left(-\frac{1}{4\pi}\frac{1}{|x-\cdot|}\right), \quad x\in\R^3\setminus \Omega.
\end{align}
Strictly speaking, this definition is not valid in that $-1/(4\pi|x-\cdot|)$ is neither
smooth nor compactly supported in $\RR^3$. However, for any compactly supported $\varphi_x\in C^{\infty}(\R^3)$ with $\varphi_x\equiv1$ in a neighborhood of $\Omega$ and $\varphi_x\equiv 0$ in a neighborhood of $x$, the function $g_{\varphi_x}(y)=-\frac{1}{4\pi}\frac{1}{|x-y|}\varphi_x(y)$ is in $C^{\infty}(\R^3)$ and compactly supported. 
Clearly  $D(g_{\varphi_x})$ is independent of the choice of $\varphi_x$,
for if $\psi_x$ is another function with the same properties then
$g_{\psi_x}-g_{\varphi_x}$ is supported in $\RR^3\setminus \Omega$ so that
$D(g_{\psi_x}-g_{\varphi_x})=0$.
Therefore \eqref{eqn:ap_Ddef} makes good sense  if we understand the latter to mean $p_D(x)=D(g_{\varphi_x})$.

In what follows, we restrict ourselves to the case where $\Omega$ has smooth 
boundary $\partial \Omega$. This is no loss of generality for the matter discussed 
in the paper, because we only consider
situations where $\Omega$ is a closed ball and we want to define balayage onto 
the boundary sphere.
The lemma below is a simple consequence of known density results for the fundamental solution of the Laplacian in $L^2(\partial \Omega)$, $C^{0}(\partial \Omega)$, and $W^{k,2}(\partial \Omega)$,  see, e.g., \cite{freeden80,freedenmichel04a,grothaus10}.

\begin{prop}\label{lem:density}
 Let $\Omega \subset \R^3$ be a compact, simply connected set with $C^{\infty}$-boundary $\partial \Omega$ and let $g_x(y)=\frac{1}{|x-y|}$. Then, the set of functions $\textnormal{span}\{g_x:x\in\R^3\setminus \Omega\}$ is dense in $C^{k}(\partial \Omega)$, for any $k\in \mathbb{N}_0$.
\end{prop}

\begin{prf}
For every $f\in C^{k}(\partial \Omega)$ and $\eps>0$, there exists $\bar{f}\in C^{\infty}(\partial \Omega)$ with $\|f-\bar{f}\|_{C^{k}(\partial \Omega)}<\eps$. In particular, $\bar{f}$ is an element of the Sobolev space $W^{k+2,2}(\partial \Omega)$. 
By \cite[Thm. 8.8]{grothaus10} we can find $N>0$, coefficients $a_i\in\mathbb{R}$, and points $x_i\in\R^3\setminus \Omega$, $i=1,\ldots,N$, such that
\begin{align*}
 \left\|\bar{f}-\sum_{i=1}^N a_i \frac{1}{|x_i-\cdot|}\right\|_{W^{k+2,2}(\partial \Omega)}<\eps.
\end{align*}
The Sobolev embedding theorem (see, e.g., \cite{adams03}) now yields that  
$W^{k+2,2}(\partial \Omega)\subset  C^{k}(\partial \Omega)$ and 
\begin{align*}
 \left\|\bar{f}-\sum_{i=1}^Na_i \frac{1}{|x_i-\cdot|}\right\|_{C^ {(k)}(\partial \Omega)}\leq M\left\|\bar{f}-\sum_{i=1}^N a_i\frac{1}{|x_i-\cdot|}\right\|_{W^{k+2,2}(\partial \Omega)}<M\eps
\end{align*}
for some constant $M>0$ depending only on $k$, which finishes the proof.
\end{prf}

\begin{prop}
 Let $ \Omega\subset \R^3$ be a compact, simply connected set with $C^{\infty}$-boundary $\partial \Omega$, and let $D$ be a distribution with support in $\Omega$. Then, 
there exists a unique distribution $\hat D$ on $\partial \Omega$ such that 
 \begin{align*}
  p_D(x)=p_{\hat D}(x), \quad x\in\R^3\setminus \Omega.
 \end{align*}
 We call $\hat D$ the balayage of $D$ onto $\partial \Omega$.
\end{prop}

\begin{prf}
First, we deal with the existence of a balayage. Since $D$ is compactly supported, it is known that there are finitely many compactly supported
continuous functions 
$\Phi_j$ and multiindices $\alpha_j\in\mathbb{N}_0^3$, $j=1,\ldots,m$,  such that $D=\sum_{j=1}^m\partial_{\alpha_j}\Phi_j$ (see, e.g., \cite{rudin91}). Due to this representation, $D$ acts on compactly supported functions $g\in C^{M}(\R^3)$, with $M=\max_{i=1,\ldots,N}|\alpha_i|$. Let $f$ be a function in $C^{\infty}(\partial \Omega)$ and $h$ its unique harmonic continuation to the interior of $\Omega$ with $h=f$ on $\partial \Omega$ \cite[Ch. 2]{lions68}. A compactly supported function $g_f\in C^{M}(\R^3)$ satisfying  $g_f=h$ in $\Omega$ can be computed as follows. The smoothness of $\partial \Omega$ implies there is an open cover $\{{U}_i\}_{i\in\mathbb{N}}$ of $\partial \Omega$ by open sets in $\RR^3$ and diffeomorphisms $\Psi_i\in C^{\infty}({U}_i,\mathbb{B}_1)$ that satisfy $\Psi_i({U}_i\cap \partial \Omega)\subset\R^2\times\{0\}$, $\Psi_i({U}_i\cap \Omega)\subset\R^3_-$, and $\Psi_i({U}_i\cap (\mathbb{R}^3\setminus \Omega))\subset\R^3_+$. Here,
$\R^3_\pm$ refer to upper and lower half spaces.
Let  $\{\varphi_i\}_{i\in\mathbb{N}}\subset C^{\infty}(\R^3)$ be a partition of unity 
 subordinated to the cover $\{{U}_i\}_{i\in\mathbb{N}}$. According to the construction in \cite[(2.21)]{lions68}, there exist functions $\bar{g}_i\in C^{M}(\R^3)$, $1\leq i\leq N$,
compactly supported in $\BB_1$,
with $\bar{g}_i=(h\varphi_i)\circ\Psi_i^{-1}$ on $\R^3_-\cap\BB_1$ for every
 $i$. The function 
$g_f=\sum_{i=1}^\infty\bar{g}_i\circ \Psi_i$ gives us the desired extension of $h$.
We now define $\hat D$ for any $f\in C^{\infty}(\partial \Omega)$ by
\begin{align}\label{eqn:hD}
 \hat D(f)=D(g_f).
\end{align}
Since any two $C^M$-smooth extensions of $h$ have the same derivatives of 
order less than or equal to $M$  on
$\Omega$, we see that $\hat D$ does not depend on the particular extension of $h$
that we use. 
Thus, it holds that
\begin{align*}
  p_D(x)=p_{\hat D}(x), \quad x\in\R^3\setminus \Omega,
\end{align*}
because when $x\notin \Omega$, then $g_x(y)$ is a harmonic function of $y$ in a neighborhood
of $\Omega$.
Uniqueness of $\hat D$ is a direct consequence of the requirement $\hat{D}(g_x)=p_{\hat D}(x)=p_D(x)=D(g_x)$ for $x\in\R^3\setminus \Omega$, and  
of Lemma \ref{lem:density} which guarantees the density of 
$\{g_x:x\in\R^3\setminus \Omega\}$ in $C^{k}(\partial \Omega)$ for 
all $k\in\mathbb{N}_0$.
\end{prf}

\section{Appendix: Differential forms and Hodge theory}
\label{sec:appendix2}
Below we gather some basic definitions and facts from Hodge theory on
a smooth simply connected surface $\mathcal{M}$ embedded in $\RR^3$,
that will be used to prove the rotation lemma in Appendix \ref{sec:appendix1}.
A detailed and more general treatment can be found,
{e.g.,} in \cite[Ch. 6]{Warner}.

Tangent spaces, smooth functions, vector fields,
metric tensor, area measure and Lebesgue spaces
are defined as in Section \ref{sec:aux}. 
Note that
$\mathcal{M}$ must be a finite 
union of topological spheres, as
follows from the classification theorem for surfaces \cite{Massey}
and the fact that $g$-holed tori are not simply connected while
projective planes cannot embed in $\RR^3$.
In particular $\mathcal{M}$ is orientable.

For $\mathcal{V}$ a real vector space of dimension 2,  let
$\mathcal{V}^*$ indicate its dual and $\mathcal{A}_2\mathcal{V}$ 
the bilinear alternating forms on $\mathcal{V}$. If $(v_1,v_2)$ 
is a basis of $\mathcal{V}$, the linear maps $v^*_1,v^*_2:\mathcal{V}\to\RR$
such that $v^*_j(v_k)=\delta_{jk}$ form a basis of $\mathcal{V}^*$, 
dual to $(v_1,v_2)$. 
The bilinear alternating form $v^*_1\wedge v_2^*$ defined by
\[v_1^*\wedge v_2^*(w_1,w_2)= \textnormal{det}\big((v_j^*(w_k))_{j,k=1,2}\big)=v_1^*(w_1)v_2^*(w_2)-v_1^*(w_2)v_2^*(w_1)\]
is a basis of the 1-dimensional space $\mathcal{A}_2\mathcal{V}$. Hereafter
we put
\[
\mathcal{E}\mathcal{V}=\RR\oplus\mathcal{V}^*\oplus\mathcal{A}_2\mathcal{V}.
\]
If $(w_1,w_2)$ is another basis of $\mathcal{V}$,
we say that $(w_1,w_2)$ has the same 
orientation as $(v_1,v_2)$ if $v_1^*\wedge v_2^*(w_1,w_2)>0$, the
opposite orientation if $v_1^*\wedge v_2^*(w_1,w_2)<0$.
We orient $\mathcal{V}$ by choosing one 
of the two equivalence classes of bases with the same orientation.
If $\mathcal{V}$ is equipped with a Euclidean scalar product
$\langle  \cdot,\cdot\rangle$, then
each $L\in\mathcal{V}^*$  is of the form $L(v)=\langle w,v\rangle$ for 
some unique
$w\in\mathcal{V}$. This way we identify $\mathcal{V}^*$ with 
$\mathcal{V}$ and $\mathcal{A}_2\mathcal{V}$ with the exterior product
$\mathcal{V}\wedge\mathcal{V}$ (the tensor product $\mathcal{V}\otimes\mathcal{V}$ quotiented by all  relations $v\otimes v=0$). Under this identification,
given a positively 
oriented orthonormal 
basis $(e_1,e_2)$ of $\mathcal{V}$, we define the star operator 
$\mathcal{E}\mathcal{V}\to\mathcal{E}\mathcal{V}$ to be the linear map such 
that $*1=e_1\wedge e_2$, $*(e_1)=e_2$, $*(e_2)=-e_1$, $* (e_1\wedge e_2)=1$.
The star operator does not depend on the positively
oriented orthonormal basis we use to define it. Clearly, $**=\textnormal{id}$ 
on $\RR\oplus\mathcal{A}_2\mathcal{V}$ and $**=-\textnormal{id}$ on 
$\mathcal{V}^*$.

We now introduce differential forms on $\mathcal{M}$.
A 0-form is a function $\mathcal{M}\to\RR$,
a 1-form is
a map associating to each $x\in\mathcal{M}$ a member 
of $T_x^*$,
a 2-form is a map associating to $x$ a member of $\mathcal{A}_2 T_x$; 
here and below,
$T_x$ indicates the tangent space to $\mathcal{M}$ at $x$.
Given a $k$-form $\omega$ and a chart $(U,\psi)$ on 
$\mathcal{M}$ with $\psi(U)=V\subset\RR^2$, one can define a $k$-form 
$\tilde{\omega}$ on $V$ by the rule
\begin{equation}
\label{locform}
\tilde{\omega}[y](v_1,\cdots,v_k)=\omega[\psi^{-1}(y)](\DD\psi^{-1}(v_1),\cdots,
\DD\psi^{-1}(v_k)),\quad y\in V,\quad v_1,\cdots, v_k\in\RR^2,
\end{equation}
which represents $\omega$ in local coordinates using the isomorphism
$\DD\psi^{-1}(y):\RR^2\to T_{\psi^{-1}(y)}$.
This way a form on $\mathcal{M}$ may be regarded as a collection of forms on 
images of
charts which define the same form $\omega$  on overlaps {via}
\eqref{locform}. Hence if we use a superscript prime to denote another 
system of local coordinates and if we set $h=\psi'\circ\psi^{-1}$ for
the corresponding change of charts, we have if $k=2$ that
\begin{equation}
\label{chvar2}
\tilde{\omega}[y](v_1,v_2)=
(\textnormal{det} (\DD h(y)))
\,\tilde{\omega}'[h(y)](v_1,v_2),
\quad y\in V\cap h^{-1}(V').
\end{equation}
A 1-form $\omega$ can be written in local coordinates as
 $\tilde{\omega}[y]=a(y)\dd y_1+b(y)\dd y_2$, where
$a$, $b$ are real functions of $y\in \psi(U)$ and $\dd y_1$, $\dd y_2$ is the 
basis of $(\RR^2)^*$ dual to the
canonical basis of $\RR^2$. If $\omega$ is a 2-form, then 
$\tilde{\omega}[y]=c(y)\dd y_1\wedge \dd y_2$ where $c$ is real-valued
on $V$. 
The wedge product is an  associative binary operation on forms,
bilinear over functions, that associates to  a $k_1$-form $\omega_1$ and a $k_2$-form 
$\omega_2$ a
$k_1+k_2$-form $\omega_1\wedge\omega_2$ such that,
in local coordinates, $*(\dd y_1)\wedge *(\dd y_2)=\dd y_1\wedge \dd y_2=-
*(\dd y_2)\wedge *(\dd y_1)$ and $*(\dd y_1)\wedge*(\dd y_1)=*(\dd y_2)\wedge*(\dd y_2)=0$.
Note that $k$-forms with $k>2$ (mapping $x\in\mathcal{M}$ 
to a $k$-linear 
alternating map on $(T_x)^k$) are 
identically zero for $T_x$ has dimension 2.
The wedge product is independent of the chart used to 
compute a local representative.
We say that a 1-form or a 2-form is smooth if its coefficients $a,b$ or $c$ are smooth functions in every chart.
We write $\Lambda^k\mathcal{M}$ for the space of smooth forms of degree $k$
on $\mathcal{M}$, and we let 
$\Lambda\mathcal{M}=\oplus_{k=0}^2\Lambda^k\mathcal{M}$ for the direct sum.

A smooth 2-form $\omega$ can be integrated over a Borel set 
$E\subset\mathcal{M}$: if $(U,\psi)$ is a chart with $\psi(U)=V$ and 
$\tilde{\omega}[y]=c(y)\dd y_1\wedge \dd y_2$, and if moreover $E\subset U$, we set
$\int_E\omega=\int_{\psi(E)}c(y)\dd\lambda(y)$ where $\lambda$ indicates Lebesgue 
measure. In the general case we cover $E$ with finitely
many domains of charts and we use 
a partition of unity; relation \eqref{chvar2} and the change of variable 
formula ensure that the definition does not depend on which charts
or partition we use. 

The exterior differential $\dd:\Lambda^k\mathcal{M}\to\Lambda^{k+1}\mathcal{M}$ is defined as follows.
If $g$ is a function, then $\dd g$ is the
usual differential, namely in local coordinates
$\widetilde{\dd g}=\partial \tilde{g}/\partial y_1\dd y_1+\partial \tilde{g}/\partial y_2 \dd y_2$. If $\tilde{\omega}=a\dd y_1+b\dd y_2$ is a 1-form in local coordinates, then
$\widetilde{\dd \omega}=
(\partial b/\partial y_1-\partial a/\partial y_2)\dd y_1\wedge \dd y_2$. The 
differential of a 2-form is zero. Differentiation is 
meaningful in that it is independent
of the chart used to compute its local representative. Moreover it holds 
that $\dd\circ \dd=0$. If $\dd \omega=0$, we say that $\omega$ is closed, and
if $\omega=\dd\nu$ for some $\nu$ we say that $\omega$ is exact.
Exact forms are closed, and the quotient space of closed $k$-forms by 
exact $k$-forms is called the $k$-th (de Rham) cohomology group 
$H^k(\mathcal{M})$. The simple connectedness of $\mathcal{M}$ means that
$H^1(\mathcal{M})=0$, i.e. every closed 1-form on $\mathcal{M}$ is exact 
\cite[Ch. 5]{Warner}. 

The Hodge-star operator maps $\Lambda^k\mathcal{M}$ to  $\Lambda^{2-k}\mathcal{M}$  for $0\leq k\leq 2$, 
by acting pointwise as the star operator on $\mathcal{E} T_x$ for each
$x\in\mathcal{M}$. If we identify
a 1-form $\omega$  with the tangent 
vector field $\mathbf{v}_\omega$ such that
$\omega[x](w)=\mathbf{v}_\omega(x)\cdot w$
for $w\in T_x$, then the Hodge star operator merely rotates 
$\mathbf{v}_\omega$ by $\pi/2$ in the tangent space at each point.
To check that it maps  smooth forms to smooth forms, 
we need only produce in a neighborhood of each  $x_0\in\mathcal{M}$
a positively oriented orthonormal basis $(e_1(x),e_2(x))$
of $T_x$ that varies smoothly with $x$. 
If $(U,\psi)$ is a chart with $x_0\in U$ and
$V=\psi(U)$,  we may   choose $e_j(\psi^{-1}(y))=D\psi^{-1}(y)
\mathbf{G}(y)^{-1/2}\kappa_j$ for $y\in V$, 
where $\mathbf{G}$ is the metric tensor and
$\kappa_1,\kappa_2$ the canonical basis of $\RR^2$.
We denote the action of
the Hodge star operator on a form $\omega$ by $*\omega$, as
no confusion should arise with the star operator acting on 
$\mathcal{E} T_x$ for fixed $x$. 
Next, one defines a pairing on 
$\Lambda^k\mathcal{M}$ by letting
\begin{equation}
\label{scf}
\langle\omega_1,\omega_2\rangle=\int_{\mathcal{M}}\omega_1\wedge*\omega_2.
\end{equation}
Identifying $T_x^*$ and $T_x$ {via} the scalar product in $\RR^3$,
it follows  
from the definitions, with the notation of \eqref{chvar2}, 
that in local coordinates $\widetilde{e_1\wedge *e_2}=\widetilde{e_2\wedge *e_1}=0$ and, in addition,
\[
\widetilde{1\wedge*1}=\widetilde{e_1\wedge *e_1}=\widetilde{e_2\wedge*e_2}=
\widetilde{(e_1\wedge e_2)\wedge*
(e_1\wedge e_2)}=\sqrt{g}\,\dd y_1\wedge \dd y_2.
\]
Hence \eqref{scf} is symmetric and positive definite,
moreover we have that 
\begin{equation}
\label{L2vecform}
\langle f,f\rangle=\|f\|^2_{L^2(\mathcal{M})}\quad\textnormal{and}\quad
\langle \omega,\omega\rangle=\|\mathbf{v}_\omega\|^2_{L^2(\mathcal{M},\RR^3)},
\qquad  f\in\Lambda^0\mathcal{M},\ \omega\in\Lambda^1\mathcal{M}. 
\end{equation}
One extends $\langle\cdot,\cdot\rangle$ to a scalar product on $\Lambda\mathcal{M}$
by requiring that forms of different degree are orthogonal.
Let $\delta:\Lambda^k\mathcal{M}\to\Lambda^{k-1}\mathcal{M}$ be
the operator defined by $\delta(\omega)=(-1)^{k(2-k)}*\dd(*\omega)$. 
Since $**\omega=(-1)^k\omega$ when $\omega\in\Lambda^k\mathcal{M}$,
it holds if $\omega_1\in\Lambda^{k-1}\mathcal{M}$ and 
 $\omega_2\in\Lambda^{k}\mathcal{M}$ that
\[\dd(\omega_1\wedge*\omega_2)=\dd \omega_1\wedge*\omega_2+(-1)^{k-1}\omega_1\wedge
\dd(*\omega_2)=
\dd \omega_1\wedge*\omega_2-\omega_1\wedge*\delta(\omega_2),\]
and since the left hand side integrates to $0$ over $\mathcal{M}$ by Stoke's 
theorem it implies that $\delta$ is the adjoint of $\dd$  in $\Lambda\mathcal{M}$
equipped with \eqref{scf}. In particular, we see from 
\eqref{L2vecform} that $\delta$ must coincide with the divergence operator 
on $\Lambda^1\mathcal{M}$ when the latter is identified with smooth 
tangent vector fields. The operator
$\Delta=\dd\delta+\delta \dd$ which maps $\Lambda^k\mathcal{M}$ into itself is 
the Laplace Beltrami operator on $\Lambda\mathcal{M}$. 
The kernel of $\Delta$ in $\Lambda^k\mathcal{M}$
is the space of harmonic $k$-forms, denoted by $\mathcal{H}^k$.
Now, a fundamental result in Hodge theory \cite[Thm. 6.8]{Warner} is 
the existence of an orthogonal sum:
\begin{equation}
\label{fundHodge}
\Lambda^k\mathcal{M}=\dd(\Lambda^{k-1}\mathcal{M})\oplus\delta (\Lambda^{k+1}\mathcal{M})\oplus\mathcal{H}^k, \qquad k=0,1,2,
\end{equation}
where orthogonality holds with respect to \eqref{scf} (by convention $\Lambda^{-1}\mathcal{M}=\{0\}$).
Using \eqref{fundHodge} and elliptic regularity theory, one 
can further show that each equivalence class in the cohomology group
$H^k(\mathcal{M})$ has a unique harmonic representative 
\cite[Thm. 6.11]{Warner}.
Since $H^1(\mathcal{M})=\{0\}$ we deduce that $\mathcal{H}^1=0$, hence
the orthogonal decomposition \eqref{fundHodge} specializes in our case to
\begin{equation}
\label{fundHodges}
\Lambda^1\mathcal{M}=\dd(\Lambda^{0}\mathcal{M})\oplus\delta 
(\Lambda^{2}\mathcal{M}). 
\end{equation}
Moreover, since $*$ is obviously surjective 
$\Lambda^2\mathcal{M}\to\Lambda^0\mathcal{M}$ (for the inverse image of a 
smooth function $f$ is $f\dd e_1\wedge \dd e_2$), we get that
\begin{equation}
\label{divref}
\textnormal{Im}(\delta:\Lambda^2\mathcal{M}\to\Lambda^1\mathcal{M})
=\textnormal{Im}(*\dd:\Lambda^0\to\Lambda^1\mathcal{M}).
\end{equation}

\section{Appendix: the Rotation Lemma}
\label{sec:appendix1}
In the notation of Section \ref{sec:aux},
we prove below  that the operator 
$J:\mathcal{T}_R\to \mathcal{T}_R$, which rotates a tangent
vector field by $\pi/2$ at every point
in the positively oriented tangent plane, isometrically maps tangential 
gradients to divergence free vector fields 
and vice-versa. This we call the rotation lemma.
The result actually holds on any smooth simply 
connected  compact surface $\mathcal{M}$ embedded in $\RR^3$, and we
deal below with this more general version  but 
restricting ourselves to the sphere would not 
simplify the proof.

Gradients, Sobolev spaces, tangent and divergence-free vector fields 
are defined as in Section \ref{sec:aux}. 
Thus, letting $\mathcal{T}$, $\mathcal{G}$ and $\mathcal{D}$ indicate
respectively tangent, gradient, and divergence 
free vector fields in  $L^2(\mathcal{M},\RR^3)$, we have  the orthogonal 
decomposition:
\begin{equation}
\label{HHMa}
\mathcal{T}=\mathcal{G}\oplus\mathcal{D}.
\end{equation} 
As pointed out in Appendix \ref{sec:appendix2},
$\mathcal{M}$ is orientable, which makes it possible to define $J$
as rotation of  a tangent vector field pointwise by $\pi/2$ in 
the positively oriented tangent plane.

\begin{lem}
\label{RGD}
For $\mathcal{M}$ a compact simply connected surface embedded in $\RR^3$,
the map $J:\T\to\T$ isometrically maps
$\G $ onto $\D$ and conversely.
\end{lem}

\begin{prf}
That $J$ is isometric is obvious for 
it preserves length pointwise.
Moreover, since $J^2=-I$, it suffices  to establish  that $J(\G)=\D$.
By \eqref{HHMa} this amounts to prove that $\T=\G\oplus J(\G)$, and 
since smooth vector fields and smooth functions are dense in $\T$ and
$W^{1,2}(\mathcal{M})$ respectively, it is enough by the  isometric character 
of $J$ to show that
\begin{equation}
\label{smoothH}
\T_S=\G_S\oplus J(\G_S),
\end{equation}
where the subscript ''$S$'' indicates the smooth elements of 
the corresponding space. Now, representing a 1-form $\omega$ as
the pointwise Euclidean scalar product with a tangent vector field 
$\mathbf{v}_\omega$ as we did in Appendix \ref{sec:appendix2},
we have for any smooth function $f:\mathcal{M}\to\RR$
that $\mathbf{v}_{\dd f}$ is just the gradient
$\nabla_{\mathcal{M}}f$ and, since we observed in the latter appendix 
that the Hodge star operator coincides with $J$ on $\mathbf{v}_\omega$,
the decomposition \eqref{smoothH} follows immediately
from
\eqref{fundHodges} and \eqref{divref}.
 \end{prf}

\end{document}